\documentclass[12pt]{amsart}
\usepackage{amssymb}
\usepackage{amsmath}
\usepackage{amscd}
\usepackage[all]{xy}
\theoremstyle{plain}
\newtheorem{theorem}{Theorem}[section]
\newtheorem{corollary}[theorem]{Corollary}
\newtheorem{proposition}[theorem]{Proposition}
\newtheorem{lemma}[theorem]{Lemma}
{\theoremstyle{remark}

\newtheorem{remark}[theorem]{Remark}}
{\theoremstyle{definition}
\newtheorem{definition}[theorem]{Definition}
\newtheorem{example}[theorem]{Example}
\newtheorem{problem}[theorem]{Problem}}

\newcommand{\benu}{\begin{enumerate}\renewcommand{\labelenumi}{{\rm (\roman{enumi})}}\renewcommand{\itemsep}{0pt}}
\newcommand{\eenu}{\end{enumerate}}

\newcommand{\N}{\mathbb{N}}
\newcommand{\Z}{\mathbb{Z}}
\newcommand{\R}{\mathbb{R}}
\newcommand{\C}{\mathbb{C}}
\newcommand{\T}{\mathbb{T}}
\newcommand{\e}{\varepsilon}
\newcommand{\cK}{{\mathcal K}}
\newcommand{\cL}{{\mathcal L}}
\newcommand{\cM}{{\mathcal M}}
\newcommand{\dcM}{{\mathcal M}'}
\newcommand{\tcM}{\widetilde{\mathcal M}}

\newcommand{\M}{\mathbb{M}}

\newcommand{\F}{{\mathcal F}}
\newcommand{\G}{{\mathcal G}}
\newcommand{\cO}{{\mathcal O}}

\newcommand{\ip}[2]{\langle{#1},{#2}\rangle}

\newcommand{\s}[3]{{{#1}^{#2}_{\textnormal{#3}}}}
\newcommand{\rs}[1]{{\textnormal{#1}}}
\newcommand{\per}{{\textnormal{per}}}
\newcommand{\aper}{{\textnormal{aper}}}
\newcommand{\drtimes}{{}_d\!\times_r}
\DeclareMathOperator{\id}{id}

\DeclareMathOperator{\Ad}{Ad}

\DeclareMathOperator{\Orb}{Orb}
\DeclareMathOperator{\Per}{Per}
\DeclareMathOperator{\Aper}{Aper}
\DeclareMathOperator{\supp}{supp}
\DeclareMathOperator{\spa}{span}

\begin{document}
\title[A class of $C^*$-algebras III, ideal structures]
{A class of \boldmath{$C^*$}-algebras generalizing 
both graph algebras and homeomorphism \boldmath{$C^*$}-algebras III, \\
ideal structures}
\author[Takeshi KATSURA]{Takeshi KATSURA}
\address{Department of Mathematics, 
Hokkaido University, Kita 10, Nishi 8, 
Kita-Ku, Sapporo, 060-0810, JAPAN}
\email{katsura@math.sci.hokudai.ac.jp}
\date{}

\begin{abstract}
We investigate the ideal structures of the $C^*$-algebras 
arising from topological graphs. 
We give the complete description of ideals of such $C^*$-algebras 
which are invariant under the so-called gauge action, 
and give the condition on topological graphs 
so that all ideals are invariant under the gauge action. 
We get conditions for our $C^*$-algebras 
to be simple, prime or primitive. 
We completely determine the prime ideals, 
and show that most of them are primitive. 
Finally, we construct a discrete graph 
such that the associated $C^*$-algebra is prime but not primitive. 
\end{abstract}

\maketitle

\setcounter{section}{-1}

\section{Introduction}

From a homeomorphism on some locally compact space, 
we can construct a $C^*$-algebra called a homeomorphism $C^*$-algebra 
(or a crossed product). 
The ideal structure of the homeomorphism $C^*$-algebra 
reflects the orbit structure of the given homeomorphism 
(see \cite{Wi,T2,T3} for example). 
On the other hand, 
the ideal structures of graph algebras have been examined 
by many researchers (for example, \cite{HR,KPRR,BPRS,BHRS,DT,HS}). 
These two lines of analysis have several similar aspects in common. 
Topological graphs introduced in \cite{Ka1} 
generalize dynamical systems and (discrete) graphs, 
and the construction of $C^*$-algebras from topological graphs 
defined in \cite{Ka1} 
generalizes the ones of homeomorphism $C^*$-algebras 
and graph algebras. 
In this paper, 
we unify the two analyses of ideal structures, 
and generalize them to the setting of topological graphs. 
The purposes of this paper 
include giving a dynamical insight into the theory of 
graphs and graph algebras. 
We mainly borrow terminologies 
from the theory of dynamical systems. 

In Section~\ref{SecPre}, 
we recall the definition of topological graphs, 
and the way to construct a $C^*$-algebra $\cO(E)$ 
from a topological graph $E$. 
In Section~\ref{adpair}, 
we introduce the notion of 
invariant sets and admissible pairs of a topological graph $E$, 
and see that admissible pairs arise 
from ideals of the $C^*$-algebra $\cO(E)$. 
Conversely, in Section~\ref{SecGII}, 
we see that an ideal of the $C^*$-algebra $\cO(E)$ 
arises from each admissible pair of $E$. 
We show that by these correspondences, 
the set of all gauge-invariant ideals of $\cO(E)$ 
corresponds bijectively to the set of all admissible pairs of $E$ 
(Theorem~\ref{gaugeideals}). 
The key ingredient for the proof of this theorem 
is the observation done in \cite[Section~3]{Ka2}. 
In Section~\ref{SecOrbit}, 
we introduce the notion of orbits 
which generalizes the one for dynamical systems. 
In Section~\ref{SecHS}, 
we study hereditary and saturated sets 
which are generalizations of the ones in 
the theory of graph algebras. 
As an application, 
we give a characterization of maximal heads, 
which is defined in the previous section. 
In Section~\ref{SecTopFree}, 
we see that some of the ideals of $\cO(E)$ 
are strongly Morita equivalent to the $C^*$-algebras 
of subgraphs of $E$. 
Using this observation, 
we show that the topological freeness is needed 
in the Cuntz-Krieger Uniqueness Theorem (Theorem~\ref{topfree}). 
In Section~\ref{SecFree}, 
we define periodic and aperiodic points of 
a topological graph, and freeness of topological graphs. 
We show that a topological graph $E$ is free 
if and only if all ideals of $\cO(E)$ are gauge-invariant 
(Theorem~\ref{GII&free}). 
In Sections~\ref{SecMin} and~\ref{SecPrime}, 
we generalize the notions of minimality 
and topological transitivity 
from dynamical systems to topological graphs, 
and give a couple of equivalent conditions 
on a topological graph $E$ so that 
the $C^*$-algebra $\cO(E)$ is simple (Theorem~\ref{ThmSimple}) 
and prime (Theorem~\ref{C*prime}), respectively. 
For the primeness of the $C^*$-algebra $\cO(E)$, 
we use the results in Section~\ref{SecPrimePair} 
where we introduce the primeness for admissible pairs, 
and completely determine prime admissible pairs. 
Using the analysis in Section~\ref{SecPrimePair}, 
we completely determine prime ideals of 
our $C^*$-algebras (Theorem~\ref{ThmPrime}). 
In Section~\ref{SecIrrep}, 
we give some ways to construct irreducible representations 
of our $C^*$-algebras, 
and show that most of their prime ideals are primitive 
(Theorem~\ref{ThmPrimitive}). 
In Section~\ref{SecPrimitive}, 
we give one sufficient condition and one necessary condition 
on a topological graph $E$ 
so that 
the $C^*$-algebra $\cO(E)$ is primitive. 
Finally, we construct a discrete graph $E$ 
such that $\cO(E)$ is prime but not primitive 
(a similar construction can be found in \cite{Ka4}). 
This $C^*$-algebra is an inductive limit of finite dimensional 
$C^*$-algebras and simultaneously residually finite dimensional. 
Thus $\cO(E)$ is an easier example 
of a prime $C^*$-algebra which is not primitive 
than the first such example found in \cite{We}. 

While this paper was under construction, 
P. S. Muhly and M. Tomforde introduced topological quivers in \cite{MT} 
which include topological graphs as special examples. 
In the same paper, 
they give the method to construct $C^*$-algebras from them 
which generalizes our construction. 
Among others, 
they analyze the ideal structures of their $C^*$-algebras, 
and so some of our results are generalized 
to their general setting. 
In particular, 
Theorems~\ref{gaugeideals},~\ref{topfree} and~\ref{ThmSimple} 
in this paper 
are valid for topological quivers 
without significant changes of statements and proofs. 

\medskip
The author would like to thank Mark Tomforde and Paul S. Muhly 
for useful discussion on Proposition~\ref{Baire} 
and on their topological quivers. 
He is also grateful to Jun Tomiyama for useful discussion 
about dynamical systems. 
He thanks the referee for careful reading. 
This work was partially supported by Research Fellowship 
for Young Scientists of the Japan Society for the Promotion of Science.

\section{Preliminaries}\label{SecPre}

\begin{definition}
A {\em topological graph} $E=(E^0,E^1,d,r)$ consists of 
two locally compact spaces $E^0$ and $E^1$, 
and two maps $d,r\colon E^1\to E^0$, 
where $d$ is locally homeomorphic 
and $r$ is continuous.
\end{definition}

In this paper, $E=(E^0,E^1,d,r)$ 
always means a topological graph. 
The triple $(E^1,d,r)$ is called a topological correspondence over $E^0$ 
in \cite{Ka1}, which can be considered as a multi-valued continuous map. 
A (topological) dynamical system is a pair 
$\Sigma=(X,\sigma)$ consisting of 
a locally compact space $X$ and a homeomorphism $\sigma$ on $X$. 
From a dynamical system $\Sigma=(X,\sigma)$, 
we get a topological graph $E_\Sigma$ by $E_\Sigma=(X,X,\id_X,\sigma)$. 
We regard topological graphs as a generalization of 
dynamical systems and borrow many notions 
from the theory of dynamical systems. 
Sometimes, we think $E^0$ as a set of vertices 
and $E^1$ as a set of edges, 
and that an edge $e\in E^1$ is directed 
from its domain $d(e)\in E^0$ 
to its range $r(e)\in E^0$. 
This viewpoint explains why $E=(E^0,E^1,d,r)$ is called 
a topological graph. 

Let us denote by $C_d(E^1)$ 
the set of continuous functions $\xi$ on $E^1$ such that 
$\ip{\xi}{\xi}(v)=\sum_{e\in d^{-1}(v)}|\xi(e)|^2<\infty$ for any $v\in E^0$ 
and $\ip{\xi}{\xi}\in C_0(E^0)$. 
For $\xi,\eta\in C_d(E^1)$ and $f\in C_0(E^0)$, 
we define $\xi f\in C_d(E^1)$ and $\ip{\xi}{\eta}\in C_0(E^0)$ by 
$$(\xi f)(e)=\xi (e)f(d(e)) \mbox{ for } e\in E^1$$
$$\ip{\xi}{\eta}(v)=\sum_{e\in d^{-1}(v)}\overline{\xi(e)}\eta(e)
\mbox{ for } v\in E^0.$$
With these operations,  
$C_d(E^1)$ is a (right) Hilbert $C_0(E^0)$-module 
(\cite[Proposition 1.10]{Ka1}).
We define a left action $\pi_r$ of $C_0(E^0)$ on $C_d(E^1)$
by $(\pi_r(f)\xi)(e)=f(r(e))\xi(e)$ for $e\in E^1$, 
$\xi\in C_d(E^1)$ and $f\in C_0(E^0)$.
Thus we get a $C^*$-correspondence $C_d(E^1)$ over $C_0(E^0)$. 

We set $d^0=r^0=\id_{E^0}$ and $d^1=d, r^1=r$. 
For $n=2$, $3$, \ldots, 
we recursively define a space $E^n$ of paths with length $n$ 
and domain and range maps $d^n,r^n\colon E^n\to E^0$ by 
$$E^n=\{(e',e)\in E^1\times E^{n-1}\mid d^1(e')=r^{n-1}(e)\},$$ 
$d^n((e',e))=d^{n-1}(e)$ and $r^n((e',e))=r^1(e')$. 
We can define a $C^*$-correspondence $C_{d^n}(E^n)$ over $C_0(E^0)$
similarly as $C_d(E^1)$. 
We have 
$C_{d^{n+m}}(E^{n+m})\cong C_{d^n}(E^n)\otimes C_{d^{m}}(E^{m})$ 
as $C^*$-correspondences over $C_0(E^0)$ for any $n,m\in\N=\{0,1,2,\ldots\}$.
As long as no confusion arises, 
we omit the superscript $n$ 
and simply write $d,r$ for $d^n,r^n$. 
Thus we get two maps $d,r\colon E^*\to E^0$ 
where $E^*=\coprod_{n=0}^\infty E^n$ 
is the {\em finite path space} of the topological graph $E$. 

\begin{definition}\label{DefTpl}
A {\em Toeplitz $E$-pair} on a $C^*$-algebra $A$ 
is a pair of maps $T=(T^0,T^1)$ 
where $T^0\colon C_0(E^0)\to A$ is a $*$-homomorphism 
and $T^1\colon C_d(E^1)\to A$ is a linear map satisfying that
\benu
\item $T^1(\xi)^*T^1(\eta)=T^0(\ip{\xi}{\eta})$ for $\xi,\eta\in C_d(E^1)$, 
\item $T^0(f)T^1(\xi)=T^1(\pi_r(f)\xi)$ 
for $f\in C_0(E^0)$ and $\xi\in C_d(E^1)$.
\eenu
\end{definition}

For a Toeplitz $E$-pair $T=(T^0,T^1)$, 
the equation $T^1(\xi)T^0(f)=T^1(\xi f)$ 
holds automatically 
from the condition (i). 
We write $C^*(T)$ for denoting the $C^*$-algebra 
generated by the images of the maps $T^0$ and $T^1$. 
For $n\geq 2$, we can define a linear map $T^n\colon C_d(E^n)\to C^*(T)$ 
by $T^n(\xi)=T^1(\xi_1)T^1(\xi_2)\cdots T^1(\xi_n)$
for $\xi=\xi_1\otimes\xi_2\otimes\cdots\otimes\xi_n\in C_d(E^n)$.
The linear space 
$$\spa\{T^n(\xi)T^m(\eta)^*\mid 
\xi\in C_d(E^n),\ \eta\in C_d(E^m),\ n,m\in\N\}$$
is dense in $C^*(T)$ (see the remark after \cite[Lemma 2.4]{Ka1}). 

We say that a Toeplitz $E$-pair $T=(T^0,T^1)$ is {\em injective} 
if $T^0$ is injective.
If a Toeplitz $E$-pair $T=(T^0,T^1)$ is injective,
then $T^n$ are isometric for all $n\in\N$.
We say a Toeplitz $E$-pair $T$ {\em admits a gauge action} 
if there exists an automorphism $\beta'_z$ on $C^*(T)$ 
with $\beta'_z(T^0(f))=T^0(f)$ and $\beta'_z(T^1(\xi))=zT^1(\xi)$
for every $z\in\T$. 
If $T$ admits a gauge action $\beta'$, 
then we have 
$$\beta'_z(T^n(\xi)T^m(\eta)^*)=z^{n-m}T^n(\xi)T^m(\eta)^*$$
for $\xi\in C_d(E^n)$ and $\eta\in C_d(E^m)$.

\begin{definition}
We define three open subsets 
$\s{E}{0}{sce},\s{E}{0}{fin}$ and $\s{E}{0}{rg}$ of $E^0$ by
$\s{E}{0}{sce}=E^0\setminus\overline{r(E^1)}$,
\begin{align*}
\s{E}{0}{fin}=\{v\in E^0\mid\mbox{ there exists}&\mbox{ a neighborhood } 
V \mbox{ of } v\\
&\mbox{ such that }
r^{-1}(V)\subset E^1 \mbox{ is compact}\},
\end{align*}
and $\s{E}{0}{rg}=\s{E}{0}{fin}\setminus\overline{\s{E}{0}{sce}}$.
We define two closed subsets $\s{E}{0}{inf}$ and $\s{E}{0}{sg}$ of $E^0$ by 
$\s{E}{0}{inf}=E^0\setminus \s{E}{0}{fin}$ and 
$\s{E}{0}{sg}=E^0\setminus \s{E}{0}{rg}$. 
\end{definition}

A vertex in $\s{E}{0}{sce}$ is called a {\em source}. 
We have $\s{E}{0}{sg}=\s{E}{0}{inf}\cup\overline{\s{E}{0}{sce}}$. 
Vertices in $\s{E}{0}{rg}$ are said to be {\em regular}, 
and those in $\s{E}{0}{sg}$ are said to be {\em singular}. 
We have that $\pi_r^{-1}(\cK(C_d(E^1)))=C_0(\s{E}{0}{fin})$ 
and $\ker\pi_r=C_0(\s{E}{0}{sce})$ (\cite[Proposition 1.24]{Ka1}).
Hence the restriction of $\pi_r$ to $C_0(\s{E}{0}{rg})$ is 
an injection into $\cK(C_d(E^1))$. 
When a topological graph $E$ is defined from a dynamical system $\Sigma$, 
we get $\s{E}{0}{rg}=E^0$. 
The set of regular vertices $\s{E}{0}{rg}$ 
is the maximal open subset 
such that the restriction 
$r\colon r^{-1}(\s{E}{0}{rg})\to \s{E}{0}{rg}$ is surjective and proper. 
This fact implies the following property of regular vertices, 
which roughly says that 
$\s{E}{0}{rg}$ is the part in which 
the topological correspondence $(E^1,d,r)$ 
is ``reversible''. 

\begin{lemma}[{\cite[Lemma 1.21]{Ka1}}]\label{RgPr}
For $v\in\s{E}{0}{rg}$, 
the set $r^{-1}(v)\subset E^1$ is a non-empty compact set, 
and for an open subset $U$ of $E^1$ with $r^{-1}(v)\subset U$ 
there exists a neighborhood $V$ of $v$ such that $r^{-1}(V)\subset U$.
\end{lemma}

For a Toeplitz $E$-pair $T=(T^0,T^1)$, 
we define a $*$-homomorphism $\varPhi\colon \cK(C_d(E^1))\to C^*(T)$ 
by $\varPhi(\theta_{\xi,\eta})=T^1(\xi)T^1(\eta)^*$ 
for $\xi,\eta\in C_d(E^1)$. 

\begin{definition}
A Toeplitz $E$-pair $T=(T^0,T^1)$ is called 
a {\em Cuntz-Krieger $E$-pair} 
if $T^0(f)=\varPhi(\pi_r(f))$ for any $f \in C_0(\s{E}{0}{rg})$. 

The $C^*$-algebra $\cO(E)$ is generated 
by the universal Cuntz-Krieger $E$-pair $t=(t^0,t^1)$. 
\end{definition}

Since $t^0$ is injective (\cite[Proposition 3.7]{Ka1}),
the $*$-homomorphism $\varphi\colon\cK(C_d(E^1))\to \cO(E)$ is injective, 
and the maps $t^n\colon C_d(E^n)\to\cO(E)$ 
are isometric for all $n\in\N$. 
The universal Cuntz-Krieger $E$-pair $t=(t^0,t^1)$ 
admits a gauge action, 
which will be denoted by $\beta\colon \T\curvearrowright \cO(E)$.

\section{Admissible pairs}\label{adpair}

In this section, 
we introduce invariant sets and admissible pairs 
for a topological graph $E$, 
and see that admissible pairs correspond to ideals 
of the $C^*$-algebra $\cO(E)$. 

\begin{definition}
A subset $X^0$ of $E^0$ is said to be {\em positively invariant} 
if $d(e)\in X^0$ implies $r(e)\in X^0$ 
for each $e\in E^1$, 
and to be {\em negatively invariant} 
if for $v\in X^0\cap \s{E}{0}{rg}$, 
there exists $e\in E^1$ with $r(e)=v$ and $d(e)\in X^0$.
A subset $X^0$ of $E^0$ is said to be {\em invariant} 
if $X^0$ is both positively and negatively invariant. 
\end{definition}

These terminologies come from regarding $E=(E^0,E^1,d,r)$ 
as a generalization of dynamical systems. 

For a closed positively invariant subset $X^0$,
we define a closed subset $X^1$ of $E^1$ 
by $X^1=d^{-1}(X^0)$. 
Since $X^0$ is positively invariant, 
we have $r(X^1)\subset X^0$.
Hence $X=(X^0,X^1,d_X,r_X)$ 
is a topological graph, 
where $d_X$ and $r_X$ are the restrictions of $d$ and $r$ 
to $X^1$. 

\begin{proposition}\label{condforNI}
For a closed positively invariant subset $X^0$ of $E^0$, 
the following conditions are equivalent;
\benu
\item $X^0$ is negatively invariant, 
\item $\s{X}{0}{sce}\cap\s{E}{0}{rg}=\emptyset$, 
\item $\s{X}{0}{sg}\subset\s{E}{0}{sg}$.
\eenu
\end{proposition}

\begin{proof}
(i)$\Rightarrow$(ii): 
If $X^0$ is negatively invariant, 
then every $v\in X^0\cap\s{E}{0}{rg}$ satisfies 
that $r_X^{-1}(v)=r^{-1}(v)\cap X^1\neq\emptyset$. 
Hence $\s{X}{0}{sce}\cap\s{E}{0}{rg}=\emptyset$.

(ii)$\Rightarrow$(iii): 
From (ii) we have $\s{X}{0}{sce}\subset\s{E}{0}{sg}$. 
Thus $\overline{\s{X}{0}{sce}}\subset\s{E}{0}{sg}$.
Clearly $\s{X}{0}{inf}\subset\s{E}{0}{inf}\subset\s{E}{0}{sg}$. 
Therefore we have $\s{X}{0}{sg}\subset\s{E}{0}{sg}$.

(iii)$\Rightarrow$(i): 
Take $v\in X^0\cap\s{E}{0}{rg}$. 
From (iii), we have $v\in\s{X}{0}{rg}$. 
Hence there exists $e\in X^1=d^{-1}(X^0)$ with $r(e)=v$ 
by Lemma~\ref{RgPr}. 
Thus $X^0$ is negatively invariant. 
\end{proof}

By Proposition~\ref{condforNI}, 
we have an inclusion $\s{X}{0}{sg}\subset \s{E}{0}{sg}\cap X^0$ 
for a closed invariant set $X^0$. 
There are many examples in which this inclusion is proper 
(see Example~\ref{AbsVal} and Proposition~\ref{row-finite}). 

\begin{definition}
A pair $\rho=(X^0,Z)$ consisting of 
two closed subsets $X^0$ and $Z$ of $E^0$ 
is called an {\em admissible pair} if $X^0$ is invariant 
and $\s{X}{0}{sg}\subset Z\subset \s{E}{0}{sg}\cap X^0$.
\end{definition}

Admissible pairs naturally arise from ideals of $\cO(E)$. 
We denote by $\G^0$ and $\G^1$ the images of 
$t^0\colon C_0(E^0)\to\cO(E)$ and $\varphi\colon\cK(C_d(E^1))\to\cO(E)$, 
respectively. 

\begin{definition}\label{XIZI}
For an ideal $I$ of $\cO(E)$,
we define closed subsets $X_I^0,Z_I$ of $E^0$ by
\begin{align*}
X_I^0&=\{v\in E^0\mid f(v)=0\mbox{ for all }f\in C_0(E^0)
\mbox{ with }t^0(f)\in I\},\\
Z_I&=\{v\in E^0\mid 
f(v)=0\mbox{ for all }f\in C_0(E^0)\text{ with }
t^0(f)\in I+\G^1\}.
\end{align*}
\end{definition}

The closed sets $X_I^0,Z_I\subset E^0$ are determined by
$$t^0(C_0(E^0\setminus X_I^0))=I\cap \G^0,\quad 
t^0(C_0(E^0\setminus Z_I))=(I+\G^1)\cap \G^0.$$
We denote by $\rho_I$ the pair $(X_I^0,Z_I)$ of closed subsets of $E^0$.
We will show that the pair $\rho_I=(X_I^0,Z_I)$ is admissible. 

\begin{proposition}\label{hereditary}
For an ideal $I$ of $\cO(E)$, 
the closed set $X_I^0$ is positively invariant.
\end{proposition}

\begin{proof}
Take $e\in E^1$ with $d(e)\in X_I^0$. 
Take $f\in C_0(E^0)$ with $t^0(f)\in I$ arbitrarily,
and we will show that $f(r(e))=0$. 
There exists $\xi\in C_d(E^1)$ with $\xi(e)=1$ and $\xi(e')=0$ 
for all $e'\in d^{-1}(d(e))\setminus\{e\}$ 
because $e$ is isolated in $d^{-1}(d(e))$. 
We have 
$$\ip{\xi}{\pi_r(f)\xi}(d(e))
=\sum_{e'\in d^{-1}(d(e))}\overline{\xi(e')}f(r(e'))\xi(e')=f(r(e)).$$
From $t^0(\ip{\xi}{\pi_r(f)\xi})=t^1(\xi)^*t^0(f)t^1(\xi)\in I$ 
and $d(e)\in X_I^0$, 
we get $\ip{\xi}{\pi_r(f)\xi}(d(e))=0$. 
Thus we have $f(r(e))=0$. 
Hence $r(e)\in X_I^0$. 
Therefore $X_I^0$ is positively invariant.
\end{proof}

\begin{lemma}\label{IcapGn}
Let $I$ be an ideal of $\cO(E)$. 
\benu
\item For $\xi\in C_d(E^1)$, $t^1(\xi)\in I$ 
if and only if $\xi\in C_d(E^1\setminus X_I^1)$. 
\item For $x\in \cK(C_d(E^1))$, 
$\varphi(x)\in I$ 
if and only if $x\xi\in C_d(E^1\setminus X_I^1)$ 
for all $\xi\in C_d(E^1)$. 
\eenu
\end{lemma}

\begin{proof}
\benu
\item For $\xi\in C_d(E^1)$, 
$t^1(\xi)\in I$ if and only if 
$t^1(\xi)^*t^1(\xi)=t^0(\ip{\xi}{\xi})$ is in $I$, 
which is equivalent to $\ip{\xi}{\xi}\in C_0(E^0\setminus X_I^0)$ 
by the definition of $X_I^0$. 
By \cite[Lemma 1.12]{Ka1}, 
$\ip{\xi}{\xi}\in C_0(E^0\setminus X_I^0)$ 
is equivalent to $\xi\in C_d(E^1\setminus X_I^1)$. 
\item Take $x\in \cK(C_d(E^1))$ with $\varphi(x)\in I$. 
For $\xi\in C_d(E^1)$, 
we have $x\xi\in C_d(E^1\setminus X_I^1)$ by (i) 
because $t^1(x\xi)=\varphi(x)t^1(\xi)\in I$. 
Conversely suppose that $x\in \cK(C_d(E^1))$ satisfies 
$x\xi\in C_d(E^1\setminus X_I^1)$ 
for all $\xi\in C_d(E^1)$. 
Then $x$ is approximated by elements 
in the form $\sum_{k=1}^K\theta_{\xi_k,\eta_k}$ 
where $\xi_k,\eta_k\in C_d(E^1\setminus X_I^1)$ 
by \cite[Lemma 1.13]{Ka1}.
Since we have 
$$\varphi\bigg(\sum_{k=1}^K\theta_{\xi_k,\eta_k}\bigg)=
\sum_{k=1}^Kt^1(\xi_k)t^1(\eta_k)^*\in I$$
by (i), we get $\varphi(x)\in I$.
\eenu
\end{proof}

\begin{proposition}\label{saturated}
For an ideal $I$ of $\cO(E)$, 
the closed set $X_I^0$ is negatively invariant. 
\end{proposition}

\begin{proof}
Take $v\in \s{E}{0}{rg}$ such that $d(e)\notin X_I^0$ for all $e\in r^{-1}(v)$,
and we will prove that $v\notin X_I^0$.
By Lemma~\ref{RgPr}, 
we can find a neighborhood $V\subset \s{E}{0}{rg}$ of $v$ so that
$r^{-1}(V)\cap X_I^1=\emptyset$.
Take $f\in C_0(V)$ with $f(v)=1$.
For $\xi\in C_d(E^1)$ and $e\in X_I^1$, 
we have $(\pi_r(f)\xi)(e)=f(r(e))\xi(e)=0$.
Hence we see that $\pi_r(f)\xi\in C_d(E^1\setminus X_I^1)$ 
for all $\xi\in C_d(E^1)$. 
By Lemma~\ref{IcapGn} (ii), 
we have $\varphi(\pi_r(f))\in I$. 
Since $f\in C_0(\s{E}{0}{rg})$, 
we have $t^0(f)=\varphi(\pi_r(f))\in I$.
Therefore $v\notin X_I^0$. 
This shows that $X_I^0$ is negatively invariant. 
\end{proof}

\begin{proposition}\label{rhoI}
For an ideal $I$ of $\cO(E)$, 
the pair $\rho_I=(X_I^0,Z_I)$ is admissible.
\end{proposition}

\begin{proof}
We have already seen that $X_I^0$ is invariant 
in Proposition~\ref{hereditary} and Proposition~\ref{saturated}.
Since $t^0(C_0(\s{E}{0}{rg}))\subset\G^1$, 
we have $Z_I\cap \s{E}{0}{rg}=\emptyset$. 
We also have $Z_I\subset X_I^0$ 
because $t^0(C_0(E^0\setminus X_I^0))\subset I$. 
Hence we get $Z_I\subset \s{E}{0}{sg}\cap X_I^0$. 
We will show that $(X_I^0)_{\rs{sg}}\subset Z_I$. 
To derive a contradiction, suppose that 
there exists $v_0\in (X_I^0)_{\rs{sg}}$ with $v_0\notin Z_I$.
Then there exist $f\in C_0(E^0)$ with $t^0(f)\in I+\G^1$ 
and $f(v_0)\neq 0$. 
Take $x\in \cK(C_d(E^1))$
such that $t^0(f)+\varphi(x)\in I$. 
Since $v_0\in (X_I^0)_{\rs{sg}}$, 
either $v_0\in\overline{(X_I^0)_{\rs{sce}}}$
or $v_0\in (X_I^0)_{\rs{inf}}$.

First we consider the case that $v_0\in\overline{(X_I^0)_{\rs{sce}}}$.
We can find $v_1\in (X_I^0)_{\rs{sce}}$ with $f(v_1)\neq 0$.
Since $v_1\notin\overline{r(X_I^1)}$, 
there exists $g\in C_0(E^0)$ with $g(v_1)\neq 0$ 
and $g(v)=0$ for $v\in r(X_I^1)$. 
For $\xi\in C_d(E^1)$ and $e\in X_I^1$,
we have $(\pi_r(g)x\xi)(e)=g(r(e))\,(x\xi)(e)=0$.
Hence $(\pi_r(g)x)\xi\in C_d(E^1\setminus X_I^1)$ 
for all $\xi\in C_d(E^1)$.
By Lemma~\ref{IcapGn} (ii),
we have $\varphi(\pi_r(g)x)\in I$.
Since we have 
$$t^0(gf)+\varphi(\pi_r(g)x)=t^0(g)\big(t^0(f)+\varphi(x)\big)\in I,$$
we get $t^0(gf)\in I$. 
This contradicts the fact that $(gf)(v_1)\neq 0$ and $v_1\in X_I^0$.

Next we consider the case that $v_0\in (X_I^0)_{\rs{inf}}$.
We can find $\e>0$ and a neighborhood $V$ of $v_0$ 
such that $|f(v)|\geq\e$ for all $v\in V$.
There exists $x'\in\cK(C_d(E^1))$ satisfying that $\|x-x'\|<\e$ 
and $x'=\sum_{i=1}^k\theta_{\xi_i,\eta_i}$ for some 
$\xi_i,\eta_i\in C_c(E^1)$. 
We set 
$K=\bigcup_{i=1}^k\supp \eta_i$, 
which is a compact subset of $E^1$.
Since $v_0\in (X_I^0)_{\rs{inf}}$, 
we have $r^{-1}(V)\cap X_I^1\not\subset K$.
Take $e\in (r^{-1}(V)\cap X_I^1)\setminus K$.
We can find $\xi\in C_d(E^1)$ such that $\|\xi\|=1$, 
$\xi(e)=1$ and $K\cap \supp\xi=\emptyset$.
Since $t^1(\pi_r(f)\xi+x\xi)=\big(t^0(f)+\varphi(x)\big)t^1(\xi)\in I$,
we have $(\pi_r(f)\xi+x\xi)(e)=0$ by Lemma~\ref{IcapGn} (i).
Since $\|x\xi\|=\|(x-x')\xi\|<\e$, we have $|x\xi(e)|<\e$.
On the other hand, we have 
$$|(\pi_r(f)\xi)(e)|
=|f(r(e))\xi(e)|\geq\e.$$
This is a contradiction.
The proof is completed.
\end{proof}

For two admissible pairs $\rho_1=(X_1^0,Z_1),\ \rho_2=(X_2^0,Z_2)$, 
we write $\rho_1\subset\rho_2$ if $X_1^0\subset X_2^0$ and $Z_1\subset Z_2$. 
By the definition, we can see the following. 

\begin{lemma}\label{inclusion}
For two ideals $I_1,I_2$ of $\cO(E)$, 
$I_1\subset I_2$ implies $\rho_{I_1}\supset \rho_{I_2}$. 
\end{lemma}

For two admissible pairs $\rho_1=(X_1^0,Z_1),\ \rho_2=(X_2^0,Z_2)$, 
we denote by $\rho_1\cup\rho_2$ 
the pair $(X_1^0\cup X_2^0,Z_1\cup Z_2)$. 
We will show that the pair $\rho_1\cup\rho_2$ is admissible.

\begin{lemma}\label{cup}
For two invariant closed sets $X_1^0,X_2^0$, 
the closed set $X^0=X_1^0\cup X_2^0$ is invariant, 
and we have 
$\s{X}{0}{sce}\subset(X_1^0)_{\rs{sce}}\cup (X_2^0)_{\rs{sce}}$, 
$\s{X}{0}{inf}=(X_1^0)_{\rs{inf}}\cup (X_2^0)_{\rs{inf}}$ 
and $\s{X}{0}{sg}\subset(X_1^0)_{\rs{sg}}\cup (X_2^0)_{\rs{sg}}$.
\end{lemma}

\begin{proof}
Take $e\in E^1$ with $d(e)\in X^0$. 
Then either $d(e)\in X_1^0$ or $d(e)\in X_2^0$. 
Since $X_1^0$ and $X_2^0$ are positively invariant, 
either $r(e)\in X_1^0$ or $r(e)\in X_2^0$ holds. 
Hence $r(e)\in X^0$. 
Thus $X^0$ is positively invariant. 

Next we will show that 
$\s{X}{0}{sce}\subset(X_1^0)_{\rs{sce}}\cup (X_2^0)_{\rs{sce}}$. 
Take $v\in\s{X}{0}{sce}$. 
There exists a neighborhood $V$ of $v\in X^0$ 
such that $r^{-1}(V)\cap X^1=\emptyset$. 
When $v\in X_1^0$, 
the set $V\cap X_1^0$ is a neighborhood of $v\in X_1^0$ and 
we have $r^{-1}(V\cap X_1^0)\cap X_1^1=\emptyset$ 
because $X_1^1\subset X^1$. 
Hence we have $v\in (X_1^0)_{\rs{sce}}$.
Similarly we get $v\in (X_2^0)_{\rs{sce}}$ when $v\in X_2^0$. 
Thus we have shown that 
$\s{X}{0}{sce}\subset(X_1^0)_{\rs{sce}}\cup (X_2^0)_{\rs{sce}}$. 
Hence we have 
$$\s{X}{0}{sce}\cap\s{E}{0}{rg}\subset 
\big((X_1^0)_{\rs{sce}}\cap\s{E}{0}{rg}\big)
\cup \big((X_2^0)_{\rs{sce}}\cap\s{E}{0}{rg}\big)=\emptyset$$ 
because $X_1^0$ and $X_2^0$ are negatively invariant. 
By Proposition~\ref{condforNI}, $X^0$ is negatively invariant. 

Next we will prove 
$\s{X}{0}{inf}=(X_1^0)_{\rs{inf}}\cup (X_2^0)_{\rs{inf}}$. 
Take $v\in (X_1^0)_{\rs{inf}}$. 
For any closed neighborhood $V$ of $v\in X^0$, 
$V\cap X_1^0$ is a closed neighborhood of $v\in X_1^0$. 
Since $v\in (X_1^0)_{\rs{inf}}$,
we have that $r^{-1}(V\cap X_1^0)\cap X_1^1$ is not compact. 
From $r^{-1}(V)\cap X^1\supset r^{-1}(V\cap X_1^0)\cap X_1^1$, 
we see that $r^{-1}(V)\cap X^1$ is not compact.
Hence $v\in \s{X}{0}{inf}$. 
Similarly if $v\in (X_2^0)_{\rs{inf}}$ then $v\in\s{X}{0}{inf}$. 
Thus we have shown that $v\in\s{X}{0}{inf}$ 
for every $v\in (X_1^0)_{\rs{inf}}\cup (X_2^0)_{\rs{inf}}$. 
Conversely take $v\notin (X_1^0)_{\rs{inf}}\cup (X_2^0)_{\rs{inf}}$, 
and we will show that $v\notin\s{X}{0}{inf}$. 
From $v\notin (X_1^0)_{\rs{inf}}$, 
we have either $v\notin X_1^0$ or $v\in (X_1^0)_{\rs{fin}}$. 
For both cases, we can find a compact neighborhood $V_1$ of $v\in E^0$ 
such that $r^{-1}(V_1\cap X_1^0)\cap X_1^1$ is compact. 
Similarly there exists a compact neighborhood $V_2$ of $v\in E^0$ 
such that $r^{-1}(V_2\cap X_2^0)\cap X_2^1$ is compact. 
Set $V=V_1\cap V_2$ 
which is a compact neighborhood of $v\in E^0$. 
We have 
\begin{align*}
r^{-1}(V)\cap X^1
&=(r^{-1}(V)\cap X_1^1)\cup (r^{-1}(V)\cap X_2^1)\\
&=(r^{-1}(V\cap X_1^0)\cap X_1^1)\cup (r^{-1}(V\cap X_2^0)\cap X_2^1)\\
&\subset (r^{-1}(V_1\cap X_1^0)\cap X_1^1)
 \cup (r^{-1}(V_2\cap X_2^0)\cap X_2^1).
\end{align*}
Hence $r^{-1}(V)\cap X^1$ is compact. 
Thus $v\notin \s{X}{0}{inf}$. 
Therefore 
we have $\s{X}{0}{inf}=(X_1^0)_{\rs{inf}}\cup (X_2^0)_{\rs{inf}}$. 
Now it is easy to see that 
$\s{X}{0}{sg}\subset(X_1^0)_{\rs{sg}}\cup (X_2^0)_{\rs{sg}}$. 
\end{proof}

Note that in general 
$X_1^0\subset X_2^0$ does not imply 
$(X_1^0)_{\rs{sg}}\subset (X_2^0)_{\rs{sg}}$ 
for two invariant closed sets $X_1^0,X_2^0$. 

\begin{example}
Let $E=(E^0,E^1,d,r)$ be the discrete graph 
given by 
$$
\begin{array}{cc}
E^0=\{v,v',w\},&
E^1=\{e_k\}_{k\in\N},\\
d(e_k)=\left\{\begin{array}{ll}
v& (k=0)\\
v'& (k\geq 1)
\end{array}\right. ,&
r(e_k)=w\ (k\in\N).
\end{array}
\qquad
\raisebox{1cm}{\xymatrix{
&\bullet\ v' \ar@<-1.8ex>[d]_{e_1} \ar@<-1.0ex>[d]^{\cdots}\\
v\ \bullet \ar[r]_{e_0} & \bullet\ w}}
$$
This example is the same as in \cite[Example 4.9]{Ka2}. 
Since $\s{E}{0}{sg}=E^0$, 
every subset of $E^0$ is negatively invariant. 
The two sets $X_1^0=\{w\}$ and $X_2^0=\{v,w\}$ are 
closed invariant sets satisfying $X_1^0\subset X_2^0$. 
However we do not have $(X_1^0)_{\rs{sg}}\subset (X_2^0)_{\rs{sg}}$ 
because $(X_1^0)_{\rs{sg}}=\{w\}$ 
and $(X_2^0)_{\rs{sg}}=\{v\}$. 
\end{example}

\begin{proposition}
For two admissible pairs $\rho_1=(X_1^0,Z_1),\ \rho_2=(X_2^0,Z_2)$, 
the pair $\rho_1\cup\rho_2$ is admissible. 
\end{proposition}

\begin{proof}
By Lemma~\ref{cup}, 
$X^0=X_1^0\cup X_2^0$ is invariant
and 
$$\s{X}{0}{sg}\subset(X_1^0)_{\rs{sg}}\cup (X_2^0)_{\rs{sg}}
\subset Z_1\cup Z_2
\subset (\s{E}{0}{sg}\cap X_1^0)\cup(\s{E}{0}{sg}\cap X_2^0)
=\s{E}{0}{sg}\cap X^0.$$
Therefore $\rho_1\cup\rho_2$ is admissible. 
\end{proof}

\begin{proposition}\label{cap}
For two ideals $I_1,I_2$ of $\cO(E)$, 
we have $\rho_{I_1\cap I_2}=\rho_{I_1}\cup\rho_{I_2}$. 
\end{proposition}

\begin{proof}
This follows from the following computations 
\begin{align*}
I_1\cap I_2\cap\G^0
&=(I_1\cap\G^0)\cap (I_2\cap\G^0)\\
\big(I_1\cap I_2+\G^1\big)\cap\G^0
&=\big((I_1+\G^1)\cap(I_2+\G^1)\big)\cap\G^0\\
&=\big((I_1+\G^1)\cap\G^0\big)\cap\big((I_2+\G^1)\cap\G^0\big)
\end{align*}
and the remark after Definition~\ref{XIZI}. 
\end{proof}

\section{Gauge-invariant ideals}\label{SecGII}

In Section~\ref{adpair},
we get admissible pairs $\rho_I$ from ideals $I$ of $\cO(E)$. 
Conversely we can construct ideals $I_{\rho}$ 
from admissible pairs $\rho$. 

We set $\F^1=\G^0+\G^1\subset \cO(E)$ which is a $C^*$-subalgebra. 
We recall the description of $\F^1$ done in \cite[Proposition 5.2]{Ka1}. 
Two $*$-homomorphisms $\pi_0^1\colon \F^1\to C_0(\s{E}{0}{sg})$ 
and $\pi_1^1\colon \F^1\to\cL(C_d(E^1))$ are defined by 
$\pi_0^1(a)=f|_{\s{E}{0}{sg}}$ and $\pi_1^1(a)=\pi_r(f)+x$
for $a=t^0(f)+\varphi(x)\in\F^1$. 
Note that $\pi_0^1$ is a surjective map whose kernel is $\G^1$, 
and that the restriction of $\pi_1^1$ to $\G^1$ 
is the inverse map of $\varphi$. 
Thus 
$$\pi_0^1\oplus\pi_1^1\colon 
\F^1\to C_0(\s{E}{0}{sg})\oplus \cL(C_d(E^1))$$ 
is injective. 
Note that we have $at^1(\xi)=t^1(\pi_1^1(a)\xi)$ 
for $a\in\F^1$ and $\xi\in C_d(E^1)$, 
and that for an ideal $I$ of $\cO(E)$ 
the closed set $Z_I$ is determined by 
$C_0(\s{E}{0}{sg}\setminus Z_I)=\pi_0^1(I\cap\F^1)$.

Let $X^0$ be a closed positively invariant subset of $E^0$. 
We get the topological graph $X=(X^0,X^1,d_X,r_X)$ 
where $X^1=d^{-1}(X^0)\subset E^1$ 
and $d_X,r_X$ are the restrictions of $d,r$ to $X^1$. 
There is a natural surjection $C_d(E^1)\ni\xi\mapsto \dot{\xi}\in C_{d_X}(X^1)$
whose kernel is $C_d(E^1\setminus X^1)$ (\cite[Lemma 1.12]{Ka1}).
Since the submodule $C_d(E^1\setminus X^1)$ is invariant 
under the action of $\cL(C_d(E^1))$,
we can define a $*$-homomorphism 
$\omega_X\colon \cL(C_d(E^1))\to \cL(C_{d_X}(X^1))$ 
by $\omega_X(x)\dot{\xi}=\dot{(x\xi)}$. 
For $f\in C_0(E^0)$, we have 
$\omega_X(\pi_r(f))=\pi_{r_{X}}(f|_{X^0})$. 
One can easily see that 
$$\ker\omega_X=\{x\in \cL(C_d(E^1))\mid 
x\xi\in C_d(E^1\setminus X^1)\mbox{ for all }
\xi\in C_d(E^1)\}.$$
It is easy to see that 
the restriction of $\omega_X$ to $\cK(C_d(E^1))$ 
is a surjection onto $\cK(C_{d_X}(X^1))$ 
(\cite[Lemma 1.14]{Ka1}). 

\begin{definition}
For an admissible pair $\rho=(X^0,Z)$, 
we define an ideal $J_{\rho}$ of $\F^1$ by
\begin{align*}
J_{\rho}&=\{a\in\F^1\mid 
\pi_0^1(a)\in C_0(\s{E}{0}{sg}\setminus Z),\ 
\pi_1^1(a)\in\ker\omega_X\}\\
&=\{t^0(f)+\varphi(x)\in\F^1\mid 
f\in C_0(E^0\setminus Z),\ \omega_X(\pi_r(f)+x)=0\}.
\end{align*}
\end{definition}

\begin{lemma}\label{JtX}
The ideal $J_{\rho}$ of $\F^1$ satisfies 
$$\pi_0^1(J_{\rho})=C_0(\s{E}{0}{sg}\setminus Z), \mbox{ and }
J_{\rho}\cap\G^0=t^0(C_0(E^0\setminus X^0)).$$
\end{lemma}

\begin{proof}
By definition, $\pi_0^1(J_{\rho})\subset C_0(\s{E}{0}{sg}\setminus Z)$.
We will show the other inclusion 
$\pi_0^1(J_{\rho})\supset C_0(\s{E}{0}{sg}\setminus Z)$.
Take $g\in C_0(\s{E}{0}{sg}\setminus Z)$. 
Choose $f\in C_0(E^0\setminus Z)$ such that $f|_{\s{E}{0}{sg}}=g$. 
Since $\s{X}{0}{sg}\subset Z$, 
we have $f|_{X^0}\in C_0(\s{X}{0}{rg})$. 
Thus we see $\omega_X(\pi_r(f))=\pi_{r_{X}}(f|_{X^0})\in\cK(C_{d_X}(X^1))$. 
Hence there exists $x\in\cK(C_d(E^1))$ 
such that $\omega_X(x)=\omega_X(\pi_r(f))$.
Thus we get $a=t^0(f)-\varphi(x)\in J_{\rho}$ with $\pi_0^1(a)=g$.
Therefore $\pi_0^1(J_{\rho})=C_0(\s{E}{0}{sg}\setminus Z)$.

For $f\in C_0(E^0)$, 
$t^0(f)\in J_{\rho}$ if and only if $f\in C_0(E^0\setminus Z)$ 
and $\pi_r(f)\in\ker\omega_X$. 
Now $\pi_r(f)\in\ker\omega_X$ if and only if 
$f(v)=0$ for all $v\in r(X^1)$. 
Since $\s{X}{0}{sg}\subset Z\subset X^0$ 
and $\s{X}{0}{rg}\subset r(X^1)\subset X^0$, 
we have $Z\cup r(X^1)=X^0$.
Therefore, $t^0(f)\in J_{\rho}$ 
if and only if $f\in C_0(E^0\setminus X^0)$.
\end{proof}

\begin{definition}
For an admissible pair $\rho=(X^0,Z)$, 
we define a subset $I_{\rho}$ of $\cO(E)$ 
to be the closure of 
$$\spa\{t^n(\xi)at^m(\eta)^*\mid a\in J_{\rho},\ 
\xi\in C_d(E^n),\ \eta\in C_d(E^m),\ n,m\in\N\}.$$
\end{definition}

We will see that the set $I_{\rho}$ is the ideal generated by $J_{\rho}$. 
To prove it, we need the following lemma. 

\begin{lemma}\label{ItX0}
For $a\in J_{\rho}$ and $\zeta\in C_d(E^k)$ with $k\geq 1$, 
there exist $\zeta'\in C_d(E^k)$ and $f\in C_0(E^0\setminus X^0)$ 
such that $at^k(\zeta)=t^k(\zeta')t^0(f)$. 
\end{lemma}

\begin{proof}
Let us set $X^k=(d^k)^{-1}(X^0)\subset E^k$. 
For $\zeta''\in C_d(E^k\setminus X^k)$, 
we can find $\zeta'\in C_d(E^k)$ and $f\in C_0(E^0\setminus X^0)$ 
with $t^k(\zeta'')=t^k(\zeta')t^0(f)$ by \cite[Lemma 1.12]{Ka1}. 
Hence it suffices to show $at^k(\zeta)\in t^k(C_d(E^k\setminus X^k))$. 
We may assume $\zeta=\xi\otimes\eta$ 
for $\xi\in C_d(E^1)$ and $\eta\in C_d(E^{k-1})$. 
We have 
$$at^k(\zeta)=t^k((\pi_1^1(a)\xi)\otimes\eta).$$
Since $\pi_1^1(a)\in\ker \omega_X$, 
we have $\pi_1^1(a)\xi\in C_d(E^1\setminus X^1)$. 
Since $X^0$ is positively invariant, 
$(e_1,\ldots,e_k)\in X^k$ implies $e_1\in X^1$. 
Hence we get 
$(\pi_1^1(a)\xi)\otimes\eta\in C_d(E^k\setminus X^k)$.
We are done. 
\end{proof}

\begin{proposition}\label{IisGaugeInv}
The set $I_{\rho}$ is a gauge-invariant ideal of $\cO(E)$. 
\end{proposition}

\begin{proof}
By definition, $I_{\rho}$ is a closed $*$-invariant subspace of $\cO(E)$,
which is invariant under the gauge action. 
To prove that $I_{\rho}$ is an ideal, 
it suffices to see that $xy\in I_{\rho}$ 
for $x=t^n(\xi)at^m(\eta)^*\in I_{\rho}$ 
and $y=t^{n'}(\xi')t^{m'}(\eta')^*\in \cO(E)$ 
where $a\in J_{\rho}$. 
When $m\geq n'$, obviously $xy\in I_{\rho}$. 
Hence we only need to show that 
$t^n(\xi)at^{k}(\zeta)t^{m'}(\eta')^*\in I_{\rho}$ 
for some $k\geq 1$ and $\zeta\in C_d(E^k)$.
By Lemma~\ref{ItX0}, 
there exist $\zeta'\in C_d(E^k)$ and $f\in C_0(E^0\setminus X^0)$
such that $at^{k}(\zeta)=t^{k}(\zeta')t^0(f)$.
Thus we get
$$t^n(\xi)at^{k}(\zeta)t^{m'}(\eta')^*
=t^{n+k}(\xi\otimes\zeta')t^0(f)t^{m'}(\eta')^*\in I_{\rho},$$
because $t^0(f)\in J_{\rho}$ by Lemma~\ref{JtX}.
This completes the proof.
\end{proof}

We will show that for an admissible pair $\rho$, 
we have $\rho_{I_{\rho}}=\rho$. 
To this end, we need a series of lemmas.

\begin{lemma}\label{ItX2}
There exists a norm-decreasing map 
$\pi_0\colon\cO(E)\to C_0(\s{E}{0}{sg})$ 
satisfying that $\pi_0(t^0(f))=f|_{\s{E}{0}{sg}}$ for $f\in C_0(E^0)$ 
and $\pi_0(t^n(\xi)t^m(\eta)^*)=0$ when $n\geq 1$ or $m\geq 1$. 
\end{lemma}

\begin{proof}
Let $\F\subset \cO(E)$ be 
the fixed point algebra of the gauge action $\beta$, 
and $\varPsi\colon \cO(E)\to\F$ 
be the conditional expectation defined 
by $\varPsi(x)=\int_{\T}\beta_z(x)dz$ for $x\in \cO(E)$. 
For $\xi\in C_d(E^n)$ and $\eta\in C_d(E^m)$, 
we have 
$$\varPsi(t^{n}(\xi)t^{m}(\eta)^*)=\delta_{n,m}t^{n}(\xi)t^{m}(\eta)^*.$$
We can define a $*$-homomorphism 
$\pi_0'\colon \F\to C_0(\s{E}{0}{sg})$ 
by $\pi_0'(t^0(f))=f|_{\s{E}{0}{sg}}$ for $f\in C_0(E^0)$ 
and $\pi_0'(t^n(\xi)t^n(\eta)^*)=0$ when $n\geq 1$ 
(see the remark after \cite[Lemma 5.1]{Ka1}). 
The composition of $\varPsi$ and $\pi_0'$ 
satisfies the desired properties.
\end{proof}

\begin{lemma}\label{ItX3}
We have $\pi_0^1(I_{\rho}\cap\F^1)=C_0(\s{E}{0}{sg}\setminus Z)$.
\end{lemma}

\begin{proof}
Since $I_{\rho}\cap\F^1\supset J_{\rho}$,
we have $\pi_0^1(I_{\rho}\cap\F^1)\supset C_0(\s{E}{0}{sg}\setminus Z)$
by Lemma~\ref{JtX}. 
We will prove the other inclusion 
$\pi_0^1(I_{\rho}\cap\F^1)\subset C_0(\s{E}{0}{sg}\setminus Z)$.
Since the restriction of the map $\pi_0$ in Lemma~\ref{ItX2} 
to $\F^1$ is $\pi_0^1$, 
we have $\pi_0^1(I_{\rho}\cap\F^1)\subset \pi_0(I_{\rho})$. 
For $x=t^{n}(\xi)at^{m}(\eta)^*\in I_{\rho}$, 
we have $\pi_0(x)=0$ when $n\geq 1$ or $m\geq 1$ 
and 
$$\pi_0(x)=(\xi|_{\s{E}{0}{sg}})\pi_0^1(a)(\eta|_{\s{E}{0}{sg}})
\in C_0(\s{E}{0}{sg}\setminus Z)$$ 
when $n=m=0$. 
Hence $\pi_0(I_{\rho})\subset C_0(\s{E}{0}{sg}\setminus Z)$.
Thus we get $\pi_0^1(I_{\rho}\cap\F^1)=C_0(\s{E}{0}{sg}\setminus Z)$.
\end{proof}

\begin{lemma}\label{ItX4}
Let $v\in X^0$. 
Either there exist $n\in\N$ and $e\in E^n$ 
with $r^n(e)=v$ and $d^n(e)\in Z$, 
or for all $n\in\N$ we can find $e\in E^n$ 
with $r^n(e)=v$ and $d^n(e)\in X^0$.
\end{lemma}

\begin{proof}
If $v\in Z$, 
the first alternative holds with $n=0$ and $e=v\in E^0$. 
If $v\in X^0\setminus Z$, 
then $v\in \s{X}{0}{rg}$
because $\s{X}{0}{sg}\subset Z$. 
Hence we can find $e_1\in X^1$ 
such that $r(e_1)=v$ and $d(e_1)\in X^0$ by Lemma~\ref{RgPr}.
If $d(e_1)\in Z$, 
the first alternative holds with $n=1$ and $e=e_1\in E^1$. 
If $d(e_1)\in X^0\setminus Z$, 
the argument above shows that 
there exists $e_2\in X^1$ 
such that $r(e_2)=d(e_1)$ and $d(e_2)\in X^0$. 
If $d(e_2)\in Z$, 
the first alternative holds with $n=2$ and $e=(e_1,e_2)\in E^2$. 
Repeating this argument, 
either we get the first alternative, 
or we can find $e_k\in X^1$ for $k=1,2,\ldots$ 
such that $r(e_1)=v$, $r(e_{k+1})=d(e_k)$ and $d(e_k)\in X^0$ 
for all $k$. 
The latter situation implies the second alternative of this lemma. 
\end{proof}

\begin{lemma}\label{ItX5}
We have $I_{\rho}\cap\G^0=t^0(C_0(E^0\setminus X^0))$.
\end{lemma}

\begin{proof}
We have $I_{\rho}\cap\G^0\supset J_{\rho}\cap\G^0=t^0(C_0(E^0\setminus X^0))$
by Lemma~\ref{JtX}.
We will prove the other inclusion. 
To the contrary, 
suppose that there exists $f\in C_0(E^0)$ such that $t^0(f)\in I_{\rho}$ 
and $f(v)=1$ for some $v\in X^0$.
By Lemma~\ref{ItX4},
either there exist $n\in\N$ and $e\in E^n$ 
with $r^n(e)=v$ and $d^n(e)\in Z$, 
or for all $n\in\N$ we can find $e\in E^n$ 
with $r^n(e)=v$ and $d^n(e)\in X^0$.

We first consider the case that there exist $n\in\N$ and $e\in E^n$ 
with $r^n(e)=v$ and $d^n(e)\in Z$. 
Take a neighborhood $U$ of $e\in E^n$ such that the restriction of $d^n$ 
to $U$ is injective, and take $\zeta\in C_c(U)$ with $\zeta(e)=1$. 
We set $g=\ip{\zeta}{\pi_{r^n}(f)\zeta}\in C_0(E^0)$. 
Then we have $t^0(g)=t^n(\zeta)^*t^0(f)t^n(\zeta)\in I_{\rho}$ and 
$$g(d^n(e))=\sum_{e'\in (d^n)^{-1}(d^n(e))}
\overline{\zeta(e')}f(r^n(e'))\zeta(e')
=\overline{\zeta(e)}f(r^n(e))\zeta(e)=1.$$
This contradicts Lemma~\ref{ItX3}.

Next we consider the case that for all $n\in\N$ we can find $e\in E^n$ 
with $r^n(e)=v$ and $d^n(e)\in X^0$. 
Take $\xi_k\in C_d(E^{n_k})$, $\eta_k\in C_d(E^{m_k})$ 
and $a_k\in J_{\rho}$ such that 
$$\bigg\|t^0(f)-\sum_{k=1}^Kt^{n_k}(\xi_k)a_kt^{m_k}(\eta_k)^*\bigg\|
<\frac12.$$
Since the conditional expectation $\varPsi\colon \cO(E)\to\F$ 
defined in the proof of Lemma~\ref{ItX2} satisfies 
$\varPsi(t^0(f))=t^0(f)$ and 
$$\varPsi(t^{n_k}(\xi_k)a_kt^{m_k}(\eta_k)^*)
=\delta_{n_k,m_k}t^{n_k}(\xi_k)a_kt^{m_k}(\eta_k)^*,$$
we have 
$$\bigg\|t^0(f)-\sum_{n_k=m_k}t^{n_k}(\xi_k)a_kt^{m_k}(\eta_k)^*\bigg\|
<\frac12.$$
Take $n\in\N$ such that $n>n_k$ for $k$ with $n_k=m_k$. 
We can find $e\in E^n$ with $r^n(e)=v$ and $d^n(e)\in X^0$.
Take a neighborhood $U$ of $e\in E^n$ such that the restriction of $d^n$ 
to $U$ is injective and take $\zeta\in C_c(U)$ with $0\leq\zeta\leq 1$ and 
$\zeta(e)=1$. 
We have $\|\zeta\|=1$.
Hence we get 
$$\bigg\|t^n(\zeta)^*
\bigg(t^0(f)-\sum_{n_k=m_k}t^{n_k}(\xi_k)a_kt^{m_k}(\eta_k)^*\bigg)
t^n(\zeta)\bigg\|<\frac12.$$
We set $g=\ip{\zeta}{\pi_{r^n}(f)\zeta}\in C_0(E^0)$. 
By Lemma~\ref{ItX0}, 
we have 
$$t^n(\zeta)^*t^{n_k}(\xi_k)a_kt^{m_k}(\eta_k)^*t^n(\zeta)
\in t^0(C_0(E^0\setminus X^0))$$ 
for $k$ with $n_k=m_k$. 
Hence $|g(d^n(e))|<1/2$. 
However, we can prove $g(d^n(e))=1$ similarly as above.
This is a contradiction. 

Thus, we have shown that $I_{\rho}\cap\G^0=t^0(C_0(E^0\setminus X^0))$.
\end{proof}

\begin{proposition}\label{XIX}
For an admissible pair $\rho$, we have $\rho_{I_{\rho}}=\rho$.
\end{proposition}

\begin{proof}
This follows from the equation 
$C_0(\s{E}{0}{sg}\setminus Z_I)=\pi_0^1(I\cap\F^1)$ 
explained in the beginning of this section, 
Lemma~\ref{ItX3} and Lemma~\ref{ItX5}.
\end{proof}

We have shown that $\rho_{I_{\rho}}=\rho$ 
for every admissible pair $\rho$. 
Next we study for which ideal $I$ of $\cO(E)$ 
the equation $I_{\rho_I}=I$ holds. 
Clearly the condition that $I$ is gauge-invariant is necessary.
This condition is shown to be sufficient (Proposition~\ref{IXI}). 
For general ideals, we have the following. 

\begin{lemma}\label{IXIsubsetI}
For an ideal $I$ of $\cO(E)$, 
we have $I\cap\F^1=J_{\rho_I}$ and $I\supset I_{\rho_I}$. 
\end{lemma}

\begin{proof}
Take $a\in I\cap\F^1$. 
Since the closed set $Z_I$ is determined by 
$C_0(\s{E}{0}{sg}\setminus Z_I)=\pi_0^1(I\cap\F^1)$, 
we have $\pi_0^1(a)\in C_0(\s{E}{0}{sg}\setminus Z_I)$. 
For $\xi\in C_d(E^1)$, 
we have $\pi_1^1(a)\xi\in C_d(E^1\setminus X_I^1)$ 
by Lemma~\ref{IcapGn} (i) 
because $t^1(\pi_1^1(a)\xi)=at^1(\xi)\in I$. 
This shows that $\pi_1^1(a)\in\ker\omega_X$. 
Hence $a\in J_{\rho_I}$. 

Conversely take $b\in J_{\rho_I}$. 
Since $\pi_0^1(b)\in C_0(\s{E}{0}{sg}\setminus Z_I)$, 
we can find $a\in I\cap\F^1$ such that 
$\pi_0^1(b)=\pi_0^1(a)$. 
Since $b-a\in \ker\pi_0^1$, 
we can find $x\in\cK(C_d(E^1))$ with $\varphi(x)=b-a$. 
We have $x=\pi_1^1(b-a)$. 
By the former part of this proof, 
we have $\pi_1^1(b-a)\in\ker\omega_X$. 
Hence $x\xi\in C_d(E^1\setminus X_I^1)$ for all $\xi\in C_d(E^1)$. 
By Lemma~\ref{IcapGn} (ii), 
we have $\varphi(x)\in I$. 
Hence $b=a+\varphi(x)\in I$. 
We have shown $I\cap\F^1=J_{\rho_I}$. 

This implies $I\supset J_{\rho_I}$. 
Since $I_{\rho_I}$ is an ideal generated by $J_{\rho_I}$, 
we have $I\supset I_{\rho_I}$. 
\end{proof}

\begin{corollary}\label{IF=J}
We have $I_{\rho}\cap\F^1=J_{\rho}$ for an admissible pair $\rho$. 
\end{corollary}

\begin{proof}
Clear from Proposition~\ref{XIX} and Lemma~\ref{IXIsubsetI}. 
\end{proof}

\begin{definition}\label{Erho}
For an admissible pair $\rho=(X^0,Z)$, 
we define a topological graph 
$E_{\rho}=(E_{\rho}^0,E_{\rho}^1,d_{\rho},r_{\rho})$ as follows. 
Set $Y_\rho=Z\cap\s{X}{0}{rg}$ and define 
$$E_{\rho}^0=X^0\mathop{\amalg}_{\partial Y_\rho}\overline{Y_\rho}\ ,\qquad
E_{\rho}^1=X^1\!\!\!\mathop{\amalg}_{d^{-1}(\partial Y_\rho)}\!\!
d^{-1}(\overline{Y_\rho}),$$
where $\partial(Y_\rho)=\overline{Y_\rho}\setminus Y_\rho$. 
The domain map $d_{\rho}\colon E_{\rho}^1\to E_{\rho}^0$ is defined 
from $d\colon X^1\to X^0$ and 
$d\colon d^{-1}(\overline{Y_\rho})\to \overline{Y_\rho}$.
The range map $r_{\rho}\colon E_{\rho}^1\to E_{\rho}^0$ is defined 
from $r\colon X^1\to X^0$ and 
$r\colon d^{-1}(\overline{Y_\rho})\to X^0$ 
(see \cite[Section 3]{Ka2} for the detail).
\end{definition}

By definition, $E_{\rho}=X$ for $\rho=(X^0, \s{X}{0}{sg})$. 

\begin{remark}
For an admissible pair $\rho$, 
the four inclusions $X^0\to E^0$, $\overline{Y_\rho}\to E^0$, 
$X^1\to E^1$ and $d^{-1}(\overline{Y_\rho})\to E^1$ 
define continuous proper maps $m^0\colon E_{\rho}^0\to E^0$ 
and $m^1\colon E_{\rho}^1\to E^1$. 
We can see that the pair $m=(m^0,m^1)$ is a regular factor map 
in the sense of \cite[Definitions 2.1 and 2.6]{Ka2}, 
and that the ideal $I_{\rho}$ is the kernel of 
the surjection $\mu\colon \cO(E)\to\cO(E_{\rho})$ 
induced by $m=(m^0,m^1)$ \cite[Proposition 2.9]{Ka2} 
(see Proposition~\ref{IXI}). 
\end{remark}

\begin{proposition}\label{E_Ipair}
For an ideal $I$ of $\cO(E)$, 
we have a natural surjection $\cO(E_{\rho_I})\to \cO(E)/I$ 
which is injective on the image of the map 
$t_{\rho_I}^0\colon C_0(E_{\rho_I}^0)\to \cO(E_{\rho_I})$. 
\end{proposition}

\begin{proof}
It suffices to show that there exists an injective Cuntz-Krieger 
$E_{\rho_I}$-pair $\widetilde{T}$ on $\cO(E)/I$ 
with $C^*(\widetilde{T})=\cO(E)/I$. 
We use the construction in \cite[Section~3]{Ka2} 
(see \cite[Remark~3.28]{Ka2}). 

Let $\omega$ be the natural surjection from $\cO(E)$ onto $\cO(E)/I$ 
and define $T^i=\omega\circ t^i$ for $i=0,1$. 
The pair $T=(T^0,T^1)$ is a Cuntz-Krieger $E$-pair 
on $\cO(E)/I$ satisfying that $C^*(T)=\cO(E)/I$ and 
$\ker T^0=C_0(E^0\setminus X_I^0)$. 
Hence we can define an injective Toeplitz $X$-pair 
$\dot{T}=(\dot{T}^0,\dot{T}^1)$ 
by $\dot{T}^0(f|_{X^0})=T^0(f)$ and $\dot{T}^1(\xi|_{X^1})=T^1(\xi)$ 
for $f\in C_0(E^0)$ and $\xi\in C_d(E^1)$. 
Let $\dot{\varPhi}\colon \cK(C_d(X_I^1))\to \cO(E)/I$ 
be the $*$-homomorphism defined by $\dot{\varPhi}(\theta_{\xi,\eta})
=\dot{T}^1(\xi)\dot{T}^1(\eta)^*$ 
for $\xi,\eta\in C_d(X_I^1)$. 
Then we have $\dot{\varPhi}(\omega_X(x))=\omega(\varphi(x))$ 
for $x\in\cK(C_d(E^1))$. 
Take $g\in C_0(X_I^0)$. 
We will show that $\dot{T}^0(g)\in \dot{\varPhi}(\cK(C_d(X_I^1)))$ 
if and only if $g\in C_0(X_I^0\setminus Z_I)$. 
Take $f\in C_0(E^0)$ with $f|_{X_I^0}=g$. 
Then $\dot{T}^0(g)\in \dot{\varPhi}(\cK(C_d(X_I^1)))$ 
if and only if there exist $x\in \cK(C_d(E^1))$ 
such that $t^0(f)-\varphi(x)\in I$. 
This is equivalent to $f\in C_0(E^0\setminus Z_I)$. 
Hence we have shown that for $g\in C_0(X_I^0)$, 
$\dot{T}^0(g)\in \dot{\varPhi}(\cK(C_d(X_I^1)))$ 
if and only if $g\in C_0(X_I^0\setminus Z_I)$. 
This implies that the set $Y_T$ defined in \cite[Definition 3.1]{Ka2}
coincides with $Z_I\cap (X_I^0)_{\rs{rg}}=Y_{\rho_I}$ 
by \cite[Lemma 3.2]{Ka2}. 
Hence by \cite[Proposition 3.8]{Ka2}, 
the Cuntz-Krieger $E_{\rho_I}$-pair $\widetilde{T}$ 
constructed as in \cite[Proposition 3.15]{Ka2} is injective. 
This completes the proof.
\end{proof}

\begin{proposition}\label{IXI}
For an ideal $I$ of $\cO(E)$, 
the following are equivalent:
\benu
\item $I$ is gauge-invariant. 
\item The surjection $\cO(E_{\rho_I})\to\cO(E)/I$ in Proposition~\ref{E_Ipair} 
is an isomorphism. 
\item $I=I_{\rho_I}$. 
\eenu
\end{proposition}

\begin{proof}
(i)$\Rightarrow$(ii): 
When $I$ is a gauge-invariant ideal, 
the injective Cuntz-Krieger $E_{\rho_I}$-pair $\widetilde{T}$ 
defined in the proof of Proposition~\ref{E_Ipair} 
admits a gauge action by \cite[Lemma 3.17]{Ka2}. 
Hence the Gauge Invariant Uniqueness Theorem (\cite[Theorem 4.5]{Ka1}) 
implies that the natural surjection $\cO(E_{\rho_I})\to\cO(E)/I$ 
is an isomorphism 
(see also \cite[Proposition 3.27]{Ka2}). 

(ii)$\Rightarrow$(iii): 
Set $J=I_{\rho_I}$ which is gauge-invariant 
and satisfies $\rho_J=\rho_I$ by Proposition~\ref{XIX}. 
Hence there exists a surjection $\cO(E_{\rho_I})\to\cO(E)/J$ 
by Proposition~\ref{E_Ipair}. 
Since $J\subset I$ by Lemma~\ref{IXIsubsetI}, 
there exists a surjection $\cO(E)/J\to\cO(E)/I$. 
The composition of the two surjections 
is nothing but the surjection in Proposition~\ref{E_Ipair}, 
which is an isomorphism by (ii). 
Hence the surjection $\cO(E)/J\to\cO(E)/I$ is 
also an isomorphism. 
This shows $I=J=I_{\rho_I}$. 

(iii)$\Rightarrow$(i): 
Clear by Proposition~\ref{IisGaugeInv}. 
\end{proof}

\begin{proposition}\label{TopFreeQ}
Let $I$ be an ideal of $\cO(E)$. 
If the topological graph $E_{\rho_I}$ is topologically free, 
then $\cO(E)/I\cong \cO(E_{\rho_I})$. 
Hence we have $I=I_{\rho_I}$ and $I$ is gauge-invariant. 
\end{proposition}

\begin{proof}
When the topological graph $E_{\rho_I}$ is topologically free, 
the Cuntz-Krieger Uniqueness Theorem (Proposition~\ref{CKUT}) implies that 
the surjection $\cO(E_{\rho_I})\to\cO(E)/I$ in Proposition~\ref{E_Ipair} 
is an isomorphism. 
Hence the conclusion follows from Proposition~\ref{IXI}. 
\end{proof}

The following strengthens Lemma~\ref{IXIsubsetI}. 

\begin{proposition}\label{IrI}
For an ideal $I$ of $\cO(E)$, 
we have $I_{\rho_I}=\bigcap_{z\in\T}\beta_z(I)$. 
\end{proposition}

\begin{proof}
The ideal $\bigcap_{z\in\T}\beta_z(I)$ is gauge-invariant, 
and 
$$\bigcap_{z\in\T}\beta_z(I)\cap\F^1=\bigcap_{z\in\T}\beta_z(I\cap\F^1)
=I\cap\F^1=J_{\rho_I}.$$
Hence $\bigcap_{z\in\T}\beta_z(I)=I_{\rho_I}$ 
by Proposition~\ref{IXI}. 
\end{proof}

We get the main result of this section. 

\begin{theorem}\label{gaugeideals}
The maps $I\mapsto\rho_I$ and $\rho\mapsto I_{\rho}$ give 
an inclusion reversing one-to-one correspondence 
between the set of all gauge-invariant ideals 
and the set of all admissible pairs.
\end{theorem}

\begin{proof}
By Proposition~\ref{XIX} and Proposition~\ref{IXI}, 
these two maps are bijections and inverses of each others. 

For two gauge-invariant ideals $I_1$ and $I_2$, 
$I_1\subset I_2$ if and only if $I_1=I_1\cap I_2$. 
This is equivalent to $\rho_{I_1}=\rho_{I_1\cap I_2}$. 
Since $\rho_{I_1\cap I_2}=\rho_{I_1}\cup \rho_{I_2}$ 
by Proposition~\ref{cap}, 
the condition $\rho_{I_1}=\rho_{I_1\cap I_2}$ 
is the same as $\rho_{I_1}\supset \rho_{I_2}$. 
Thus $I_1\subset I_2$ if and only if 
$\rho_{I_1}\supset \rho_{I_2}$ for 
two gauge-invariant ideals $I_1,I_2$. 
\end{proof}

\begin{remark}
Admissible pairs introduced in \cite[Definition 8.16]{MT} 
are complements of our admissible pairs, 
and \cite[Theorem 8.22]{MT} includes the theorem above. 
These theorems are generalized in \cite[Theorem 8.6]{Ka3} 
for $C^*$-algebras arising from general $C^*$-correspondences. 
\end{remark}

A discrete graph $E$ is said to be {\em row-finite} 
if $E^0=\s{E}{0}{fin}$. 
For a row-finite discrete graph, 
the set of all gauge-invariant ideals 
is parameterized by invariant sets 
(see, for example, \cite[Theorem 4.1]{BPRS}). 
For a general topological graph $E$, 
the condition $E^0=\s{E}{0}{fin}$ 
does not suffice. 

\begin{example}\label{AbsVal}
Define a topological graph $E=(E^0,E^1,d,r)$ 
by $E^0=E^1=\R$, $d=\id$ and $r(x)=|x|$. 
Then we have $E^0=\s{E}{0}{fin}$. 
We have $\s{E}{0}{sce}=(-\infty,0)$ 
and $\s{E}{0}{sg}=\overline{\s{E}{0}{sce}}=(-\infty,0]$. 
Take $X^0=\{0\}$ which is a closed invariant set. 
We have $X^1=\{0\}$ and a topological graph $X=(X^0,X^1,d|_{X^1},r|_{X^1})$ 
has one vertex and one edge which is a loop on the vertex. 
We have $\s{X}{0}{rg}=\{0\}$ and $\s{X}{0}{sg}=\emptyset$. 
However $0\in \s{E}{0}{sg}$. 
Hence $\s{X}{0}{sg}\neq \s{E}{0}{sg}\cap X^0$. 
Thus there exists a topological graph $E$ with $E^0=\s{E}{0}{fin}$ 
such that invariant subsets are not sufficient 
to distinguish gauge-invariant ideals of $\cO(E)$. 
\end{example}

The correct definition of 
row-finiteness for topological graphs 
may be the following. 

\begin{definition}
A topological graph $E=(E^0,E^1,d,r)$ 
is said to be {\em row-finite} 
if $r(E^1)=\s{E}{0}{rg}$. 
\end{definition}

It is clear that this definition of row-finiteness 
generalizes the one for discrete graphs. 
It is also clear that topological graphs defined from 
dynamical systems are always row-finite. 
Since the topological graph in Example~\ref{AbsVal} 
is not row-finite, the condition $E^0=\s{E}{0}{fin}$ 
does not imply the row-finiteness in general. 
Conversely, row-finite topological graphs need not 
satisfy the condition $E^0=\s{E}{0}{fin}$. 

\begin{example}
Let us define a topological graph $E=(E^0,E^1,d,r)$ 
by $E^0=[0,1]$, $E^1=(0,1)$ and $d,r\colon E^1\to E^0$ are 
the natural embedding. 
We have $\s{E}{0}{rg}=\s{E}{0}{fin}=r(E^1)=(0,1)$ 
and $\s{E}{0}{sg}=\s{E}{0}{inf}=\{0,1\}$. 
Hence $E$ is row-finite, but $E^0\neq\s{E}{0}{fin}$. 
\end{example}

Note that $E^0=\s{E}{0}{fin}$ if and only if $r$ is proper, 
and these equivalent conditions imply that 
the image $r(E^1)$ of $r$ is closed. 
Conversely, for a row-finite topological graph $E$, 
the condition that $r(E^1)$ is closed 
implies $E^0=\s{E}{0}{fin}$. 

\begin{proposition}\label{row-finite}
Let $E$ be a row-finite topological graph. 
For a closed invariant set $X^0\subset E^0$, 
we have $\s{X}{0}{sg}=X^0\cap \s{E}{0}{sg}$ 
and the topological graph 
$X=(X^0,X^1,d_X,r_X)$ is also row-finite. 
\end{proposition}

\begin{proof}
Note that a topological graph $E$ is row-finite 
if and only if $\s{E}{0}{sg}=E^0\setminus r(E^1)$. 
Take a row-finite topological graph $E$, 
and a closed invariant set $X^0\subset E^0$. 
By Proposition~\ref{condforNI}, 
we have 
$$\s{X}{0}{sg}\subset X^0\cap \s{E}{0}{sg}=X^0\setminus r(E^1)
\subset X^0\setminus r_X(X^1)\subset X^0\setminus \s{X}{0}{rg}
=\s{X}{0}{sg}.$$
Hence we have $\s{X}{0}{sg}=X^0\cap \s{E}{0}{sg}$, 
and the topological graph $X$ is row-finite 
because $\s{X}{0}{sg}=X^0\setminus r_X(X^1)$. 
\end{proof}

\begin{corollary}
For a row-finite topological graph $E$, 
the map $I\mapsto X_I^0$ is a bijection 
from the set of all gauge-invariant ideals 
to the set of all closed invariant subsets of $E^0$. 
\end{corollary}

\begin{proof}
This follows from Theorem~\ref{gaugeideals} 
and Proposition~\ref{row-finite}. 
\end{proof}

\section{Orbits and maximal heads}\label{SecOrbit}

In this section, 
we generalize a notion of orbit spaces 
from ordinary dynamical systems 
to topological graphs, and study them. 

\begin{definition}
We define a {\em positive orbit space} $\Orb^+(v)$ of $v\in E^0$ by 
$$\Orb^+(v)=\{r^n(e)\in E^0\mid e\in (d^n)^{-1}(v), n\in\N\}.$$
\end{definition}

It is easy to see the following. 

\begin{proposition}\label{PosInvOrb}
A subset $X^0$ of $E^0$ is positively invariant 
if and only if $\Orb^+(v)\subset X^0$ for every $v\in X^0$. 
\end{proposition}

We examine for which $v\in E^0$, 
$\overline{\Orb^+(v)}$ becomes an invariant subset of $E^0$. 

\begin{lemma}\label{CloInv}
If a subset $X$ of $E^0$ is positively invariant or negatively invariant, 
then so is the closure $\overline{X}$. 
Hence $\overline{X}$ is invariant 
for an invariant subset $X\subset E^0$. 
\end{lemma}

\begin{proof}
Let $X$ be a positively invariant subset of $E^0$. 
Take $e\in E^1$ with $d(e)\in \overline{X}$. 
We can find a net $\{v_{\lambda}\}_{\lambda\in\Lambda}\subset X$ 
converging to $d(e)$.
Since $d$ is locally homeomorphic, 
we can find a net $\{e_{\lambda}\}\subset E^1$ 
such that $\lim e_{\lambda}=e$ and $d(e_{\lambda})=v_{\lambda}$ 
for sufficiently large $\lambda$. 
Since $X$ is positively invariant and 
$d(e_{\lambda})=v_{\lambda}\in X$ eventually, 
we have $r(e_{\lambda})\in X$ eventually. 
Since $r(e)=\lim r(e_{\lambda})$, 
we get $r(e)\in \overline{X}$.
Thus $\overline{X}$ is positively invariant.

Let $X$ be a negatively invariant subset of $E^0$. 
Take $v\in \s{E}{0}{rg}\cap\overline{X}$. 
Take a compact neighborhood $V$ of $v$ with $V\subset \s{E}{0}{rg}$. 
We can find a net $\{v_\lambda\}$ in $V$ 
such that $v_\lambda\in X$ and $\lim v_\lambda=v$. 
Since $v_\lambda\in X\cap V\subset X\cap \s{E}{0}{rg}$
and $X$ is negatively invariant, 
for each $\lambda$ there exists $e_\lambda\in E^1$ 
such that $r(e_\lambda)=v_\lambda$ and $d(e_\lambda)\in X$. 
Since $e_\lambda\in r^{-1}(V)$ for every $\lambda$ 
and $r^{-1}(V)\subset E^1$ is compact, 
we get a subnet $\{e_{\lambda_i}\}$ of $\{e_\lambda\}$ 
which converges to some $e\in r^{-1}(V)$. 
We have $r(e)=\lim r(e_{\lambda_i})=\lim v_{\lambda_i}=v$ 
and $d(e)=\lim d(e_{\lambda_i})\in\overline{X}$. 
Hence $\overline{X}$ is negatively invariant. 
The latter statement is clear by the former. 
\end{proof}

Recall that an element $e\in E^n$ for $n\geq 1$ 
is called a {\em loop} if $d^n(e)=r^n(e)$, 
and the vertex $d^n(e)=r^n(e)$ 
is called the {\em base point} of the loop $e$. 

\begin{proposition}\label{X=Orb}
For any $v\in E^0$, 
the closed set $\overline{\Orb^+(v)}$ 
is positively invariant. 
The set $\overline{\Orb^+(v)}$ is negatively invariant 
whenever at least one of the following three conditions is satisfied:
\benu
\item $v\in\s{E}{0}{sg}$. 
\item $v$ is not isolated in $\Orb^+(v)$.
\item $v\in E^0$ is a base point of a loop. 
\eenu
\end{proposition}

\begin{proof}
Set $X^0=\overline{\Orb^+(v)}$. 
Since $\Orb^+(v)$ is clearly positively invariant, 
$X^0$ is positively invariant by Lemma~\ref{CloInv}. 
The set $X^0$ is negatively invariant 
if and only if $\s{X}{0}{sce}\cap\s{E}{0}{rg}=\emptyset$ 
by Proposition~\ref{condforNI}. 
For $v'\in \Orb^+(v)\setminus\{v\}$, 
there exists $e\in E^1$ such that $r(e)=v'$ and $d(e)\in\Orb^+(v)$. 
Hence we get $\Orb^+(v)\setminus\{v\}\subset r(X^1)$. 
Therefore we have $\s{X}{0}{sce}=\{v\}$ 
when $v$ is isolated in $\Orb^+(v)$ and $v\notin r(X^1)$, 
and we have $\s{X}{0}{sce}=\emptyset$ otherwise. 
Hence $\s{X}{0}{sce}\cap\s{E}{0}{rg}\neq\emptyset$ 
if and only if 
all of the following three conditions hold; 
\benu
\item $v\in \s{E}{0}{rg}$, 
\item $v$ is isolated in $\Orb^+(v)$, 
\item $v\notin r(X^1)$. 
\eenu
The proof is completed 
if we show that $v\notin r(X^1)$ is 
equivalent to saying that $v\in E^0$ is a base point of no loops, 
under the assumption that $v$ is isolated in $\Orb^+(v)$. 
If $v\in E^0$ is a base point of a loop, 
we can find $e\in E^1$ such that $r(e)=v, d(e)\in\Orb^+(v)$. 
Hence $v\in r(X^1)$. 
Conversely assume that $v\in r(X^1)$. 
Take $e\in X^1$ with $r(e)=v$. 
We can find a net $\{v_\lambda\}\subset \Orb^+(v)$ 
converging to $d(e)\in X^0$. 
Since $d$ is locally homeomorphic, 
we can find a net $\{e_{\lambda}\}\subset E^1$ 
such that $\lim e_{\lambda}=e$ and $d(e_{\lambda})=v_{\lambda}$ 
for sufficiently large $\lambda$. 
We have $v=\lim r(e_{\lambda})$ and $r(e_{\lambda})\in\Orb^+(v)$. 
Since we assume that $\{v\}$ is isolated in $\Orb^+(v)$, 
we have $v=r(e_{\lambda})$ 
for sufficiently large $\lambda$. 
This means that $v\in E^0$ is a base point of a loop. 
The proof is completed. 
\end{proof}

\begin{remark}
It is not difficult to see that 
the positive orbit space $\Orb^+(v)$ of $v\in E^0$ 
is negatively invariant if and only if 
either the condition (i) or (iii) in Proposition~\ref{X=Orb} 
is satisfied (see Proposition~\ref{Orbit}). 
\end{remark}

Next we define negative orbit spaces of vertices. 
An {\em infinite path} is a sequence $e=(e_1,e_2,\ldots,e_n,\ldots)$ 
with $e_k\in E^1$ and $d(e_k)=r(e_{k+1})$ for each $k=1,2,\ldots$. 
The set of all infinite paths is denoted by $E^\infty$. 
For an infinite path $e=(e_1,e_2,\ldots,e_n,\ldots)\in E^\infty$, 
we define its range $r(e)\in E^0$ to be $r(e_1)$. 

\begin{definition}
For $n\in\N\cup\{\infty\}$, 
a path $e\in E^n$ is called a {\em negative orbit} of $v\in E^0$ 
if $r^n(e)=v$ and $d^n(e)\in\s{E}{0}{sg}$ when $n<\infty$. 
\end{definition}

Lemma~\ref{RgPr} ensures that 
each $v\in E^0$ has at least one negative orbit, 
but $v$ may have many negative orbits in general. 
If $v\in\s{E}{0}{sg}$, 
then $v$ itself is a negative orbit of $v$. 
We denote a negative orbit of $v$ by $e=(e_1,e_2,\ldots,e_n)\in E^n$ 
for $n\in\N\cup\{\infty\}$. 
When $n=\infty$ this expression is understood 
as $e=(e_1,e_2,\ldots)\in E^\infty$, 
and when $n=0$ this means that $e=v\in E^0$. 

\begin{definition}
For each negative orbit $e=(e_1,e_2,\ldots,e_n)\in E^n$ of $v\in E^0$, 
a {\em negative orbit space} $\Orb^-(v,e)$ 
is defined by 
$$\Orb^-(v,e)=\{v,d(e_1),d(e_2),\ldots,d(e_n)\}\subset E^0.$$ 
\end{definition}

\begin{proposition}\label{NegInvOrb}
A subset $X$ of $E^0$ is negatively invariant 
if and only if for each $v\in X$, 
there exists a negative orbit $e$ of $v$ 
such that $\Orb^-(v,e)\subset X$. 
\end{proposition}

\begin{proof}
Let $X$ be a negatively invariant subset of $E^0$. 
Take $v\in X$. 
When $v\in \s{E}{0}{sg}$, $v\in E^0$ is a negative orbit of $v$ 
satisfying $\Orb^-(v,v)=\{v\}\subset X$. 
When $v\in \s{E}{0}{rg}$, we can find $e_1\in E^1$ 
such that $r(e_1)=v$ and $d(e_1)\in X$ 
because $X$ is negatively invariant. 
If $d(e_1)\in \s{E}{0}{sg}$, 
then $e_1\in E^1$ is a negative orbit of $v$ 
satisfying $\Orb^-(v,e_1)=\{v,d(e_1)\}\subset X$. 
Otherwise we get $e_2\in E^1$ 
such that $r(e_2)=d(e_1)$ and $d(e_2)\in X$. 
Again if $d(e_2)\in \s{E}{0}{sg}$, 
then $e=(e_1,e_2)\in E^2$ is a negative orbit of $v$ 
satisfying $\Orb^-(v,e)=\{v,d(e_1),d(e_2)\}\subset X$. 
In such a manner, 
we get $e_k\in E^1$ for $k=1,2,\ldots,n$ 
with $r(e_1)=v$ and $d(e_k)=r(e_{k+1})\in \s{E}{0}{rg}\cap X$ 
for $1\leq k\leq n-1$ 
unless we get $d(e_n)\in \s{E}{0}{sg}\cap X$. 
If we have $d(e_n)\in \s{E}{0}{sg}\cap X$ for some $n\in\N$, 
then $e=(e_1,\ldots,e_n)\in E^n$ is a negative orbit of $v$ 
with $\Orb^-(v,e)\subset X$. 
Otherwise we get $e_k\in E^1$ 
with $r(e_1)=v$ and $d(e_k)=r(e_{k+1})\in X$ 
for $k=1$, $2$, \ldots . 
In this case, 
the infinite path $e=(e_1,e_2,\ldots)\in E^\infty$ 
is a negative orbit of $v$ with $\Orb^-(v,e)\subset X$. 

Conversely assume that for each $v\in X$, 
there exists a negative orbit $e$ of $v$ 
such that $\Orb^-(v,e)\subset X$. 
Take $v\in \s{E}{0}{rg}\cap X$. 
The vertex $v\in X$ has a negative orbit $e\in E^n$ 
with $\Orb^-(v,e)\subset X$. 
Since $v\in \s{E}{0}{rg}$, 
we have $n\geq 1$. 
Hence there exists $e_1\in E^1$ which satisfies that $r(e_1)=v$ 
and $d(e_1)\in \Orb^-(v,e)\subset X$. 
Therefore $X$ is negatively invariant. 
\end{proof}

\begin{definition}
We define the {\em orbit space} $\Orb(v,e)$ of $v\in E^0$ 
with respect to a negative orbit $e$ of $v$ by 
$$\Orb(v,e)=\bigcup_{v'\in\Orb^-(v,e)}\Orb^+(v').$$
\end{definition}

For a negative orbit $e\in E^n$ of $v\in E^0$ with $0\leq n<\infty$, 
we have $\Orb(v,e)=\Orb^+(d(e))$. 
When a negative orbit $e\in E^\infty$ of $v$ 
is defined by $e=(e',e',\ldots)$ 
for a loop $e'\in E^*$ whose base point is $v$, 
the orbit space $\Orb(v,e)$ coincides with $\Orb^+(v)$. 

\begin{proposition}\label{Orbit}
An orbit space $\Orb(v,e)$ is an invariant set 
for every $v\in E^0$ and every negative orbit $e$ of $v$. 
\end{proposition}

\begin{proof}
Take $v\in E^0$ and a negative orbit 
$e=(e_1,e_2,\ldots,e_n)$ of $v$ ($n\in\N\cup\{\infty\}$). 
Since $\Orb(v,e)$ is a union of positive orbit spaces, 
it is positively invariant. 
When $n=\infty$, for every $v'\in \Orb(v,e)$ 
there exists $e'\in r^{-1}(v')\subset E^1$ with $d(e')\in \Orb(v,e)$. 
Hence $\Orb(v,e)$ is negatively invariant. 
When $n<\infty$, for every $v'\in \Orb(v,e)$ except $d(e)$ 
there exists $e'\in r^{-1}(v')\subset E^1$ with $d(e')\in \Orb(v,e)$. 
Noting that $d(e)\in \s{E}{0}{sg}$, 
we see that $\Orb(v,e)$ is negatively invariant. 
Hence $\Orb(v,e)$ is invariant. 
\end{proof}

\begin{proposition}\label{InvOrb}
A subset $X$ of $E^0$ is invariant if and only if 
for each $v\in X$ there exists a negative orbit $e$ of $v$ 
such that $\Orb(v,e)\subset X$. 
\end{proposition}

\begin{proof}
Combine Proposition~\ref{PosInvOrb}
and Proposition~\ref{NegInvOrb}. 
\end{proof}

\begin{definition}
A subset $X^0$ of $E^0$ is called a {\em maximal head} 
if $X^0$ is a non-empty closed invariant set satisfying that 
for any $v_1,v_2\in X^0$ and 
any neighborhoods $V_1$, $V_2$ of $v_1,v_2$ respectively,
there exists $v\in X^0$ 
with $\Orb^+(v)\cap V_1\neq\emptyset$ 
and $\Orb^+(v)\cap V_2\neq\emptyset$. 
\end{definition}

Equivalently, 
a non-empty closed invariant set $X^0$ is a maximal head 
if and only if for all $v,v'\in X^0$ 
there exist nets $\{e_\lambda\},\{e_\lambda'\}\subset X^*$ 
of finite paths 
such that $d(e_\lambda)=d(e_\lambda')\in X^0$ for all $\lambda$, 
$\lim r(e_\lambda)=v$ and $\lim r(e_\lambda')=v'$. 

\begin{proposition}\label{Orbit=MH}
For every $v\in E^0$ and every negative orbit $e$ of $v$, 
the closed set $\overline{\Orb(v,e)}$ is a maximal head. 
\end{proposition}

\begin{proof}
Take $v\in E^0$ and a negative orbit 
$e=(e_1,\ldots,e_n)\in E^n$ of $v$ where $n\in\N\cup\{\infty\}$.
By Lemma~\ref{CloInv} and Proposition~\ref{Orbit}, 
the set $\overline{\Orb(v,e)}$ is invariant. 
Set $v_0=v$ and $v_k=d(e_k)$ for $k\in\{1,2,\ldots,n\}$. 
Then we have $\Orb^-(v,e)=\{v_k\}_{k=0}^n$. 
Take $w_1,w_2\in \overline{\Orb(v,e)}$ and 
neighborhoods $V_1$, $V_2$ of $w_1,w_2$ respectively. 
There exist $w_1'\in V_1\cap \Orb(v,e)$ 
and $w_2'\in V_2\cap \Orb(v,e)$. 
Since $\Orb(v,e)=\bigcup_{k=0}^n\Orb^+(v_k)$ 
and $\Orb^+(v_k)\subset \Orb^+(v_{k+1})$ for $0\leq k\leq n-1$, 
we can find an integer $m$ with $0\leq m\leq n$ 
such that $w_1',w_2'\in \Orb^+(v_m)$. 
Thus $v_m\in \overline{\Orb(v,e)}$ 
satisfies that $\Orb^+(v_m)\cap V_1\neq\emptyset$ and 
$\Orb^+(v_m)\cap V_2\neq\emptyset$. 
Therefore $\overline{\Orb(v,e)}$ is a maximal head. 
\end{proof}

The converse of Proposition~\ref{Orbit=MH} is true 
when $E^0$ is second countable. 

\begin{proposition}\label{Orbit=MH2}
If $E^0$ is second countable, 
every maximal head is of the form $\overline{\Orb(v,e)}$ 
for some negative orbit $e$ of $v\in E^0$. 
\end{proposition}

\begin{proof}
Take a maximal head $X^0$ of $E$. 
Let $\{V_k\}_{k=0}^\infty$ be a countable open basis of $X^0$. 
Take a non-empty open subset $V_0'$ of $X^0$ arbitrarily. 
Since $X^0$ is a maximal head, 
we can find $e_0\in X^{n_0}$ and $e_0'\in X^{n_0'}$ for $n_0,n_0'\in\N$ 
such that $d(e_0)=d(e_0')\in X^0$, 
$r(e_0)\in V_0$ and $r(e_0')\in V_0'$. 
Choose compact neighborhoods $U_0\subset X^{n_0}$ and 
$U_0'\subset X^{n_0'}$ of $e_0\in X^{n_0}$ and $e_0'\in X^{n_0'}$ 
such that $U_0\subset r^{-1}(V_0)$ and $U_0'\subset r^{-1}(V_0')$. 
Since $d(U_0)\cap d(U_0')$ is a neighborhood of $d(e_0)=d(e_0')$, 
we can find a non-empty open subset $V_1'$ of $X^0$ such that 
$V_1'\subset d(U_0)\cap d(U_0')$. 
Inductively we can find non-empty compact sets $U_k\subset X^{n_k}$ and 
$U_k'\subset X^{n_k'}$ for some $n_k,n_k'\in\N$ 
such that $r(U_k)\subset V_k$ and 
$r(U_k')\subset d(U_{k-1})\cap d(U_{k-1}')$ 
for $k\in\N$. 
We denote by $U'$ the compact set 
$U'_0\times U'_1\times\cdots\times U'_k\times\cdots$. 
For each $k\in\N$, 
we define a closed subset $F_k$ 
of the compact set $U'$ by 
$$F_k=\big\{(e_0,\ldots,e_k,\ldots)\in 
U'\ \big|\ d(e_{j-1})=r(e_j) 
\mbox{ for }j=1,2,\ldots,k\big\}$$
It is easy to see that $\{F_k\}_{k\in\N}$ 
is a decreasing sequence of non-empty closed subsets. 
By the compactness of $U'$, 
we can find an element $(e_0,\ldots,e_k,\ldots)\in\bigcap_{k\in\N}F_k$. 
Then $e=(e_0,\ldots,e_k,\ldots)$ is a negative orbit 
of $v=r(e)\in X^0$ 
(note that the length of $e$ is finite when $n_k=0$ eventually). 
We will prove that $X^0=\overline{\Orb(v,e)}$. 
Since $\Orb^-(v,e)\subset X^0$, 
we have $\overline{\Orb(v,e)}\subset X^0$. 
For each $k\in\N$, 
we can find $e_k'\in E^*$ such that $d(e_k')=d(e_k)\in \Orb^-(v,e)$ 
and $r(e_k')\in V_k$. 
Hence we have $\Orb(v,e)\cap V_k\neq\emptyset$ for every $k\in\N$. 
Therefore we have $X^0=\overline{\Orb(v,e)}$. 
\end{proof}

By Proposition~\ref{Orbit=MH2}, 
for discrete graphs $E=(E^0,E^1,d,r)$ with countable $E^0$, 
every maximal head is of the form $\Orb(v,e)$ 
for some $v\in E^0$ and some negative orbit $e$ of $v$. 
This is no longer true for discrete graphs $E=(E^0,E^1,d,r)$ 
such that $E^0$ is uncountable, as the following example shows. 

\begin{example}\label{NiceExample}
Let $X$ be an uncountable set, 
and $E^0$ be the set of finite subsets of $X$ 
with discrete topology. 
Let 
$$E^1=\big\{(x;v)\ \big|\ 
\text{$v\in E^0$ and $x\in v$}\big\}.$$
We define $d,r\colon E^1\to E^0$ by $d((x;v))=v$ 
and $r((x;v))=v\setminus \{x\}$ for $(x;v)\in E^1$. 
For $v\in E^0$, 
we have $\Orb^+(v)=\{w\in E^0\mid w\subset v\}$. 
Hence for $v_1,v_2\in E^0$, $v_0=v_1\cup v_2\in E^0$ 
satisfies $v_1,v_2\in \Orb^+(v_0)$. 
This shows that $E^0$ is a maximal head. 
For a negative orbit $e=(e_1,e_2,\ldots,e_n)\in E^n$ 
of $v_0\in E^0$ where $n\in\N\cup\{\infty\}$, 
$\Orb(v_0,e)=\{v\in E^0\mid v\subset \Omega\}$ 
where $\Omega=\bigcup_{k=1}^n d(e_k)$. 
Since $\Omega$ is countable, $\Orb(v_0,e)\neq E^0$ 
\end{example}

\begin{remark}
In Example~\ref{ExErgode}, 
we see another example of topological graphs $E$, 
which comes from a dynamical system, 
such that $E^0$ is a maximal head, but is not 
in the form $\overline{\Orb(v,e)}$ 
for a negative orbit $e$ of $v\in E^0$. 
\end{remark}

\section{Hereditary and saturated sets}\label{SecHS}

In this section, we study the complement of invariant subsets, 
and get a characterization of maximal heads. 

\begin{definition}
A subset $V$ of $E^0$ is said to be {\em hereditary} 
if $V$ satisfies $d(r^{-1}(V))\subset V$, 
and said to be {\em saturated} 
if $v\in\s{E}{0}{rg}$ satisfying $d(r^{-1}(v))\subset V$ is in $V$. 
\end{definition}

It is easy to see the following. 

\begin{proposition}\label{complement}
For a subset $X$ of $E^0$, 
$X$ is positively invariant if and only if 
the complement $V=E^0\setminus X$ of $X$ is hereditary, 
and $X$ is negatively invariant if and only if 
its complement $V$ is saturated. 
\end{proposition}

\begin{definition}\label{DefHS}
Let us take a subset $V$ of $E^0$. 
We define a subset $H(V)$ of $E^0$ by 
$$H(V)=\bigcup_{n=0}^\infty d^n\big((r^n)^{-1}(V)\big),$$ 
and $S(V)\subset E^0$ by $S(V)=\bigcup_{k=0}^\infty V_k$, 
where $V_k\subset E^0$ is defined by $V_0=V$ and 
$$V_{k+1}=V_k\cup \{v\in\s{E}{0}{rg}\mid d(r^{-1}(v))\subset V_k\}.$$
\end{definition}

\begin{proposition}\label{hersat1}
For a subset $V$ of $E^0$, 
$H(V)$ is the smallest hereditary set containing $V$. 
If $V$ is open, then so is $H(V)$. 
\end{proposition}

\begin{proof}
Clear by the definition of $H(V)$. 
\end{proof}

\begin{proposition}\label{hersat2}
Let $V$ be an open subset of $E^0$. 
Then $S(V)$ is open and 
the smallest saturated set containing $V$. 
If $V$ is hereditary, then so is $S(V)$. 
\end{proposition}

\begin{proof}
Let $\{V_k\}$ be subsets of $E^0$ as in Definition~\ref{DefHS} 
with $S(V)=\bigcup_{k=0}^\infty V_k$. 
If $V$ is open, 
we can prove that 
$V_k$ is open for $k\in\N$ inductively 
by Lemma~\ref{RgPr}. 
Hence $S(V)$ is open. 
By the definition of $S(V)$, 
it is clear that any saturated set containing $V$ 
contains $S(V)$. 
We will show that $S(V)$ is saturated. 
Take $v\in\s{E}{0}{rg}$ with $d(r^{-1}(v))\subset S(V)$. 
Since $d(r^{-1}(v))$ is compact by Lemma~\ref{RgPr}, 
there exists $k\in\N$ such that $d(r^{-1}(v))\subset V_k$. 
Hence we have $v\in V_{k+1}\subset S(V)$. 
Thus $S(V)$ is a saturated set containing $V$. 

Now suppose that $V$ is hereditary. 
Since we have $d(r^{-1}(V))\subset V$, 
we can prove that $d(r^{-1}(V_{k+1}))\subset V_k$ for $k\in\N$ 
inductively. 
Hence we have $d(r^{-1}(S(V)))\subset S(V)$. 
We are done. 
\end{proof}

By Proposition~\ref{hersat1} and Proposition~\ref{hersat2}, 
we get the following. 

\begin{proposition}\label{hersat4}
For an open set $V$, 
the closed set $X=E^0\setminus S(H(V))$ is the largest invariant 
set which does not intersect $V$. 
\end{proposition}

We finish this section by giving a characterization 
of maximal heads. 

\begin{lemma}\label{hersat5}
Let $V,X$ be two subsets of $E^0$ 
such that $X$ is negatively invariant. 
Then $S(V)\cap X\neq\emptyset$ if and only if $V\cap X\neq\emptyset$. 
\end{lemma}

\begin{proof}
Since $V\subset S(V)$, we have $S(V)\cap X\neq\emptyset$ 
if $V\cap X\neq\emptyset$. 
Conversely suppose $V\cap X=\emptyset$. 
Then $V$ is contained in the complement $E^0\setminus X$ 
of $X$ which is saturated by Proposition~\ref{complement}. 
Hence $S(V)\subset E^0\setminus X$. 
This shows $S(V)\cap X=\emptyset$. 
We are done. 
\end{proof}

\begin{lemma}\label{hersat7}
Let $V_1,V_2$ be two subsets of $E^0$, 
and $X^0$ be an invariant subset. 
Then $X^0\cap S(H(V_1))\cap S(H(V_2))\neq\emptyset$ 
if and only if $X^0\cap H(V_1)\cap H(V_2)\neq\emptyset$. 
\end{lemma}

\begin{proof}
By Proposition~\ref{hersat2}, 
$S(H(V_2))$ is hereditary. 
Hence the intersection $X^0\cap S(H(V_2))$ is negatively invariant. 
Thus Lemma~\ref{hersat5} shows 
$X^0\cap S(H(V_1))\cap S(H(V_2))\neq\emptyset$ 
if and only if $X^0\cap H(V_1)\cap S(H(V_2))\neq\emptyset$. 
By the same reason, 
$X^0\cap H(V_1)\cap S(H(V_2))\neq\emptyset$ 
if and only if $X^0\cap H(V_1)\cap H(V_2)\neq\emptyset$. 
We are done.
\end{proof}

\begin{proposition}\label{maxhead}
An invariant set $X^0$ of $E^0$ is a maximal head 
if and only if 
$X^0_1\cup X^0_2\supset X^0$ 
implies either $X^0_1\supset X^0$ or $X^0_2\supset X^0$ 
for two invariant sets $X^0_1,X^0_2$. 
\end{proposition}

\begin{proof}
Let $X^0$ be a maximal head.
Take invariant sets $X^0_1$ and $X^0_2$ satisfying 
$X^0_1\not\supset X^0$ and $X^0_2\not\supset X^0$, 
and we will show that $X^0_1\cup X^0_2\not\supset X^0$.
There exist $v_1,v_2\in X^0$ 
with $v_1\notin X^0_1$ and $v_2\notin X^0_2$.
Since $X^0$ is a maximal head, 
there exist $v\in X^0$ and $e_1,e_2\in E^*$ with 
$d(e_1)=d(e_2)=v$, $r(e_1)\notin X_1^0$ and $r(e_2)\notin X_2^0$. 
Since $X^0_1$ and $X^0_2$ are positively invariant, we have $v\notin X^0_1$
and $v\notin X^0_2$.
Therefore, $X^0_1\cup X^0_2\not\supset X^0$. 

Conversely assume that a closed invariant set $X^0$ 
satisfies that $X^0_1\cup X^0_2\supset X^0$ 
implies either $X^0_1\supset X^0$ or $X^0_2\supset X^0$ 
for two invariant sets $X^0_1,X^0_2$. 
Take $v_1,v_2\in X^0$ and 
neighborhoods $V_1$, $V_2$ of $v_1,v_2$ respectively. 
For $i=1,2$, 
set $X^0_i=E^0\setminus S(H(V_i))$ 
which are invariant by Proposition~\ref{hersat4}.
We have $X^0_1\not\supset X^0$ and $X^0_2\not\supset X^0$. 
By the assumption, 
we have $X^0_1\cup X^0_2\not\supset X^0$.
This implies $X^0\cap S(H(V_1))\cap S(H(V_2))\neq\emptyset$. 
By Lemma~\ref{hersat7}, 
we have $X^0\cap H(V_1)\cap H(V_2)\neq\emptyset$. 
Then $v\in X^0\cap H(V_1)\cap H(V_2)$ 
satisfies $\Orb^+(v)\cap V_1\neq\emptyset$ 
and $\Orb^+(v)\cap V_2\neq\emptyset$. 
This shows that $X^0$ is a maximal head. 
\end{proof}

\section{Ideals and topological freeness}\label{SecTopFree}

\begin{proposition}\label{IrV1}
For an open subset $V$ of $E^0$, 
the ideal $I$ of $\cO(E)$ 
generated by $t^0(C_0(V))\subset\cO(E)$ 
is gauge-invariant and satisfies 
$\rho_I=\big(E^0\setminus S(H(V)), \s{E}{0}{sg}\setminus S(H(V))\big)$. 
\end{proposition}

\begin{proof}
Since $t^0(C_0(V))$ is invariant under the gauge action, 
the ideal $I$ is gauge-invariant. 
Set $X^0=E^0\setminus S(H(V))$ which is a closed invariant set. 
We set $\rho=(X^0,\s{E}{0}{sg}\cap X^0)$ 
which is an admissible pair. 
We will show $\rho_I=\rho$. 
Since $V\subset S(H(V))=E^0\setminus X^0$, 
we have 
$$t^0(C_0(V))\subset t^0(C_0(E^0\setminus X^0))\subset I_{\rho}.$$ 
Thus $I\subset I_{\rho}$. 
By Lemma~\ref{inclusion} and Proposition~\ref{XIX}, 
we have $\rho_I\supset \rho_{I_{\rho}}=\rho$. 
Thus $X_I^0\supset X^0$ and $Z_I\supset \s{E}{0}{sg}\cap X^0$. 
Since $t^0(C_0(V))\subset I$, 
we have $X_I^0\subset E^0\setminus V$. 
Thus Proposition~\ref{hersat4} shows $X_I^0\subset E^0\setminus S(H(V))=X^0$. 
Therefore $X_I^0=X^0$. 
Hence we get $Z_I=\s{E}{0}{sg}\cap X^0$. 
Thus we have $\rho_I=\rho$. 
We are done. 
\end{proof}

\begin{remark}
For an open subset $V$ of $E^0$, 
the condition $d(r^{-1}(V))\subset V$ considered 
in \cite[Proposition 5.9]{Ka2} is nothing but 
the hereditariness, 
and the condition that 
each $v\in E^0\setminus V$ is regular 
and satisfies $d^n((r^n)^{-1}(v))\subset V$ for some $n\in\N$ 
is equivalent to $S(V)=E^0$. 
Thus Proposition~\ref{IrV1} shows that this condition 
is a necessary and sufficient condition for $A_V$ to be full, 
as predicted in \cite[Remark 5.10]{Ka2}. 
\end{remark}

Now we recall some arguments from \cite[Section 5]{Ka2}. 
Let $F^0$ be a hereditary open subset of $E^0$, 
and define $F^1=r^{-1}(F^0)$. 
The set $F^1$ is an open subset of $E^1$, 
and satisfies $d(F_1),r(F_1)\subset F^0$. 
Thus $F=(F^0,F^1,d|_{F^1},r|_{F^1})$ is a subgraph of $E$ 
in the sense of \cite[Definition 5.1]{Ka2}. 

\begin{proposition}\label{IrV2}
Under the notations above, 
the $C^*$-subalgebra $A$ of $\cO(E)$ 
generated by $t^0(C_0(F^0))$ and $t^1(C_d(F^1))$ 
is a hereditary and full subalgebra of 
the ideal $I$ generated by $t^0(C_0(F^0))$, 
and there exists a natural isomorphism between 
$A$ and $\cO(F)$ which commutes with the natural injections 
from $C_0(F^0)$ and $C_d(F^1)$ to $A$ and $\cO(F)$. 
\end{proposition}

\begin{proof}
For $\xi\in C_c(F^1)$, 
we can find $f\in C_0(F^0)$ 
such that $\pi_r(f)\xi=\xi$. 
Hence $t^1(\xi)=t^0(f)t^1(\xi)\in I$. 
This shows $t^1(C_d(F^1))\subset I$. 
Therefore $A$ is a full subalgebra of $I$. 
The rest of the statements follows from 
\cite[Propositions 5.5, 5.9]{Ka2}. 
\end{proof}

By the two propositions above, 
we get the following proposition. 

\begin{proposition}
Let $I$ be a gauge-invariant ideal of $\cO(E)$ 
with $Z_I=\s{E}{0}{sg}\cap X_I^0$. 
Take a hereditary open subset $F^0$ with 
$S(F^0)=E^0\setminus X_I^0$ 
(for example take $F^0=E^0\setminus X_I^0$). 
Set a subgraph $F=(F^0,F^1,d|_{F^1},r|_{F^1})$ of $E$ 
by $F^1=r^{-1}(F^0)$. 
Then $I$ is strongly Morita equivalent to $\cO(F)$. 
\end{proposition}

\begin{proof}
By Proposition~\ref{IrV1} and Theorem~\ref{gaugeideals}, 
$I$ coincides with the ideal generated by $t^0(C_0(F^0))$, 
which is strongly Morita equivalent to $\cO(F)$ 
by Proposition~\ref{IrV2}. 
\end{proof}

In \cite{Ka5}, 
we will see that 
for an arbitrary gauge-invariant ideal $I$ of $\cO(E)$, 
one can find a topological graph $F$ such that 
$I$ is strongly Morita equivalent to $\cO(F)$. 
We will also see that a gauge-invariant ideal $I$ itself 
can be expressed as $\cO(F')$ for some topological graph $F'$ 
which is less natural than the topological graph $F$ above. 
It is hopeless (or useless) to express an ideal of $\cO(E)$ 
which is not gauge-invariant 
as a $C^*$-algebra of some topological graph 
(see Example~\ref{easyexam}). 

In the rest of this section, 
we apply Proposition~\ref{IrV2} 
to prove that topological freeness 
is needed in the Cuntz-Krieger Uniqueness Theorem 
(\cite[Theorem 5.12]{Ka1}). 
We recall the definition of topological freeness, 
and the statement of the theorem. 

\begin{definition}
A loop $e=(e_1,\ldots,e_n)$ is said to be {\em simple} 
if $r(e_i)\neq r(e_j)$ for $i\neq j$, 
and said to be {\em without entrances} 
if $r^{-1}(r(e_k))=\{e_{k}\}$ for $k=1,\ldots,n$.
\end{definition}

It is easy to see that if $v\in E^0$ is a base point of a loop, 
then $v$ is a base point of a simple loop. 
It is also easy to see that 
if $v\in E^0$ is a base point of a loop without entrances, 
then there exists a unique simple loop 
whose base point is $v$, and this loop is also without entrances. 
(A loop $e$ without entrances is in the form $e=(e',\ldots,e')$ 
where $e'$ is a simple loop without entrances.)

\begin{definition}[{\cite[Definition 5.4]{Ka1}}]
A topological graph $E$ is said to be {\em topologically free} 
if the set of base points of loops without entrances 
has an empty interior. 
\end{definition}

The following theorem is called 
the Cuntz-Krieger Uniqueness Theorem. 

\begin{proposition}[{\cite[Theorem 5.12]{Ka1}}]\label{CKUT}
Let $E$ be a topologically free topological graph. 
Then the natural surjection $\cO(E)\to C^*(T)$ 
is an isomorphism for every injective Cuntz-Krieger $E$-pair $T$. 
\end{proposition}

We will prove its converse. 
Namely, when $E$ is not topologically free, 
we will find an injective Cuntz-Krieger $E$-pair $T$ 
such that the natural surjection $\cO(E)\to C^*(T)$ is not an isomorphism. 
To find such injective Cuntz-Krieger $E$-pair $T$, 
it suffices to find a non-zero ideal $I$ of $\cO(E)$ 
with $\bigcap_{z\in\T}\beta_z(I)=0$ by the following lemma. 

\begin{lemma}\label{IG=0}
For an ideal $I$ of the $C^*$-algebra $\cO(E)$, 
the following conditions are equivalent: 
\benu
\item $I$ is the kernel of the $*$-homomorphism $\cO(E)\to C^*(T)$ 
induced by an injective Cuntz-Krieger $E$-pair $T$. 
\item $I\cap t^0(C_0(E^0))=0$. 
\item $X_I^0=E^0$. 
\item $\rho_I=(E^0,\s{E}{0}{sg})$. 
\item $\bigcap_{z\in\T}\beta_z(I)=0$. 
\eenu
\end{lemma}

\begin{proof}
Clearly we have (i)$\iff$(ii)$\iff$(iii)$\iff$(iv). 
The equivalence (iv)$\iff$(v) follows from Proposition~\ref{IrI}. 
\end{proof}

We start with the following example. 

\begin{example}\label{easyexam}
Let $n$ be a positive integer. 
Let $E_n=(E_n^0,E_n^1,d_n,r_n)$ be a discrete graph 
such that $E_n^0=E_n^1=\Z/n\Z$, $d_n=\id_{\Z/n\Z}$, 
and $r_n(k)=k+1$. 
The graph $E_n$ consists of one loop without entrances. 
Thus every vertex of $E_n$ is a base point of a loop without entrances. 
Therefore $E_n$ is not topologically free. 

Let us denote by $\{1,2,\ldots,n\}$ the elements of $\Z/n\Z$. 
The matrix units of 
the $C^*$-algebra $\M_n$ of all $n\times n$ matrices 
are denoted by $\{u_{k,l}\}_{k,l\in \Z/n\Z}$, 
and the generating unitary of $C(\T)$ 
is denote by $w\in C(\T)$. 
We will prove that $\cO(E_n)\cong C(\T)\otimes \M_n$ 
(cf.\ the proof of Proposition~\ref{C(T,K)}). 
We define a $*$-homomorphism $T^0\colon C(E_n^0)\to C(\T)\otimes \M_n$ 
and a linear map $T^1\colon C_{d_n}(E_n^1)\to C(\T)\otimes \M_n$ by 
$$T^0(f)=\sum_{k\in \Z/n\Z}f(k)\otimes u_{k,k},\quad
T^1(\xi)=\sum_{k=1}^{n-1}\xi(k)\otimes u_{k+1,k}
+(\xi(n)w)\otimes u_{1,n},$$
for $f\in C(E_n^0)$ and $\xi\in C_{d_n}(E_n^1)$. 
It is routine to check that 
the pair $T=(T^0,T^1)$ is an injective Cuntz-Krieger $E_n$-pair 
which admits a gauge action 
and satisfies $C^*(T)=C(\T)\otimes \M_n$. 
Hence by \cite[Theorem 4.5]{Ka1}, 
we have $\cO(E_n)\cong C(\T)\otimes \M_n$. 
\end{example}

\begin{lemma}\label{non-inv-ideal1}
Let $n$ be a positive integer, 
and $E_n$ be the discrete graph in Example~\ref{easyexam}. 
Then there exists a non-zero ideal $I$ of $\cO(E_n)$ 
with $I\cap t^0(C(E_n^0))=0$. 
\end{lemma}

\begin{proof}
Take $z\in\T$. 
Then the ideal $I$ of $\cO(E)$ 
corresponding to the ideal $C_0(\T\setminus\{z\})\otimes \M_n$ 
of $C(\T)\otimes \M_n\cong\cO(E_n)$ 
satisfies that $I\neq 0$ and $I\cap t^0(C(E_n^0))=0$. 
\end{proof}

\begin{lemma}\label{non-inv-ideal2}
Let $E_n=(E_n^0,E_n^1,d_n,r_n)$ be 
the discrete graph in Example~\ref{easyexam}, 
and $X$ be a locally compact space. 
Let $F=(F^0,F^1,d_F,r_F)$ be a topological graph 
such that there exist 
homeomorphisms $m^0\colon F^0\to E_n^0\times X$ 
and $m^1\colon F^1\to E_n^1\times X$ 
with $(d_n\times\id_X)\circ m^1=m^0\circ d_F$ and 
$(r_n\times\id_X)\circ m^1=m^0\circ r_F$. 
Then there exists a non-zero ideal $I'$ of $\cO(F)$ 
with $I'\cap t_F^0(C_0(F^0))=0$ 
where $t_F=(t_F^0,t_F^1)$ is 
the universal Cuntz-Krieger $F$-pair on $\cO(F)$. 
\end{lemma}

\begin{proof}
It is easy to see that 
the two maps $m^0,m^1$ induce isomorphisms 
$C_0(F^0)\cong C(E_n^0)\otimes C_0(X)$ and 
$C_{d_F}(F^1)\cong C_{d_n}(E_n^1)\otimes C_0(X)$, 
and these isomorphisms induce 
an isomorphism $\cO(F)\cong \cO(E_n)\otimes C_0(X)$ 
(see \cite[Proposition 7.7]{Ka2}). 
By Lemma~\ref{non-inv-ideal1}, 
there exists a non-zero ideal $I$ of $\cO(E_n)$ 
such that $I\cap t^0(C(E_n^0))=0$. 
Then the ideal $I'$ of $\cO(F)$ 
corresponding to the ideal $I\otimes C_0(X)$ 
of $\cO(E_n)\otimes C_0(X)$ 
satisfies the desired conditions. 
\end{proof}

\begin{proposition}\label{Baire}
A topological graph $E$ is not topologically free 
if and only if there exist a non-empty open subset $V$ of $E^0$ 
and a positive integer $n$ such that 
all vertices in $V$ are base points of simple loops in $E^n$ 
without entrances, 
and that $\sigma=r\circ (d|_{r^{-1}(V)})^{-1}$ is 
a well-defined continuous map from $V$ to $V$ 
with $\sigma^n=\id_{V}$. 
\end{proposition}

\begin{proof}
If there exist a non-empty open subset $V$ of $E^0$ 
and a positive integer $n$ such that 
all vertices in $V$ are base points of simple loops in $E^n$ 
without entrances, 
then $E$ is not topologically free by definition. 

Suppose that $E$ is not topologically free. 
Then we can find a non-empty open set $W_0$ 
such that all vertices in $W_0$ are base points of loops without entrances. 
Set $V_0=H(W_0)$. 
It is not difficult to see 
that $V_0$ is a non-empty open set with $d(r^{-1}(V_0))=V_0$ 
such that all vertices in $V_0$ are base points of loops without entrances. 
We will show that the restriction of $d$ to $r^{-1}(V_0)$ 
is an injection onto $V_0$. 
Take $e_1,e_1'\in r^{-1}(V_0)$ with $d(e_1)=d(e_1')$. 
Let $(e_1,\ldots,e_n)$ 
be a simple loop without entrances 
whose base point is $r(e_1)\in V_0$, 
and $(e_1',\ldots,e_m')$ the one of $r(e_1')\in V_0$. 
We may assume that $n\geq m$ without loss of generality. 
Since $r(e_2)=d(e_1)=d(e_1')=r(e_2')$, 
we have $e_2=e_2'$. 
We also have $e_3=e_3'$ 
because $r(e_3)=d(e_2)=d(e_2')=r(e_3')$. 
Recursively, 
we can see that $e_k=e_k'$ for $k=2,3,\ldots,m$. 
If $n>m$, then we have $e_{m+1}=e_1'$, 
which is impossible because $d(e_1)=d(e_1')$ 
and $(e_1,\ldots,e_n)$ is a simple loop. 
Hence we have $n=m$ and $e_1=e_1'$. 
Thus the restriction of $d$ to $r^{-1}(V_0)$ 
is a bijection onto $V_0$. 
Since $d$ is a local homeomorphism, 
its restriction to the open set $r^{-1}(V_0)$ is a homeomorphism 
from $r^{-1}(V_0)$ to $V_0$. 
Hence we can define a continuous map 
$\sigma\colon V_0\to V_0$ by $\sigma=r\circ (d|_{r^{-1}(V_0)})^{-1}$. 
By Baire's category theorem (see, for example, \cite[Proposition 2.2]{T1}), 
there exist a non-empty open subset $V$ of $V_0$ and 
a positive integer $n$ such that every vertex $v\in V$ satisfies 
$\sigma^k(v)\neq v$ for $k=1,2,\ldots,n-1$ and $\sigma^n(v)=v$. 
This shows that $v\in V$ is a base point of a 
simple loop in $E^n$ without entrances. 
The proof is completed. 
\end{proof}

\begin{proposition}\label{notL}
If $E$ is not topologically free, 
then there exists a non-zero ideal $I$ of $\cO(E)$ 
with $\bigcap_{z\in\T}\beta_z(I)=0$. 
\end{proposition}

\begin{proof}
By Proposition~\ref{Baire}, 
there exist a non-empty open subset $V$ of $E^0$ 
and a positive integer $n$ such that 
all vertices in $V$ are base points of simple loops in $E^n$ 
without entrances, 
and $\sigma=r\circ (d|_{r^{-1}(V)})^{-1}$ is 
a well-defined continuous map from $V$ to $V$ 
with $\sigma^n=\id_{V}$. 
Take $v\in V$. 
We can find a compact neighborhood $V'$ of $v$ such that 
$\sigma^k(v)\notin V'$ for $k=1,2,\ldots,n-1$. 
Take an open set $V''$ such that $v\in V''\subset V'$ 
and define $X=V''\setminus\bigcup_{k=1}^{n-1}\sigma^k(V')$ 
which is a non-empty open subset of $E^0$. 
We have $X\cap \sigma^k(X)=\emptyset$ for $k=1,2,\ldots,n-1$. 
Define $F^0=\bigcup_{k=0}^{n-1}\sigma^k(X)$. 
Then $F^0$ is a hereditary open subset of $E^0$. 
Define a subgraph $F=(F^0,F^1,d|_{F^1},r|_{F^1})$ of $E$ 
such that $F^1=r^{-1}(F^0)$. 
We can apply Lemma~\ref{non-inv-ideal2} 
to the topological graph $F$, 
and hence get a non-zero ideal $I'$ of $\cO(F)$ 
with $I'\cap t_F^0(C_0(F^0))=0$ 
where $t_F^0\colon C_0(F^0)\to \cO(F)$ 
is the natural injection. 
By Proposition~\ref{IrV2}, 
there exists a natural isomorphism 
from $\cO(F)$ to the $C^*$-subalgebra $A$ of $\cO(E)$ 
generated by $t^0(C_0(F^0))$ and $t^1(C_d(F^1))$ 
which commutes with the natural injections 
from $C_0(F^0)$ and $C_d(F^1)$ to $\cO(F)$ and $A$. 
Let $I''\subset A$ be the image of $I'$ 
under this natural isomorphism. 
Thus $I''$ is an ideal of $A$. 
By Lemma~\ref{IG=0}, 
we have $\bigcap_{z\in\T}\beta'_z(I')=0$ 
where $\beta'$ is the gauge action of $\cO(F)$. 
Since the natural isomorphism $\cO(F)\to A\subset \cO(E)$ 
preserves the images of $C_0(F^0)$ and $C_d(F^1)$, 
it is equivariant under the two gauge actions. 
Hence we get $\bigcap_{z\in\T}\beta_z(I'')=0$. 
By Proposition~\ref{IrV2}, 
$A$ is a hereditary and full subalgebra of the ideal $J$
generated by $t^0(C_0(F^0))$. 
Since the map $I\mapsto I\cap A$ is a bijection 
from the set of ideals of $J$ to the ones of $A$, 
we can find an ideal $I$ of $J$ with $I\cap A=I''$. 
The ideal $\bigcap_{z\in\T}\beta_z(I)$ of $J$ satisfies 
$$\bigcap_{z\in\T}\beta_z(I)\cap A
=\bigcap_{z\in\T}\beta_z(I\cap A)
=\bigcap_{z\in\T}\beta_z(I'')=0.$$ 
Hence we get $\bigcap_{z\in\T}\beta_z(I)=0$. 
Thus we get a non-zero ideal $I$ of $\cO(E)$ 
that satisfies $\bigcap_{z\in\T}\beta_z(I)=0$. 
\end{proof}

Now the following theorem strengthens 
the Cuntz-Krieger Uniqueness Theorem. 

\begin{theorem}\label{topfree}
The following conditions for a topological graph $E$ are equivalent:
\benu
\item $E$ is topologically free.
\item The natural surjection $\cO(E)\to C^*(T)$ 
is an isomorphism for every injective Cuntz-Krieger $E$-pair $T$. 
\item $\rho_I=(E^0,\s{E}{0}{sg})$ implies $I=0$. 
\item Any non-zero ideal $I$ of $\cO(E)$ satisfies 
$I\cap t^0(C_0(E^0))\neq 0$.
\eenu
\end{theorem}

\begin{proof}
Clear from Proposition~\ref{CKUT}, 
Lemma~\ref{IG=0} and Proposition~\ref{notL}. 
\end{proof}

\begin{proposition}\label{forprimeness}
If a topological graph $E$ is not topologically free, 
then there exist non-zero ideals $I_1,I_2$ of $\cO(E)$ 
such that $I_1\cap I_2=0$. 
\end{proposition}

\begin{proof}
We can easily see the conclusion 
when $E$ is the topological graph $E_n$ in Example~\ref{easyexam}. 
For general topological graphs 
which are not topologically free, 
we can prove it similarly to the proofs of 
Lemma~\ref{non-inv-ideal2} and 
Proposition~\ref{notL}. 
\end{proof}

\section{Free topological graphs}\label{SecFree}

In this section, 
we give the condition on topological graphs $E$ 
so that all ideals of $\cO(E)$ are gauge-invariant. 

\begin{definition}\label{cond}
For a positive integer $n$, 
we denote by $\Per_n(E)$ the set of vertices $v$ satisfying 
the following three conditions;
\benu
\item there exists a simple loop 
$(l_1,\ldots,l_n)\in E^n$ whose base point is $v$, 
\item 
if $e\in E^1$ satisfies $d(e)\in \Orb^+(v)$ and 
$r(e)=r(l_k)$ for some $k$, 
then we have $e=l_k$, 
\item $v$ is isolated in $\Orb^+(v)$. 
\eenu
We set $\Per(E)=\bigcup_{n=1}^\infty \Per_n(E)$ 
and $\Aper(E)=E^0\setminus\Per(E)$. 
\end{definition}

An element in $\Per(E)$ is called a {\em periodic point} while
an element in $\Aper(E)$ is called an {\em aperiodic point}.
The conditions (i) and (ii) above 
mean that $v\in E^0$ is a base point of exactly one simple loop, 
and the condition (iii) says that there exists no ``approximated loop''
whose base point is $v$. 
When topological graphs come from homeomorphisms, 
these notions coincide with the ordinary ones 
in dynamical systems (see, for example, \cite{T1,T2}). 
Note that $\bigcup_{k=1}^n\Per_k(E)$ is not necessarily closed 
unlike the case of ordinary dynamical systems. 
If $E$ is not topologically free, 
then there exists a non-empty open subset $U$ of $E^0$ 
such that every vertex in $U$ is a base point of a loop 
without entrances. 
Thus we have $U\subset \Per(E)$. 
Hence if $\overline{\Aper(E)}=E^0$, 
then $E$ is topologically free. 
The converse is not true in general
(consider discrete graphs). 

\begin{definition}
A topological graph $E$ is said to be {\em free} 
if $\Aper(E)=E^0$. 
\end{definition}

This definition is a generalization of freeness 
in ordinary dynamical systems.
This is also a generalization of {\em Condition K} 
in the theory of graph algebras (see, for example, \cite{KPRR}). 
In \cite[Definition 9.1]{MT}, 
Muhly and Tomforde define Condition K 
for a topological quiver, 
which coincides with our freeness by Proposition~\ref{free} 
(see also \cite[Proposition 9.9]{MT}). 
Free topological graphs are topologically free. 
In fact, we get a stronger statement (Proposition~\ref{free}). 

\begin{lemma}\label{cond->condL}
Let us take $v\in \Per_n(E)$. 
Let $(l_1,\ldots,l_n)\in E^n$ 
be the unique simple loop 
whose base point is $v$. 
Then we have the following. 
\benu 
\item The closed set $X^0=\overline{\Orb^+(v)}$ is invariant. 
\item $r(l_k)$ is isolated in $X^0$ for $k=1,\ldots,n$. 
\item In the topological graph $X=(X^0,X^1,d|_{X^1},r|_{X^1})$, 
the loop $(l_1,\ldots,l_n)\in X^n$ has no entrances. 
\item The topological graph $X$ is not topologically free.
\eenu
\end{lemma}

\begin{proof}
\benu 
\item 
Since $v$ is the base point of 
the loop $(l_1,\ldots,l_n)\in E^n$, 
Proposition~\ref{X=Orb} implies that $X^0$ is invariant. 
\item 
By the assumption, $r(l_1)=v$ is isolated in $X^0$. 
Take a net $\{v_\lambda\}$ in $\Orb^+(v)$ 
which converges to $r(l_2)$, 
and we will show that $v_\lambda=r(l_2)$ eventually. 
Since $r(l_2)=d(l_1)$ and $d$ is locally homeomorphic, 
we can find a net $\{e_{\lambda}\}\subset E^1$ 
such that $\lim e_{\lambda}=l_1$ and $d(e_{\lambda})=v_\lambda$. 
We have $\lim r(e_{\lambda})=r(l_1)$ and $r(e_{\lambda})\in\Orb^+(v)$. 
Since $r(l_1)$ is isolated in $X^0$, 
we have $r(e_{\lambda})=r(l_1)$ eventually. 
For such $\lambda$, 
we have $e_{\lambda}=l_1$  
by the condition (ii) in Definition~\ref{cond}. 
Therefore we have $v_\lambda=r(l_2)$ eventually. 
This proves that $r(l_2)$ is isolated in $X^0$. 
Recursively we can prove that $r(l_k)$ is isolated in $X^0$ 
for $k=3,\ldots,n$. 
\item 
Take $e\in X^1$ with $r(e)=r(l_k)$ for some $k\in\{1,2,\ldots,n\}$, 
and we will prove $e=l_k$. 
Since $d(e)\in X^0$, 
there exists a net $\{v_\lambda\}\subset\Orb^+(v)$ 
converging to $d(e)$. 
Since $d$ is locally homeomorphic, 
we can find a net $\{e_{\lambda}\}\subset E^1$ 
such that $\lim e_{\lambda}=e$ and $d(e_{\lambda})=v_\lambda$ 
eventually. 
We have $\lim r(e_{\lambda})=r(e)=r(l_k)$ and $r(e_{\lambda})\in\Orb^+(v)$. 
Since $r(l_k)$ is isolated in $X^0$ by (ii), 
we have $r(e_{\lambda})=r(l_k)$ eventually. 
By the condition (ii) in Definition~\ref{cond}, 
we have $e_{\lambda}=l_k$ eventually. 
Hence we have $e=l_k$. 
\item Since $\{v\}$ is an open subset of $X^0$, 
the proof completes by (iii). 
\eenu
\end{proof}

For $v\in \Per_n(E)$, 
Lemma~\ref{cond->condL}~(ii) implies 
that $r(l_k)\in \Per_n(E)$, 
for the unique simple loop $(l_1,\ldots,l_n)\in E^n$ 
whose base point is $v$. 

Recall that for an admissible pair $\rho=(X^0,Z)$ 
we define a topological graph $E_{\rho}$ 
in Definition~\ref{Erho}, 
and we have $E_{\rho}=X$ for $\rho=(X^0,\s{X}{0}{sg})$. 

\begin{lemma}\label{condL->cond}
If there exists an admissible pair $\rho=(X^0,Z)$ 
such that $E_{\rho}$ is not topologically free, 
then $E$ is not free. 
\end{lemma}

\begin{proof}
Let us take an open subset $V$ of 
$E_{\rho}^0=X^0\mathop{\amalg}_{\partial Y_\rho}\overline{Y_\rho}$ 
such that every $v\in V$ is a base point of a loop without entrances. 
Since every vertex in $E_{\rho}^0\setminus X^0$ is a source, 
$V$ is contained in $X^0$. 
Take $v\in V$ arbitrarily 
and let $(l_1,\ldots,l_n)\in E_{\rho}^n$ 
be the unique simple loop 
without entrances whose base point is $v$.
We will prove that $v\in V\subset E^0$ is in $\Per_n(E)$.
Since $d(e)\in X^0\subset E_{\rho}^0$ implies $e\in X^1\subset E_{\rho}^1$ 
for $e\in E_{\rho}^1$, 
we have $l_k\in X^1$ for $k=1,\ldots,n$. 
Let us take $e\in E^1$ with $d(e)\in\Orb^+(v)$ and 
$r(e)=r(l_k)$ for some $k\in\{1,\ldots,n\}$. 
Since $d(e)\in\Orb^+(v)\subset X^0$, 
we have $e\in X^1\subset E_{\rho}^1$. 
Since the loop $(l_1,\ldots,l_n)$ has no entrances 
in the topological graph $E_{\rho}$, 
we have $e=l_k$. 
Thus the condition (ii) in Definition~\ref{cond} is satisfied. 
We will show that $v$ is isolated in $\Orb^+(v)$. 
Let us take a net $\{v_\lambda\}$ in $\Orb^+(v)$ 
which converges to $v$. 
Since $v\in V$, $\Orb^+(v)\subset X^0$ 
and $V$ is open in $X^0$, 
$v_\lambda\in V$ eventually. 
For such $\lambda$, 
$v_\lambda$ is a base point of a loop in $X^n$ 
without entrances. 
Therefore 
since there exists a path from $v$ to $v_\lambda$, 
we can find a path from $v_\lambda$ to $v$. 
Hence $v_\lambda=r(l_k)$ for some $k\in \{1,\ldots,n\}$. 
Thus we get $v_\lambda=v$ eventually. 
This implies that $v$ is isolated in $\Orb^+(v)$. 
Hence we have $v\in \Per_n(E)$ and so $\Per(E)\neq\emptyset$. 
Therefore $E$ is not free.
\end{proof}

\begin{proposition}\label{free}
A topological graph $E$ is free
if and only if $E_{\rho}$ is topologically free 
for every admissible pair $\rho$.
\end{proposition}

\begin{proof}
This follows from Lemma~\ref{cond->condL} and Lemma~\ref{condL->cond}. 
\end{proof}

\begin{theorem}\label{GII&free}
For a topological graph $E$, 
every ideal of $\cO(E)$ is gauge-invariant 
if and only if $E$ is free. 
Thus if $E$ is free, 
the set of all ideals corresponds 
bijectively to the set of all admissible pairs 
by the maps $I\mapsto\rho_I$ and $\rho\mapsto I_\rho$.
\end{theorem}

\begin{proof}
By Proposition~\ref{IXI} and Theorem~\ref{topfree}, 
every ideal of $\cO(E)$ is gauge-invariant 
if and only if $E_\rho$ is topologically free 
for every admissible pair.
It is equivalent to the freeness of $E$ 
by Proposition~\ref{free}. 
The latter statement follows from 
the former and Theorem~\ref{gaugeideals}. 
\end{proof}

\section{Minimal topological graph and simplicity of $\cO(E)$}\label{SecMin}

We give a couple of conditions on $E$ 
all of which are equivalent 
to saying that $\cO(E)$ becomes simple. 
We start with a detailed analysis on 
topological graphs with a periodic point. 

Let $E$ be a topological graph, 
and $v_0$ be a periodic point of $E$. 
Let $l=(l_1,l_2,\ldots,l_n)$ be the unique simple loop 
whose base point is $v_0$. 
We set $X^0=\overline{\Orb^+(v_0)}$ 
which is a closed invariant set of $E^0$. 
We define a topological graph $X=(X^0,X^1,d_X,r_X)$ 
so that $X^1=d^{-1}(X^0)$ and $d_X,r_X$ 
are the restrictions of $d,r$ to $X^1$. 
We simply write $d,r$ for $d_X,r_X$. 
This causes no confusion 
because in the sequel 
we only deal with the topological graph $X$, 
and do not use the topological graph $E$. 
We set $V=\{d(l_1),d(l_2),\ldots,d(l_n)\}$ 
which is an open subset of $X^0$. 

\begin{lemma}\label{infpath}
For an infinite path $e=(e_1,e_2,\ldots)\in X^\infty$ 
of the graph $X$, 
the following conditions are equivalent: 
\benu
\item There exists $k$ such that $d(e_k)\in V$. 
\item There exists $k_0$ such that $d(e_k)\in V$ for every $k\geq k_0$. 
\item $e$ is in the form $(e',l,l,\ldots)$ with some $e'\in X^*$. 
\eenu
\end{lemma}

\begin{proof}
This follows from the fact that $l=(l_1,l_2,\ldots,l_n)$ 
is a loop without entrances in $X$. 
\end{proof}

We denote by $\Lambda_l$ the set of infinite paths $e\in X^\infty$ 
satisfying the equivalent conditions in Lemma~\ref{infpath}. 

The set $V$ coincides with $H(\{v_0\})$ considered 
in the topological graph $X$. 
Let $F^0=S(V)\subset X^0$. 
Since $V$ is open in $X^0$, 
$F^0$ is an open hereditary 
and saturated subset of $X^0$ by Proposition~\ref{hersat2}. 
We define a subgraph $F$ of $X$ 
by $F=(F^0,F^1,d|_{F^1},r|_{F^1})$ where $F^1=r^{-1}(F^0)$. 

\begin{proposition}\label{SHV}
For $v\in X^0$, 
the following are equivalent: 
\benu
\item $v\in F^0$. 
\item $\{v\}$ is open and every negative orbit of $v$ is in $\Lambda_l$. 
\item $\{v\}$ is open and has only finitely many negative orbits. 
\item $\{v\}$ is open and the set $\{e\in\Lambda_l\mid r(e)=v\}$ is finite. 
\eenu
\end{proposition}

\begin{proof}
(i)$\Rightarrow$(ii): 
Let $W$ be the set of $v\in X^0$ satisfying (ii). 
Clearly $V\subset W$. 
We will show that $W$ is saturated. 
Take $v\in \s{X}{0}{rg}$ with $d(r^{-1}(v))\subset W$. 
It is clear that 
every negative orbit of $v$ is in $\Lambda_l$. 
Since $W$ is open, 
Lemma~\ref{RgPr} gives 
a neighborhood $V'$ of $v$ 
such that $r^{-1}(V')\subset d^{-1}(W)$. 
By replacing $V'$ with a smaller set, 
we may assume that $V'$ 
is a compact neighborhood of $v$ 
with $V'\subset \s{X}{0}{rg}$. 
Then $r^{-1}(V')$ 
is compact and satisfies $r(r^{-1}(V'))=V'$ 
(see \cite[Proposition 2.8]{Ka1}). 
Since $d^{-1}(W)$ is a discrete set, 
$r^{-1}(V')$ is a finite set. 
Hence $V'$ is also finite, 
and this shows that $\{v\}$ is open. 
Thus $v\in W$. 
We have shown that $W$ is a saturated set containing $V$. 
By Proposition~\ref{hersat2}, 
we have $F^0\subset W$. 

(ii)$\Rightarrow$(iii): 
Take $v\in X^0$ whose negative orbits are in $\Lambda_l$. 
We have $v\in \s{X}{0}{rg}$. 
Hence $U_1=r^{-1}(v)$ is a compact set. 
We also have $d(U_1)\subset \s{X}{0}{rg}$. 
Since $r\colon r^{-1}(\s{X}{0}{rg})\to \s{X}{0}{rg}$ 
is a proper map, 
$U_2=r^{-1}(d(U_1))$ is a compact set. 
Hence the subsets $U_k\subset X^1$, 
defined by $U_{k+1}=r^{-1}(d(U_k))$ recursively, 
are compact. 
The set of negative orbits of $v$ coincides with the compact set 
$$\Omega=\big\{e=(e_1,e_2,\ldots)\ \big|\ 
e_k\in U_k, d(e_k)=r(e_{k+1})\big\}
\subset U_1\times U_2\times \cdots.$$
For $k=1$, $2$, \ldots, 
we define a subset $\Omega_k$ of $\Omega$ by 
$$\Omega_k=\{e=(e_1,e_2,\ldots)\in \Omega\mid d(e_k)\in V\}.$$
For $e=(e_1,e_2,\ldots)\in \Omega$, 
$d(e_k)\in V$ implies $d(e_{k+1})\in V$. 
Hence we have $\Omega_k\subset\Omega_{k+1}$. 
We also see that 
$\Omega_k$ coincides with 
$\{e\in E^k\mid r^k(e)=v, d^k(e)\in V\}$. 
Since $V$ is discrete and $(r^k)^{-1}(v)$ is compact, 
this set is a finite set. 
Hence $\Omega_k$ is also a finite set. 
Since $V$ is open, $\Omega_k$ is an open subset of $\Omega$ 
for every $k$. 
By (ii), we have $\Omega=\bigcup_{k=1}^\infty \Omega_k$. 
By the compactness of $\Omega$, 
we have $\Omega=\Omega_k$ for some $k$. 
This shows that $\Omega$ is finite. 
Thus the set of negative orbits of $\{v\}$ is a finite set. 

(iii)$\Rightarrow$(iv): Obvious. 

(iv)$\Rightarrow$(i): 
Take $v\in X^0$ satisfying (iv). 
We first show that $v'$ also satisfies (iv) 
for $v'\in d(r^{-1}(v))$. 
It is clear that 
the set $\{e'\in\Lambda_l\mid r(e')=v'\}$ is finite 
for all $v'\in d(r^{-1}(v))$. 
Since $\{v\}$ is open, 
$d(r^{-1}(v))\subset X^0$ is open. 
By the condition (iv), 
$d(r^{-1}(v))\cap \Orb^+(v_0)$ 
is finite. 
Since $\Orb^+(v_0)$ is dense in $X^0$, 
$d(r^{-1}(v))$ is a finite subset of $\Orb^+(v_0)$. 
This shows that $\{v'\}$ is open for every $v'\in d(r^{-1}(v))$. 
Hence $v'$ also satisfies (iv) 
for all $v'\in d(r^{-1}(v))$. 
Next we show $v\in \s{X}{0}{rg}$ 
for $v\in X^0$ satisfying (iv). 
Since $\{v\}$ is open in $X^0$, 
we have $v\in\Orb^+(v_0)$. 
Hence the set $r^{-1}(v)$ is non-empty. 
Since we have already seen that 
$d(r^{-1}(v))$ is a finite subset, 
the non-empty set $r^{-1}(v)$ is finite. 
This shows $v\in\s{X}{0}{rg}$ by \cite[Proposition 2.8]{Ka1}. 

For $e=(e_1,e_2,\ldots)\in \Lambda_l$, 
we define $|e|\in\N$ 
by $|e|=0$ if $r(e)\in V$ 
and $|e|=\max\{k\mid d(e_k)\notin V\}$ otherwise. 
For $v$ satisfying (iv), 
we define $n(v)\in\N$ by 
$$n(v)=\max\big\{|e|\in\N\ \big|\ 
e\in\Lambda_l\text{ with }r(e)=v\big\}.$$
We will show $v\in F^0$ by induction on $n(v)$. 
It is clear when $n(v)=0$. 
Suppose that we have shown $v\in F^0$ 
for $v$ satisfying (iv) with $n(v)\leq k$. 
Take $v$ satisfying (iv) with $n(v)=k+1$. 
As we saw above, 
each $v'\in d(r^{-1}(v))$ satisfies (iv). 
Clearly $n(v')\leq k$ for $v'\in d(r^{-1}(v))$. 
Hence we have $d(r^{-1}(v))\subset F^0$ 
by the assumption of the induction. 
Since $v\in \s{X}{0}{rg}$ and $F^0$ is saturated, 
we have $v\in F^0$. 
We are done. 
\end{proof}

\begin{corollary}\label{CorSHV}
We have $F^0\subset\s{X}{0}{rg}$, 
and $F^0$ is a discrete subset of $X^0$. 
\end{corollary}

Proposition~\ref{SHV} motivates 
the following definition. 

\begin{definition}
We say that a topological graph $E$ is 
{\em generated by a loop} $l=(l_1,l_2,\ldots,l_n)$ 
if $E^0$ is discrete and 
every negative orbit is in the form 
$(e',l,l,\ldots)\in E^\infty$ with some $e'\in E^*$. 
\end{definition}

We have already seen the following. 

\begin{proposition}\label{GenByLoop}
Let $E$ be a topological graph. 
For $v\in\Per(E)$, 
the topological graph $F$ defined above is generated by the loop $l$. 
\end{proposition}

Every topological graph generated by a loop 
arises in this way. 

\begin{proposition}
Let $E$ be a topological graph 
generated by a loop $l=(l_1,l_2,\ldots,l_n)$. 
Then the base point $v_0$ of the loop $l$ is in $\Per_n(E)$, 
and $E^0=\Orb^+(v_0)=S(H(\{v_0\}))$. 
\end{proposition}

\begin{proof}
It is clear from the definition 
that $l$ has no entrances. 
Hence $v_0\in \Per_n(E)$. 
It is also clear to see $E^0=\Orb^+(v_0)$. 
Now $E^0=S(H(\{v_0\}))$ 
follows from Proposition~\ref{SHV}. 
\end{proof}

For a Hilbert space $H$, 
we denote by $K(H)$ (resp.\ $B(H)$) 
the $C^*$-algebra of all compact operators 
(resp.\ all bounded operators) on $H$. 

\begin{proposition}\label{C(T,K)}
Let $E$ be a topological graph generated by a loop. 
Then $\cO(E)$ is isomorphic to $C(\T)\otimes K(\ell^2(E^\infty))$. 
\end{proposition}

\begin{proof}
Let $E$ be a topological graph generated 
by a loop $l=(l_1,l_2,\ldots,l_n)$. 
We denote by $l^\infty\in E^\infty$ the infinite path $(l,l,\ldots)$. 
Let us denote by $w\in C(\T)$ 
the generating unitary of $C(\T)$ 
defined by $w(z)=z$ for $z\in \T$. 
Let us denote by $\{u_{e,e'}\}_{e,e'\in E^\infty}$ 
the canonical matrix units of $K(\ell^2(E^\infty))$. 
By Proposition~\ref{SHV}, 
for $v\in E^0$ 
the set $E^\infty_v=\{e\in E^\infty\mid r(e)=v\}$ is finite. 
Since $E^0$ and $E^1$ are discrete, 
the linear span of the characteristic functions 
$\{\delta_v\}_{v\in E^0}$ is dense in $C_0(E^0)$, 
and the one of $\{\delta_e\}_{e\in E^1}$ is dense in $C_d(E^1)$. 
We define a $*$-homomorphism 
$T^0\colon C_0(E^0)\to C(\T)\otimes K(\ell^2(E^\infty))$ 
and a linear map $T^1\colon C_d(E^1)\to C(\T)\otimes K(\ell^2(E^\infty))$ 
so that for $v\in E^0$ and $e\in E^1\setminus\{l_1\}$, 
$$T^0(\delta_v)=1\otimes \sum_{e'\in E^\infty_v}u_{e',e'}, \quad 
T^1(\delta_e)=1\otimes \sum_{e'\in E^\infty_{d(e)}}u_{(e,e'),e'}$$
and $T^1(\delta_{l_1})=w\otimes u_{(l_1,l_2,\ldots),(l_2,l_3,\ldots)}$. 
It is routine to check 
\begin{align*}
T^1(\delta_e)^*T^1(\delta_{e'})&=\delta_{e,e'}T^0(\delta_{d(e)}),\quad
T^0(\delta_v)T^1(\delta_e)=\delta_{v,r(e)}T^1(\delta_e)\\
T^0(\delta_v)=&\sum_{e\in r^{-1}(v)}T^1(\delta_e)T^1(\delta_e)^*
\end{align*}
for $e,e'\in E^1$ and $v\in E^0$. 
This shows that $T=(T^0,T^1)$ is a Cuntz-Krieger $E$-pair.
Since $E^\infty_v\neq\emptyset$ for all $v\in E^0$, 
$T$ is injective. 
We will show that $T$ admits a gauge action. 
For $e=(e_1,e_2,\ldots)\in E^\infty$, 
we define $|e|\in\N$ 
by $|e|=\min\{k \mid e_k=l_1\}$. 
For $z\in\T$, 
we define an automorphism $\sigma_z$ of $C(\T)$ 
by $\sigma_z(f)(z')=f(z^nz')$ for $z'\in\T$, 
and an automorphism $\Ad(u_z)$ of $K(\ell^2(E^\infty))$ 
by $\Ad(u_z)(x)=u_zxu_z^*$ 
where $u_z\in B(\ell^2(E^\infty))$ is a unitary defined 
by $u_z=\sum_{e\in E^\infty}z^{|e|}u_{e,e}$ which converges strongly. 
Then the automorphism $\beta'_z=\sigma_z\otimes \Ad(u_z)$ of 
$C(\T)\otimes K(\ell^2(E^\infty))$ satisfies 
$$\beta'_z(1\otimes u_{e,e'})=z^{|e|-|e'|}(1\otimes u_{e,e'}),\quad 
\beta'_z(w\otimes u_{l^\infty,l^\infty})=z^n(w\otimes u_{l^\infty,l^\infty}).$$
Thus the action $\beta'$ is a gauge action for the pair $T$. 
Hence $T$ induces 
an injective map $\cO(E)\to C(\T)\otimes K(\ell^2(E^\infty))$. 
We will show that it is surjective. 
For $e=(e_1,e_2,\ldots)\in E^\infty$, 
we have 
$$T^1(\delta_{e_1})T^1(\delta_{e_2})\cdots T^1(\delta_{e_{k-1}})
T^0(\delta_{d(l_1)})=1\otimes u_{e,l^\infty}$$
where $k=|e|$. 
We also have 
$$T^1(\delta_{l_1})T^1(\delta_{l_2})\cdots T^1(\delta_{l_n})
=w\otimes u_{l^\infty,l^\infty}.$$
Since $C(\T)\otimes K(\ell^2(E^\infty))$ 
is generated by $\{1\otimes u_{e,l^\infty}\}_{e\in E^\infty}$ 
and $w\otimes u_{l^\infty,l^\infty}$, 
we see that the map 
$\cO(E)\to C(\T)\otimes K(\ell^2(E^\infty))$ is surjective. 
Therefore $\cO(E)$ is isomorphic to $C(\T)\otimes K(\ell^2(E^\infty))$. 
\end{proof}

From now on, 
we study for which topological graph $E$, 
the $C^*$-algebra $\cO(E)$ is simple. 
We introduce the following notion. 

\begin{definition}
A topological graph $E$ is said to be {\em minimal} 
if there exist no closed invariant sets other than $\emptyset$ or $E^0$. 
\end{definition}

\begin{proposition}\label{minimal}
For a topological graph $E$, 
the following conditions are equivalent: 
\benu
\item $E$ is minimal. 
\item The orbit space $\Orb(v,e)$ is dense in $E^0$ 
for every $v\in E^0$ and every negative orbit $e$ of $v$. 
\item For every non-empty open set $V\subset E^0$, 
we have $S(H(V))=E^0$. 
\eenu
\end{proposition}

\begin{proof}
(i)$\Rightarrow$(ii): 
For a negative orbit $e$ of $v\in E^0$, 
the closed subset $\overline{\Orb(v,e)}$ is invariant. 
By (i), we have $\overline{\Orb(v,e)}=E^0$. 

(ii)$\Rightarrow$(i): 
Let $X^0$ be a non-empty closed invariant subset of $E^0$. 
Take $v\in X^0$. 
By Proposition~\ref{InvOrb}, 
there exists a negative orbit $e$ of $v$ 
such that $\Orb(v,e)\subset X^0$. 
By (ii), we have $\overline{\Orb(v,e)}=E^0$. 
Hence $X^0=E^0$. 
Thus $E$ is minimal. 

(i)$\Rightarrow$(iii): 
Take a non-empty open set $V\subset E^0$. 
Then $E^0\setminus S(H(V))$ is 
a closed invariant subset 
such that $E^0\setminus S(H(V))\neq E^0$. 
By the minimality, we have $E^0\setminus S(H(V))=\emptyset$. 
Thus we have $S(H(V))=E^0$. 

(iii)$\Rightarrow$(i): 
Let $X^0$ be a closed invariant subset of $E^0$ 
with $X^0\neq E^0$. 
Then $V=E^0\setminus X^0$ is 
a non-empty hereditary and saturated set. 
By (iii), we have $V=S(H(V))=E^0$. 
Hence we have $X^0=\emptyset$. 
Thus $E$ is minimal. 
\end{proof}

\begin{remark}
The notion of minimality extends the one of ordinary dynamical systems 
for which Proposition~\ref{minimal} is well-known. 
When a graph is discrete, 
the condition (ii) above is equivalent to 
the cofinality in the sense of \cite{BPRS} 
with the extra condition that 
for two vertices $v\in E^0$ and $w\in\s{E}{0}{sg}$ 
there exists a path $e\in E^*$ such that $d(e)=w$ and $r(e)=v$. 
For a discrete graph, 
the equivalence of (ii) and (iii) in Proposition~\ref{minimal} 
had been certainly known (see, for example, Introduction of \cite{Pa}).
\end{remark}

\begin{lemma}
A topological graph generated by a loop is minimal. 
\end{lemma}

\begin{proof}
Let $E$ be a topological graph generated by a loop $l$, 
and $v_0\in E^0$ be the base point of the loop $l$. 
By definition, 
for any $v\in E^0$ and any negative orbit $e$ of $v$, 
the negative orbit space $\Orb^-(v,e)$ contains $v_0$. 
Since $\Orb^+(v_0)=E^0$, $E$ is minimal 
by Proposition~\ref{minimal}. 
\end{proof}

\begin{theorem}\label{ThmSimple}
For a topological graph $E$, 
the following conditions are equivalent: 
\benu
\item The $C^*$-algebra $\cO(E)$ is simple. 
\item $E$ is minimal and topologically free. 
\item $E$ is minimal and free. 
\item $E$ is minimal and not generated by a loop. 
\eenu
\end{theorem}

\begin{proof}
(i)$\Rightarrow$(iv): 
If $E$ is not minimal, 
then there exists a non-trivial gauge-invariant ideal. 
Hence $\cO(E)$ is not simple. 
If $E$ is generated by a loop, 
$\cO(E)$ is not simple by Proposition~\ref{C(T,K)}. 

(iv)$\Rightarrow$(iii): 
Suppose that $E$ is minimal and not free. 
Take $v_0\in\Per(E)$. 
Since $E$ is minimal, 
$\overline{\Orb^+(v_0)}=E^0$ by Proposition~\ref{minimal}. 
Hence $\{v_0\}$ is open in $E^0$. 
Using Proposition~\ref{minimal} again, 
we get $S(H(\{v_0\}))=E^0$. 
Hence $E$ is generated by a loop 
by Proposition~\ref{GenByLoop}. 

(iii)$\Rightarrow$(ii): 
Obvious. 

(ii)$\Rightarrow$(i): 
Assume that $E$ is minimal and topologically free. 
Take an ideal $I$ of $\cO(E)$ with $I\neq \cO(E)$. 
Then $X_I$ is a closed invariant set other than $\emptyset$. 
By the minimality, we have $X_I=E^0$. 
By Theorem~\ref{topfree} we have $I=0$. 
Thus $\cO(E)$ is simple. 
\end{proof}

\begin{corollary}
When $E^0$ is not discrete, 
$\cO(E)$ is simple if and only if $E$ is minimal. 
\end{corollary}

\begin{remark}
The above theorem generalizes the result 
on simplicity of graph algebras 
(\cite[Theorem 12]{S}, \cite[Corollary 2.15]{DT}, 
see also \cite[Theorem 4]{Pa}). 
In the case that $d$ is injective, 
a topological graph $E=(E^0,E^1,d,r)$ is generated by a loop 
if and only if $E$ is minimal and $E^0$ is a finite set. 
Hence in this case the condition (iv) in the above theorem 
is equivalent to 
\begin{itemize}
\item[(iv)'] $E$ is minimal and $E^0$ is an infinite set. 
\end{itemize}
Thus Theorem~\ref{ThmSimple} generalizes a criterion for simplicity 
of homeomorphism $C^*$-algebras due to Zeller-Meier \cite{Z}. 
\end{remark}

\section{Primeness for admissible pairs}\label{SecPrimePair}

In this section, we give a necessary condition for an ideal to be prime 
in terms of admissible pairs. 
We will use it after in order to determine all prime ideals 
(Theorem~\ref{ThmPrime}). 
Recall that an ideal $I$ of a $C^*$-algebra $A$ is 
said to be {\em prime} 
if for two ideals $I_1,I_2$ of $A$, 
$I_1\cap I_2\subset I$ implies either $I_1\subset I$ or $I_2\subset I$. 
We define primeness for admissible pairs. 

\begin{definition}
An admissible pair $\rho$ is called {\em prime} 
if $\rho_1\cup \rho_2\supset\rho$ 
implies either $\rho_1\supset\rho$ 
or $\rho_2\supset\rho$
for two admissible pairs 
$\rho_1,\ \rho_2$.
\end{definition}

It is well-known that an ideal $I$ is prime if and only if 
the equality $I_1\cap I_2=I$ implies either $I_1=I$ or $I_2=I$ 
for two ideals $I_1,I_2$ 
(see the proof of (iii)$\Rightarrow$(iv) of Proposition~\ref{primepair}).
The following is the counterpart of this fact 
for prime admissible pairs.

\begin{proposition}\label{primepair}
For an admissible pair $\rho$, 
the following are equivalent: 
\benu
\item $\rho$ is prime.
\item 
For two admissible pairs $\rho_1$, 
$\rho_2$,
the equality $\rho_1\cup \rho_2=\rho$
implies either $\rho_1=\rho$ 
or $\rho_2=\rho$.
\item For two gauge invariant ideals $I_1,I_2$ of $\cO(E)$,
the equality $I_1\cap I_2=I_{\rho}$ implies either
$I_1=I_{\rho}$ or $I_2=I_{\rho}$.
\item For two gauge invariant ideals $I_1,I_2$ of $\cO(E)$, 
the inclusion $I_1\cap I_2\subset I_{\rho}$ implies 
either $I_1\subset I_{\rho}$ or $I_2\subset I_{\rho}$.
\eenu
\end{proposition}

\begin{proof}
(i)$\Rightarrow$(ii): 
Take two admissible pairs $\rho_1$, 
$\rho_2$ with $\rho_1\cup \rho_2=\rho$.
By (i), we have either $\rho_1\supset \rho$ 
or $\rho_2\supset \rho$.
Hence we get either $\rho_1=\rho$ 
or $\rho_2=\rho$.

(ii)$\Rightarrow$(iii): 
Take two gauge invariant ideals $I_1,I_2$
with $I_1\cap I_2=I_{\rho}$.
We have $\rho_{I_1}\cup \rho_{I_2}=\rho$.
By (ii), we have either $\rho_{I_1}=\rho$ or
$\rho_{I_2}=\rho$.
By Proposition~\ref{IXI},
we have either $I_1=I_{\rho}$ or $I_2=I_{\rho}$.

(iii)$\Rightarrow$(iv): 
Take two gauge invariant ideals $I_1,I_2$
with $I_1\cap I_2\subset I_{\rho}$.
Then we have 
$$(I_1+I_{\rho})\cap (I_2+I_{\rho})=
(I_1\cap I_2) +I_{\rho}=I_{\rho}.$$
By (iii), either $I_1+I_{\rho}=I_{\rho}$ or 
$I_2+I_{\rho}=I_{\rho}$ holds.
Hence we get 
either $I_1\subset I_{\rho}$ or $I_2\subset I_{\rho}$.

(iv)$\Rightarrow$(i): 
Take two admissible pairs 
$\rho_1$, $\rho_2$ with
$\rho_1\cup \rho_2\supset \rho$.
The two gauge-invariant ideals 
$I_{\rho_1}$ and $I_{\rho_2}$ satisfy
$$I_{\rho_1}\cap I_{\rho_2}=I_{\rho_1\cup \rho_2}\subset 
I_{\rho}.$$
Hence we have either $I_{\rho_1}\subset I_{\rho}$ 
or $I_{\rho_2}\subset I_{\rho}$.
By Theorem~\ref{gaugeideals}, 
we get either $\rho_1\supset \rho$ or $\rho_2\supset\rho$.
Thus $\rho$ is prime.
\end{proof}

We will use the implication (ii)$\Rightarrow$(i) 
to determine which admissible pairs are prime.

\begin{proposition}\label{prime}
If an ideal $I$ of $\cO(E)$ is prime, 
then $\rho_I$ is a prime admissible pair.
\end{proposition}

\begin{proof}
The proof goes similarly 
as the proof of (iv)$\Rightarrow$(i) in Proposition~\ref{primepair}, 
hence we omit it.
\end{proof}

In general, the converse of Proposition~\ref{prime} is not true 
(see Propositions~\ref{prime2} and \ref{NotTopFree}). 
We will determine all prime admissible pairs 
(Proposition~\ref{pp}). 

\begin{lemma}\label{primepair0}
If an admissible pair $\rho=(X^0,Z)$ is prime, 
then either $Z=\s{X}{0}{sg}$ 
or $Z=\s{X}{0}{sg}\cup\{v\}$ for some $v\in \s{X}{0}{rg}$.
\end{lemma}

\begin{proof}
Let $\rho=(X^0,Z)$ be a prime admissible pair.
To derive a contradiction, assume $Z\setminus \s{X}{0}{sg}$ has 
two elements $v_1,v_2$.
Take open sets $V_1\ni v_1$, $V_2\ni v_2$ with $V_1\cap V_2=\emptyset$,
$V_1\cap \s{X}{0}{sg}=\emptyset$ and $V_2\cap \s{X}{0}{sg}=\emptyset$.
Then $\rho_i=(X^0,Z\setminus V_i)\ (i=1,2)$ 
are admissible pairs satisfying 
$\rho=\rho_1\cup \rho_2$.
However, we have $\rho\not\subset \rho_1$ and 
$\rho\not\subset \rho_2$.
This contradicts the primeness of $\rho$.
\end{proof}

First we consider the case $Z=\s{X}{0}{sg}$. 

\begin{lemma}\label{primepair1}
An admissible pair $\rho=(X^0,\s{X}{0}{sg})$ is prime if and only if
$X^0$ is a maximal head.
\end{lemma}

\begin{proof}
Suppose that $\rho=(X^0,\s{X}{0}{sg})$ is a prime admissible pair.
Take invariant sets $X^0_1,X^0_2$ 
with $X^0\subset X^0_1\cup X^0_2$.
Set $\rho_i=(X^0_i, X^0_i\cap \s{E}{0}{sg})$ for $i=1,2$. 
Since 
$$\s{X}{0}{sg}\subset X^0\cap \s{E}{0}{sg}=
(X^0_1\cap \s{E}{0}{sg})\cup (X^0_2\cap \s{E}{0}{sg}),$$
We have $\rho\subset \rho_1\cup \rho_2$. 
Since $\rho$ is prime, 
either $\rho\subset \rho_1$ 
or $\rho\subset \rho_2$ holds. 
Hence we have either $X^0\subset X^0_1$ or $X^0\subset X^0_2$. 
By Proposition~\ref{maxhead}, $X^0$ is a maximal head. 
Conversely assume that $X^0$ is a maximal head. 
Take two admissible pairs $\rho_1=(X^0_1,Z_1)$, $\rho_2=(X^0_2,Z_2)$ with 
$\rho_1\cup \rho_2=\rho$.
By Proposition~\ref{maxhead}, 
either $X^0\subset X^0_1$ or $X^0\subset X^0_2$ holds.
We may assume $X^0\subset X^0_1$.
Then $X^0=X^0_1$.
Hence $\s{X}{0}{sg}=(X_1^0)_{\rs{sg}}\subset Z_1\subset\s{X}{0}{sg}$.
We get $\rho_1=\rho$.
By Proposition~\ref{primepair}, 
$\rho=(X^0,\s{X}{0}{sg})$ is a prime admissible pair.
\end{proof}

Next we consider the case $Z=\s{X}{0}{sg}\cup\{v\}$. 

\begin{lemma}\label{primepair2}
An admissible pair $\rho=(X^0,\s{X}{0}{sg}\cup\{v\})$ is prime 
for some $v\in \s{X}{0}{rg}\cap\s{E}{0}{sg}$ 
if and only if $X^0=\overline{\Orb^+(v)}$.
\end{lemma}

\begin{proof}
Suppose that an admissible pair 
$\rho=(X^0,\s{X}{0}{sg}\cup\{v\})$ is prime.
Since $v\in \s{E}{0}{sg}$, 
the pair 
$\rho_1=(\overline{\Orb^+(v)},\overline{\Orb^+(v)}\cap\s{E}{0}{sg})$
is an admissible pair by Proposition~\ref{X=Orb}. 
If we set an admissible pair $\rho_2=(X^0,\s{X}{0}{sg})$, 
then $\rho_1$ and $\rho_2$ satisfy $\rho\subset \rho_1\cup \rho_2$.
Since $\rho$ is prime and $\rho\not\subset\rho_2$, 
we have $\rho\subset\rho_1$. 
Hence $\overline{\Orb^+(v)}\subset X^0\subset\overline{\Orb^+(v)}$.
Thus, we get $X^0=\overline{\Orb^+(v)}$.
Conversely, assume $X^0=\overline{\Orb^+(v)}$.
Take two admissible pairs $\rho_1=(X^0_1,Z_1)$, 
$\rho_2=(X^0_2,Z_2)$ with 
$\rho_1\cup \rho_2=\rho$.
We may assume $v\in Z_1$.
Then we have $X^0=\overline{\Orb^+(v)}\subset X_1^0\subset X^0$.
Hence $X^0_1=X^0$.
We have $\s{X}{0}{sg}\cup\{v\}=(X_1^0)_{\rs{sg}}\cup\{v\}\subset 
Z_1\subset \s{X}{0}{sg}\cup\{v\}$.
Therefore $\rho_1=\rho$.
By Proposition~\ref{primepair}, 
$\rho$ is a prime admissible pair.
\end{proof}

Note that if $v\in\s{E}{0}{sg}$ 
then $X^0=\overline{\Orb^+(v)}$ is invariant 
by Proposition~\ref{X=Orb}. 
We define 
$$BV(E)=\{v\in\s{E}{0}{sg}\mid 
  v\in \s{X}{0}{rg}\mbox{ where }X^0=\overline{\Orb^+(v)}\}.$$
Elements in $BV(E)$ are called {\em breaking vertices}.
When $E$ is discrete, 
breaking vertices are infinite receivers (\cite{BHRS}). 
In general, breaking vertices may not be in $\s{E}{0}{inf}$. 

\begin{example}
Let $E$ be the topological graph in Example~\ref{AbsVal}.
The vertex $0\in\R$ is a breaking vertex. 
We have $\s{E}{0}{inf}=\emptyset$ and 
$0\in\overline{\s{E}{0}{sce}}$. 
\end{example}

For $v\in BV(E)$, 
we define $\rho_v=(X^0,\s{X}{0}{sg}\cup\{v\})$ 
where $X^0=\overline{\Orb^+(v)}$. 
We denote by $\cM(E)$ the set of all maximal heads. 
For $X^0\in \cM(E)$, 
we define $\rho_{X^0}=(X^0,\s{X}{0}{sg})$. 

\begin{proposition}\label{pp}
The map $x\mapsto \rho_x$ gives a bijection 
from $BV(E)\amalg \cM(E)$ to 
the set of all prime admissible pairs. 
\end{proposition}

\begin{proof}
The injectivity of the map is easy to see from the definitions, 
and the surjectivity follows from 
Lemma~\ref{primepair0}, Lemma~\ref{primepair1} 
and Lemma~\ref{primepair2}. 
\end{proof}

\section{Primeness of $\cO(E)$}\label{SecPrime}

A $C^*$-algebra is said to be {\em prime} 
if $0$ is a prime ideal. 
Using the results in the previous section, 
we give conditions on $E$ 
for $\cO(E)$ to be a prime $C^*$-algebra. 

\begin{definition}
A topological graph $E$ is called {\em topologically transitive} 
if we have $H(V_1)\cap H(V_2)\neq\emptyset$ 
for any two non-empty open sets $V_1,V_2\subset E^0$. 
\end{definition}

There exist many equivalent conditions 
of topological transitivity. 

\begin{proposition}\label{TopTran}
For a topological graph $E$, consider the following conditions. 
\benu
\item $E$ is topologically transitive. 
\item For two non-empty open sets $V_1,V_2\subset E^0$, 
we have $S(H(V_1))\cap S(H(V_2))\neq\emptyset$. 
\item The admissible pair $(E^0,\s{E}{0}{sg})$ is prime.
\item The set $E^0$ is a maximal head. 
\item There exists $v\in E^0$ and a negative orbit 
$e$ of $v$ such that $\Orb(v,e)$ is dense in $E^0$. 
\eenu
Then the conditions {\em (i), (ii), (iii)} and {\em (iv)} are equivalent 
and implied by {\em (v)}.
When $E^0$ is second countable, 
the five conditions are equivalent.
\end{proposition}

\begin{proof}
Lemma~\ref{hersat7} gives (i)$\Leftrightarrow$(ii). 
By the definition of maximal heads, we have (i)$\Leftrightarrow$(iv). 
By Lemma~\ref{primepair1}, we have (iii)$\Leftrightarrow$(iv). 
By Proposition~\ref{Orbit=MH}, (v) implies (iv). 
When $E^0$ is second countable, we have (iv)$\Rightarrow$(v) 
by Proposition~\ref{Orbit=MH2}. 
\end{proof}

By this proposition, 
we can see that a minimal topological graph is 
topologically transitive. 

\begin{theorem}\label{C*prime}
A $C^*$-algebra $\cO(E)$ is prime 
if and only if $E$ is topologically free 
and topologically transitive. 
\end{theorem}

\begin{proof}
Suppose that the $C^*$-algebra $\cO(E)$ is prime. 
By Proposition~\ref{forprimeness}, 
$E$ is topologically free. 
The admissible pair $\rho_0=(E^0,\s{E}{0}{sg})$ is prime 
by Proposition~\ref{prime}. 
Hence $E$ is topologically transitive. 
Conversely assume that $E$ is topologically free 
and topologically transitive. 
Take two ideals $I_1,I_2$ of $\cO(E)$ with $I_1\cap I_2=0$, 
and we will show that either $I_1=0$ or $I_2=0$. 
Since $I_1\cap I_2=0$, 
we have $\rho_{I_1}\cup \rho_{I_1}=\rho_{0}=(E^0,\s{E}{0}{sg})$. 
By Proposition~\ref{TopTran}, 
$(E^0,\s{E}{0}{sg})$ is prime. 
Hence either $(E^0,\s{E}{0}{sg})=\rho_{I_1}$ 
or $(E^0,\s{E}{0}{sg})=\rho_{I_2}$ holds. 
We have either $I_1=0$ or $I_2=0$ by Theorem~\ref{topfree}. 
Thus we show that $\cO(E)$ is prime.
\end{proof}

In Proposition~\ref{NotTopFree}, 
we will see that 
a topologically transitive topological graph 
is not topologically free only when 
there exists $v\in\Per(E)$ such that $E^0=\overline{\Orb^+(v)}$.

\section{Prime ideals}\label{SecPrimeIdeal}

In this section, 
we completely determine the set of prime ideals of 
the $C^*$-algebra $\cO(E)$ of a topological graph $E$. 
In Proposition~\ref{prime}, 
we see that for a prime ideal $P$ of $\cO(E)$, 
the admissible pair $\rho_P$ is necessarily prime. 
The following proposition determines 
when the converse of this fact is true 
for gauge-invariant ideals. 

\begin{proposition}\label{prime2}
For a prime admissible pair $\rho$, 
$I_\rho$ is a prime ideal 
if and only if $E_\rho$ is topologically free. 
\end{proposition}

\begin{proof}
If $E_\rho$ is not topologically free, 
then $I_\rho$ is not prime by Proposition~\ref{forprimeness}. 
Suppose $E_\rho$ is topologically free. 
Take two ideals $I_1,I_2$ of $\cO(E)$ with $I_1\cap I_2=I_\rho$. 
By Proposition~\ref{XIX} and Proposition~\ref{cap}, 
we have $\rho=\rho_{I_\rho}=\rho_{I_1\cap I_2}=\rho_{I_1}\cup \rho_{I_2}$. 
Since $\rho$ is a prime admissible pair, 
either $\rho_{I_1}=\rho$ or $\rho_{I_2}=\rho$ by Proposition~\ref{primepair}. 
Without loss of generality, 
we may assume $\rho_{I_1}=\rho$. 
Then we have $E_{\rho_{I_1}}=E_\rho$ which is topologically free. 
Hence by Proposition~\ref{TopFreeQ} we have $I_1=I_\rho$. 
Thus $I_\rho$ is a prime ideal. 
\end{proof}

We study which prime admissible pair $\rho$ satisfies 
that $E_\rho$ is topologically free. 
Recall that in Proposition~\ref{pp} we saw that 
all prime admissible pairs are in the form 
$\rho_v=(\overline{\Orb^+(v)},\overline{\Orb^+(v)}_{\rs{sg}}\cup\{v\})$ 
for $v\in BV(E)$, 
or $\rho_{X^0}=(X^0,\s{X}{0}{sg})$ 
for $X^0\in \cM(E)$. 

\begin{proposition}\label{PropTFBV}
For $v\in BV(E)$, 
the topological graph $E_{\rho_v}$ is topologically free. 
\end{proposition}

\begin{proof}
From Definition~\ref{Erho}, 
we see that 
$E_{\rho_v}^0$ is the disjoint union of $\overline{\Orb^+(v)}$ 
and an extra vertex $\bar{v}$. 
Since the vertex $\bar{v}$ receives no edge, 
it is not a base point of a loop. 
For every vertex $w\in \Orb^+(v)\subset E_{\rho_v}^0$, 
there exists $e\in E_{\rho_v}^*$ 
such that $d(e)=\bar{v}$ and $r(e)=w$. 
Hence $w$ cannot be a base point of a loop without entrances. 
Since $\Orb^+(v)\cup\{\bar{v}\}$ is dense in $E_{\rho_v}^0$, 
the set of base points of loops without entrances
has an empty interior. 
Thus the topological graph $E_{\rho_v}$ is topologically free. 
\end{proof}

We denote 
$$\cM_{\per}(E)
=\big\{\overline{\Orb^+(v)}\in\cM(E)\ \big|\ v\in\Per(E)\big\},$$ 
and $\cM_{\aper}(E)=\cM(E)\setminus \cM_{\per}(E)$. 

\begin{proposition}\label{NotTopFree}
For $X^0\in \cM(E)$, 
the topological graph $E_{\rho_{X^0}}=X$ is topologically free
if and only if $X^0\in\cM_{\aper}(E)$. 
\end{proposition}

\begin{proof}
If $X^0=\overline{\Orb^+(v)}$ for some $v\in\Per(E)$,  
then $X$ is not topologically free by Lemma~\ref{cond->condL} (iv). 
Suppose that $X^0$ is a maximal head such that 
the topological graph $X$ is not topologically free, 
and we will show that $X^0\in \cM_{\per}(E)$. 
By Proposition~\ref{Baire}, 
there exist a non-empty open subset $V$ of $X^0$ 
and a positive integer $n$ such that 
all vertices in $V$ are base points of simple loops in $X^n$ 
without entrances, 
and that $\sigma=r\circ (d|_{r^{-1}(V)})^{-1}$ is 
a well-defined continuous map from $V$ to $V$ 
with $\sigma^n=\id_{V}$. 
Take $v\in V$ arbitrarily. 
We will show that $V=\{v,\sigma(v),\ldots,\sigma^{n-1}(v)\}$. 
Take $w\in V$. 
Since $X^0$ is a maximal head, 
we have two nets $\{e_\lambda\},\{e_\lambda'\}\subset X^*$ 
such that $d(e_\lambda)=d(e_\lambda')\in X^0$ and 
$\lim r(e_\lambda)=v$ and $\lim r(e_\lambda')=w$. 
We may assume that $r(e_\lambda),r(e_\lambda')\in V$ 
for every $\lambda$. 
Since $r(e_\lambda')\in V$ is a base point of a simple loop in $X^n$
without entrances, 
we can find a path from $r(e_\lambda')$ to 
$d(e_\lambda')=d(e_\lambda)$. 
Thus there exists a path from $r(e_\lambda')$ to $r(e_\lambda)$. 
Since $r(e_\lambda)\in V$ is a base point of a simple loop in $X^n$
without entrances, 
we can find $k_\lambda\in\{0,1,\ldots,n-1\}$ 
such that $r(e_\lambda')=\sigma^{k_\lambda}(r(e_\lambda))$ 
for each $\lambda$. 
Then, we can find $k\in\{0,1,\ldots,n-1\}$ 
with $k_{\lambda}=k$ frequently. 
We have 
\begin{align*}
w&=\lim r(e_{\lambda}')
=\lim \sigma^k(r(e_{\lambda}))
=\sigma^k(\lim r(e_{\lambda}))
=\sigma^k(v).
\end{align*}
Thus we have shown that $V=\{v,\sigma(v),\ldots,\sigma^{n-1}(v)\}$. 
Since $V$ is open in $X^0$, 
$\{v\}$ is also open in $X^0$. 
Hence $\{v\}$ is isolated in $\Orb^+(v)$ 
because $\Orb^+(v)\subset X^0$, 
Since the condition (ii) in Definition~\ref{cond} 
is clearly satisfied, 
we have $v\in\Per_n(E)$. 
The proof will complete once we will show that $X^0=\overline{\Orb^+(v)}$. 
Clearly $X^0\supset \overline{\Orb^+(v)}$. 
Take $w\in X^0$. 
By noting that $\{v\}$ is a neighborhood of $v\in X^0$, 
we can find a net $\{e_\lambda\}\subset X^*$ 
such that $\lim r(e_\lambda)=w$ and $v\in\Orb^+(d(e_\lambda))$ 
because $X^0$ is a maximal head. 
Since $v$ is a base point of a loop without entrances in $X$, 
$v\in\Orb^+(d(e_\lambda))$ 
implies that $d(e_\lambda)\in\Orb^+(v)$. 
Hence we have $r(e_\lambda)\in\Orb^+(v)$. 
This implies that $w\in \overline{\Orb^+(v)}$. 
Therefore we have $X^0=\overline{\Orb^+(v)}$ for $v\in\Per_n(E)$.
Thus $X^0\in \cM_{\per}(E)$. 
\end{proof}

\begin{definition}
For $v\in BV(E)$ and $X^0\in\cM_{\aper}(E)$, 
we define $P_v=I_{\rho_v}$ and $P_{X^0}=I_{\rho_{X^0}}$. 
\end{definition}

\begin{proposition}\label{Pxuniq}
The ideal $P_x$ for $x\in BV(E)\amalg\cM_{\aper}(E)$ 
is a prime ideal. 
If an ideal $I$ satisfies $\rho_I=\rho_x$ 
for $x\in BV(E)\amalg\cM_{\aper}(E)$, 
then $I=P_x$. 
\end{proposition}

\begin{proof}
This follows from 
Proposition~\ref{prime2}, Proposition~\ref{PropTFBV} and 
Proposition~\ref{NotTopFree}. 
\end{proof}

In order to list all the prime ideals, 
the only remaining thing to do 
is to investigate the prime ideals $P$ 
with $\rho_P=(X^0,\s{X}{0}{sg})$ for $X^0\in \cM_{\per}(E)$. 

We have the surjection 
$\Per(E)\ni v\to \overline{\Orb^+(v)}\in \cM_{\per}(E)$. 
We first see this surjection carefully. 
For $v\in E^0$, 
let us denote by $[v]\subset E^0$ 
the equivalence class of $v$ with respect to 
the equivalence relation on $E^0$ defined so that 
$v$ and $v'$ are equivalent 
if and only if $v'\in \Orb^+(v)$ and $v\in \Orb^+(v')$, 
which is equivalent to $\Orb^+(v)=\Orb^+(v')$. 
Thus $[v]$ is the union of $\{v\}$ and 
the set of vertices which lie on a loop 
whose base point is $v$. 

\begin{lemma}
For $v,v'\in \Per(E)$, 
$\overline{\Orb^+(v)}=\overline{\Orb^+(v')}$ 
if and only if $[v]=[v']$. 
\end{lemma}

\begin{proof}
It is clear that if $[v]=[v']$ then 
$\overline{\Orb^+(v)}=\overline{\Orb^+(v')}$. 
Conversely, suppose that $v,v'\in \Per(E)$ satisfy 
$\overline{\Orb^+(v)}=\overline{\Orb^+(v')}$. 
By definition, 
$\{v\}$ is open in $\overline{\Orb^+(v)}$. 
Hence $\{v\}$ is open in $\overline{\Orb^+(v')}$. 
This implies $v\in \Orb^+(v')$. 
Similarly we have $v'\in \Orb^+(v)$. 
Thus $[v]=[v']$. 
\end{proof}

For $v\in \Per(E)$, 
we have $[v]=\{d(l_1), d(l_2), \ldots, d(l_n)\}$ 
where $l=(l_1,l_2,\ldots,l_n)$ is the unique simple loop 
whose base point is $v$. 
Hence the subset $\Per(E)\subset E^0$ is closed under 
the equivalence relation above. 
Let $[\Per(E)]$ be the set of equivalence classes in $\Per(E)$ 
of the equivalence relation above. 
For a positive integer $n$, 
the map $\Per_n(E)\ni v\mapsto [v]\in [\Per(E)]$ is $n:1$. 
By the lemma above, 
we get the following. 

\begin{corollary}
The map $[\Per(E)]\ni [v]\to \overline{\Orb^+(v)}\in \cM_{\per}(E)$ 
is a well-defined bijection whose inverse is given by 
$\cM_{\per}(E)\ni X^0\to V\in [\Per(E)]$ where 
$$V=\{v\in \Per(E)\mid \overline{\Orb^+(v)}=X^0\}.$$ 
\end{corollary}

Take $X^0\in \cM_{\per}(E)$, 
and set $V=\{v\in \Per(E)\mid \overline{\Orb^+(v)}=X^0\}\in [\Per(E)]$. 
This subset $V$ is the same as the one considered 
in the beginning of Section~\ref{SecMin}. 
Hence by Corollary~\ref{CorSHV}, 
$S(V)$ is an open hereditary and saturated subset of $X^0$ 
contained in $\s{X}{0}{rg}$. 
Define a topological graph $X=(X^0,X^1,d_X,r_X)$ 
so that $X^1=d^{-1}(X^0)$ and $d_X,r_X$ are the restrictions of $d,r$. 
We can define $T^0\colon C_0(X^0)\to \cO(E)/I_{\rho_{X^0}}$ 
and $T^1\colon C_{d_X}(X^1)\to \cO(E)/I_{\rho_{X^0}}$ 
by $T^0(f|_{X^0})=\omega(t^0(f))$ 
and $T^1(\xi|_{X^1})=\omega(t^1(\xi))$ 
for $f\in C_0(E^0)$ and $\xi\in C_d(E^1)$ 
where $\omega\colon \cO(E)\to \cO(E)/I_{\rho_{X^0}}$ 
is the natural quotient map. 
By the proof of Proposition~\ref{E_Ipair}, 
the pair $T=(T^0,T^1)$ induces the isomorphism 
$\cO(X)\to \cO(E)/I_{\rho_{X^0}}$. 
Let us define a subgraph $F=(F^0,F^1,d|_{F^1},r|_{F^1})$ 
of $X$ by $F^0=S(V)$, $F^1=(r_X)^{-1}(F^0)$. 
This subgraph $F$ is nothing but the discrete graph 
considered in Proposition~\ref{GenByLoop}. 
Thus characteristic functions $\delta_v,\delta_e$ 
of $v\in F^0$ and $e\in F^1$ are in $C_0(X^0)$ and $C_{d_X}(X^1)$, 
respectively. 

Choose $v_0\in V$, 
and let $l=(l_1,l_2,\ldots,l_n)$ be the unique simple loop 
whose base point is $v_0$. 
Thus we have $X^0=\overline{\Orb^+(v_0)}$ 
and $V=\{d(l_1), d(l_2), \ldots, d(l_n)\}$. 
Let us set $p_0=T^0(\delta_{v_0})\in \cO(E)/I_{\rho_{X^0}}$ and 
$$u_0=T^1(\delta_{l_1})T^1(\delta_{l_2})\cdots T^1(\delta_{l_n})
\in \cO(E)/I_{\rho_{X^0}}.$$
Note that the element $u_0$ can also be expressed as $T^n(\delta_l)$ 
for using $\delta_l\in C_{d_X}(X^n)$, 
and that we have $u_0^*u_0=u_0u_0^*=p_0$ 
(see the proof of Proposition~\ref{FullHered}). 

\begin{definition}\label{PVw}
For each $w\in \T$, 
we define an ideal $P_{V,w}$ of $\cO(E)$ 
such that $P_{V,w}$ is the inverse image of 
the ideal generated by $u_0-wp_0\in \cO(E)/I_{\rho_{X^0}}$ 
by the natural quotient map $\cO(E)\to \cO(E)/I_{\rho_{X^0}}$. 
\end{definition}

Although $p_0,u_0$ depend on the choice of $v_0\in V$, 
the ideal $P_{V,w}$ does not depend. 
In fact, for $k=1,2,\ldots,n$ 
the elements $p_k$ and $u_k$ defined from $d(l_k)\in V$ as above 
satisfy 
$$p_k=v_k^*p_0v_k,\quad u_k=v_k^*u_0v_k,\quad p_0=v_kp_kv_k^*,\quad 
\text{and}\quad u_0=v_ku_kv_k^*$$ 
for 
$v_k=T^1(\delta_{l_1})T^1(\delta_{l_2})\cdots T^1(\delta_{l_k})
\in \cO(E)/I_{\rho_{X^0}}$. 
This justifies the notation $P_{V,w}$. 

\begin{proposition}
For $z,w\in\T$, 
we have $\beta_{z}(P_{V,w})=P_{V,z^{-n}w}$. 
\end{proposition}

\begin{proof}
This follows from the fact $\beta'_z(p_0)=p_0$ and 
$\beta'_z(u_0)=z^nu_0$ for $z\in\T$ 
where the action $\beta'\colon \T\curvearrowright \cO(E)/I_{\rho_{X^0}}$ 
is induced by the gauge action $\beta$ of $\cO(X)$. 
\end{proof}

We are going to show that $\{P_{V,w}\}_{w\in \T}$ is 
the list of all prime ideals $P$ with $\rho_P=\rho_{X^0}$. 
In the next section, 
we will express the ideal $P_{V,1}$ 
as a kernel of a certain irreducible representation. 
We define an admissible pair $\rho_{X^0}'$ of $E$ 
by $\rho_{X^0}'=(X^0\setminus S(V),\s{X}{0}{sg})$. 
Then we have $\rho_{X^0}'\subset \rho_{X^0}=(X^0,\s{X}{0}{sg})$. 

\begin{lemma}\label{rhoX0}
For an ideal $I$ of $\cO(E)$, 
$\rho_I=\rho_{X^0}$ 
if and only if $I_{\rho_{X^0}}\subset I$ and $I_{\rho_{X^0}'}\not\subset I$. 
\end{lemma}

\begin{proof}
Since $I_{\rho_{X^0}}\subset I$ if and only if $\rho_I\subset \rho_{X^0}$, 
and $I_{\rho_{X^0}'}\not\subset I$ 
if and only if $\rho_I\not\subset \rho_{X^0}'$, 
it suffices to show that 
$\rho_I=\rho_{X^0}$ is equivalent to 
$\rho_I\subset \rho_{X^0}$ and $\rho_I\not\subset \rho_{X^0}'$ 
for an ideal $I$. 
Suppose that the pair $\rho_I=(X_I^0,Z_I)$
satisfies $\rho_I\subset \rho_{X^0}=(X^0,\s{X}{0}{sg})$ 
and $\rho_I\not\subset \rho_{X^0}'=(X^0\setminus S(V),\s{X}{0}{sg})$. 
We have $S(V)\cap X_I^0\neq\emptyset$. 
This implies $V\cap X_I^0\neq\emptyset$ 
by Lemma~\ref{hersat5}. 
Thus $X^0=\overline{\Orb^+(v_0)}\subset X_I^0\subset X^0$. 
Hence $X_I^0=X^0$. 
This implies $\s{X}{0}{sg}=(X_I^0)_{\rs{sg}}\subset Z_I\subset \s{X}{0}{sg}$. 
Thus $Z_I=\s{X}{0}{sg}$. 
We have shown that $\rho_I=\rho_{X^0}$ 
when $\rho_I=(X_I^0,Z_I)$
satisfies $\rho_I\subset \rho_{X^0}=(X^0,\s{X}{0}{sg})$ 
and $\rho_I\not\subset \rho_{X^0}'=(X^0\setminus S(V),\s{X}{0}{sg})$. 
The converse is clear because $\rho_{X^0}\not\subset \rho_{X^0}'$. 
\end{proof}

\begin{lemma}\label{LemPrim1}
The map $P\mapsto (P\cap I_{\rho_{X^0}'})/I_{\rho_{X^0}}$ 
is a bijection from the set of 
prime ideals $P$ of $\cO(E)$ with $\rho_P=\rho_{X^0}$ 
to the set of prime ideals of $I_{\rho_{X^0}'}/I_{\rho_{X^0}}$. 
\end{lemma}

\begin{proof}
It is well-known and routine to check that 
the map $P\mapsto (P\cap I_{\rho_{X^0}'})/I_{\rho_{X^0}}$ 
is a bijection from the set of 
prime ideals $P$ of $\cO(E)$ with 
$I_{\rho_{X^0}}\subset P$ and $I_{\rho_{X^0}'}\not\subset P$ 
to the set of prime ideals of $I_{\rho_{X^0}'}/I_{\rho_{X^0}}$ 
(see \cite[Proposition A.27]{RW} 
for the analogous statement for primitive ideals). 
Thus the conclusion follows from Lemma~\ref{rhoX0}. 
\end{proof}

\begin{lemma}\label{FullHered}
The $C^*$-subalgebra $C^*(u_0)$ generated by $u_0$ 
is a hereditary and full subalgebra 
of the ideal $I_{\rho_{X^0}'}/I_{\rho_{X^0}}$. 
\end{lemma}

\begin{proof}
If we identify $\cO(X)$ and $\cO(E)/I_{\rho_{X^0}}$ 
by the isomorphism induced by the pair $T=(T^0,T^1)$, 
then the admissible pair of $X$ corresponding to 
the ideal $I_{\rho_{X^0}'}/I_{\rho_{X^0}}$ of $\cO(E)/I_{\rho_{X^0}}$ 
is $(X^0\setminus S(V),\s{X}{0}{sg})$. 
Since $\s{X}{0}{sg}\cap S(V)=\emptyset$, 
Proposition~\ref{IrV1} shows that 
$I_{\rho_{X^0}'}/I_{\rho_{X^0}}$ 
is the ideal generated by 
$T^0(C_0(F^0))\subset \cO(E)/I_{\rho_{X^0}}$. 
By Proposition~\ref{IrV2}, 
the $C^*$-subalgebra $A$ of $\cO(E)/I_{\rho_{X^0}}$ 
generated by $T^0(C_0(F^0))$ and $T^1(C_d(F^1))$ 
is a hereditary and full subalgebra 
of the ideal $I_{\rho_{X^0}'}/I_{\rho_{X^0}}$, 
and is naturally isomorphic to $\cO(F)$. 
By Proposition~\ref{C(T,K)} and its proof, 
there exists an isomorphism 
from $A\cong \cO(F)$ to $C(\T)\otimes K(\ell^2(F^\infty))$ 
which sends $p_0$ and $u_0$ to $1\otimes e_{l^\infty,l^\infty}$ 
and $w\otimes e_{l^\infty,l^\infty}$ respectively, 
where $e_{l^\infty,l^\infty}$ is a minimal projection of 
$K(\ell^2(F^\infty))$ 
and $w$ is the generating unitary of $C(\T)$. 
Thus the $C^*$-subalgebra $C^*(u_0)$ of $A$ 
generated by $u_0$ is sent onto $C(\T)\otimes e_{l^\infty,l^\infty}$ 
by the isomorphism $A\cong C(\T)\otimes K(\ell^2(F^\infty))$ 
described in the proof of Proposition~\ref{C(T,K)}. 
Hence $C^*(u_0)$ is a hereditary and full subalgebra of $A$, 
and hence of $I_{\rho_{X^0}'}/I_{\rho_{X^0}}$. 
\end{proof}

From this lemma, 
we see that $P_{V,w}\subset I_{\rho_{X^0}'}$ holds 
for all $w\in\T$. 
By Lemma~\ref{FullHered}, 
the map $P\mapsto P\cap C^*(u_0)$ 
is a bijection from the set of prime ideals $P$ 
of $I_{\rho_{X^0}'}/I_{\rho_{X^0}}$ 
to the set of prime ideals of $C^*(u_0)$. 
As seen in the proof of Lemma~\ref{FullHered}, 
the $C^*$-algebra $C^*(u_0)$ is isomorphic to $C(\T)$, 
and hence the set of prime ideals of $C^*(u_0)$ 
are $\{P_w\}_{w\in\T}$ where $P_w\subset C^*(u_0)$ 
is the ideal of $C^*(u_0)$ 
generated by $u_0-wp_0\in C^*(u_0)$. 
Combining these facts with Lemma~\ref{LemPrim1}, 
we get the following. 

\begin{proposition}\label{QDef}
Let $X^0\in \cM_{\per}(E)$, 
and set $V=\{v\in \Per(E)\mid \overline{\Orb^+(v)}=X^0\}$. 
Then 
the set of all prime ideals $P$ of $\cO(E)$ satisfying $\rho_{P}=\rho_{X^0}$ 
is parameterized by $\T$ as $\{P_{V,w}\}_{w\in\T}$ 
such that $u_0-wp_0\in P_{V,w}/I_{\rho_{X^0}}$ for $w\in\T$. 
\end{proposition}

From the analysis above, 
we get the following theorem. 

\begin{theorem}\label{ThmPrime}
The set of all prime ideals of $\cO(E)$ is 
the union of the following three disjoint sets; 
\benu
\item $\{P_v\mid v\in BV(E)\}$, 
\item $\{P_{X^0}\mid X^0\in \cM_{\aper}(E)\}$, 
\item $\{P_{V,w}\mid V\in[\Per(E)], w\in\T\}$. 
\eenu
The prime ideals in {\rm (i)} and {\rm (ii)} are gauge-invariant, 
and for $v\in \Per_n(E)\subset \Per(E)$ and $w\in\T$ we have 
$\{z\in\T\mid \beta_z(P_{[v],w})=P_{[v],w}\}=\{z\in\T\mid z^n=1\}$. 
\end{theorem}

\section{Irreducible representations and primitive ideals}\label{SecIrrep}

An ideal of a $C^*$-algebra is said to be {\em primitive} 
if it is a kernel of some irreducible representation. 
Every primitive ideal is prime (\cite[Proposition A.17 (b)]{RW}), 
and the converse is true when the $C^*$-algebra is separable 
(\cite[Proposition A.49]{RW}). 
In this section, 
we try to list all primitive ideals of $\cO(E)$. 

We define a subset $\dcM(E)\subset \cM(E)$ by 
$$\dcM(E)=\big\{\overline{\Orb(v,e)}\ \big|\ 
\text{$v\in E^0$ and $e$ is a negative orbit of $v$}\big\}.$$
Then we have $\cM_{\per}(E)\subset \dcM(E)$. 
We define $\dcM_{\aper}(E)\subset \cM_{\aper}(E)$ by 
$$\dcM_{\aper}(E)=\dcM(E)\cap \cM_{\aper}(E)
=\dcM(E)\setminus \cM_{\per}(E).$$ 
The following is the main theorem of this section. 

\begin{theorem}\label{ThmPrimitive}
The following ideals of $\cO(E)$ are primitive; 
\begin{itemize}
\item[(i)\phantom{'}] $\{P_v\mid v\in BV(E)\}$, 
\item[(ii)'] $\{P_{X^0}\mid X^0\in \dcM_{\aper}(E)\}$, 
\item[(iii)\phantom{'}] $\{P_{V,w}\mid V\in[\Per(E)], w\in\T\}$. 
\end{itemize}
\end{theorem}

\begin{remark}
The author was not able to determine all primitive ideals. 
To determine all primitive ideals, 
it suffices to determine the subset $\tcM_{\aper}(E)\subset \cM_{\aper}(E)$ 
defined by 
$$\tcM_{\aper}(E)=\big\{X^0\in \cM_{\aper}(E)\ \big|\ 
\text{$P_{X^0}$ is a primitive ideal}\big\}$$
by Theorem~\ref{ThmPrime} and Theorem~\ref{ThmPrimitive}. 
The theorem above implies that we have 
$\dcM_{\aper}(E)\subset \tcM_{\aper}(E)\subset \cM(E)$. 
In Section~\ref{SecPrimitive}, 
we will see that these two inclusions can be proper. 
The author does not know how to describe $\tcM_{\aper}(E)$ 
in terms of the topological graph $E$. 
\end{remark}

The $C^*$-algebra $\cO(E)$ is separable 
if and only if both $E^0$ and $E^1$ are second countable 
(\cite[Proposition 6.3]{Ka1}). 
Hence in this case every prime ideal of $\cO(E)$ is primitive. 
The following generalizes this fact slightly. 

\begin{corollary}
When $E^0$ is second countable, 
every prime ideal of $\cO(E)$ is primitive. 
\end{corollary}

\begin{proof}
By Lemma~\ref{Orbit=MH2}, 
we have $\dcM_{\aper}(E)=\cM_{\aper}(E)$ 
when $E^0$ is second countable. 
Hence the conclusion follows from 
Theorem~\ref{ThmPrime} and Theorem~\ref{ThmPrimitive}. 
\end{proof}

To prove Theorem~\ref{ThmPrimitive}, 
we need the following lemma. 

\begin{lemma}\label{LemOrb}
For $X^0\in \dcM(E)$, 
either $X^0=\overline{\Orb(r(l),l)}$ for an infinite path $l\in E^\infty$, 
or $X^0=\overline{\Orb^+(v_0)}$ for $v_0\in\s{E}{0}{sg}\setminus BV(E)$. 
\end{lemma}

\begin{proof}
Take $v\in E^0$ and a negative orbit $e$ of $v$ 
such that $X^0=\overline{\Orb(v,e)}$. 
When $e\in E^\infty$ we need to do nothing. 
Suppose $e\in E^*$. 
Then we have $X^0=\overline{\Orb^+(v_0)}$ 
where $v_0=d(e)\in \s{E}{0}{sg}$. 
When $v_0\notin BV(E)$, we are done. 
When $v_0\in BV(E)$ then $v_0\in \s{X}{0}{rg}$. 
Hence we can find $l_1\in d^{-1}(X^0)$ with $r(l_1)=v_0$ 
by Lemma~\ref{RgPr}. 
Set $v_1=d(l_1)\in X^0$. 
Then we have $X^0=\overline{\Orb^+(v_1)}$. 
When $v_1\in \s{E}{0}{sg}\setminus BV(E)$, we are done. 
Otherwise, we have $v_1\in \s{X}{0}{rg}$. 
Hence we can find $l_2\in d^{-1}(X^0)$ with $r(l_2)=v_1$. 
Then we have $X^0=\overline{\Orb^+(v_2)}$ 
where $v_2=d(l_2)\in X^0$. 
By repeating this argument, 
either we can find $v_n\in \s{E}{0}{sg}\setminus BV(E)$ 
with $X^0=\overline{\Orb^+(v_n)}$ 
or we get $l=(l_1,l_2,\ldots)\in E^\infty$ 
such that $X^0=\overline{\Orb^+(r(l),l)}$. 
We are done. 
\end{proof}

By this lemma, Theorem~\ref{ThmPrimitive} 
follows from the next two propositions. 

\begin{proposition}\label{PropIrrepFin}
For $v_0\in \s{E}{0}{sg}$, 
there exists an irreducible representation 
$\psi_{v_0}\colon\cO(E)\to B(H_{v_0})$ such that 
$\rho_{\ker\psi_{v_0}}=
(X^0,\s{X}{0}{sg}\cup\{v_0\})$ 
where $X^0=\overline{\Orb^+(v_0)}$. 
\end{proposition}

\begin{proposition}\label{PropIrrepInf}
For $l\in E^\infty$, 
there exists an irreducible representation 
$\psi_l\colon\cO(E)\to B(H_l)$ such that 
$\rho_{\ker\psi_l}=(X^0,\s{X}{0}{sg})$ 
where $X^0=\overline{\Orb(r(l),l)}$. 
\end{proposition}

\begin{proof}[Proof of Theorem~\ref{ThmPrimitive}]
For $v_0\in BV(E)\subset \s{E}{0}{sg}$, 
the ideal $\ker\psi_{v_0}$ in 
Proposition~\ref{PropIrrepFin} 
coincides with $P_{v_0}$ by Proposition~\ref{Pxuniq}. 
Hence $P_{v_0}$ is primitive. 

For $v_0\in\s{E}{0}{sg}\setminus BV(E)$, 
we have $(X^0,\s{X}{0}{sg}\cup\{v_0\})=\rho_{X^0}$ 
where $X^0=\overline{\Orb^+(v_0)}$. 
Hence in a similar way as above 
using Proposition~\ref{Pxuniq} 
with the help of 
Lemma~\ref{LemOrb}, Proposition~\ref{PropIrrepFin} 
and Proposition~\ref{PropIrrepInf}, 
we can prove that $P_{X^0}$ is primitive 
for $X^0\in \dcM_{\aper}(E)$. 

Let us take $v_0\in\Per_n(E)\subset\Per(E)$, 
and set $X^0=\overline{\Orb^+(v_0)}$ and $V=[v_0]$. 
Let $l\in E^\infty$ be the unique negative orbit of $v_0$. 
Then we have $X^0=\overline{\Orb(r(l),l)}$. 
Hence by Proposition~\ref{PropIrrepInf}, 
we get a primitive ideal $P$ with $\rho_P=\rho_{X^0}$. 
Since $P$ is prime, 
$P=P_{V,w_0}$ for some $w_0\in\T$ 
by Proposition~\ref{QDef}. 
For any $w\in\T$, 
$P_{V,w}=\beta_{z}(P)$ 
with $z\in\T$ satisfying $w=z^{-n}w_0$. 
Hence $P_{V,w}$ is primitive for all $w\in\T$. 
This completes the proof of Theorem~\ref{ThmPrimitive} 
modulo the proofs of Proposition~\ref{PropIrrepFin} and 
Proposition~\ref{PropIrrepInf}. 
\end{proof}

We will prove Proposition~\ref{PropIrrepFin} and 
Proposition~\ref{PropIrrepInf} 
to finish the proof of Theorem~\ref{ThmPrimitive}. 

Take $v_0\in\s{E}{0}{sg}$. 
We define a set $\Lambda_{v_0}$ 
by $\Lambda_{v_0}=\{\lambda\in E^*\mid d(\lambda)=v_0\}$. 
For $\lambda\in \Lambda_{v_0}$ with $\lambda\in E^k$, 
we set its length by $|\lambda|=k$. 
We define $E^1\drtimes \Lambda_{v_0}$ by 
$$E^1\drtimes \Lambda_{v_0}=
\big\{(e,\lambda)\in E^1\times \Lambda_{v_0}\ \big|\ 
d(e)=r(\lambda)\big\}.$$ 
For $(e,\lambda)\in E^1\drtimes \Lambda_{v_0}$, 
we define $e\lambda\in \Lambda_{v_0}$ by 
$ev_0=e$ and $e\lambda=(e,e_1,\ldots,e_k)$ 
for $\lambda=(e_1,\ldots,e_k)\in \Lambda_{v_0}$. 
Then we have $|e\lambda|=|\lambda|+1$ and the map 
$$E^1\drtimes \Lambda_{v_0}\ni (e,\lambda)\mapsto 
e\lambda \in \Lambda_{v_0}\setminus\{v_0\}$$
is a bijection. 
Let $H_{v_0}$ be the Hilbert space whose complete orthonormal system 
is given by $\{\delta_{\lambda}\}_{\lambda\in\Lambda_{v_0}}$. 

\begin{definition}
We define a $*$-homomorphism 
$T_{v_0}^0\colon C_0(E^0)\to B(H_{v_0})$ and 
a linear map 
$T_{v_0}^1\colon C_d(E^1)\to B(H_{v_0})$ 
by 
$$T_{v_0}^0(f)\delta_\lambda=f(r(\lambda))\delta_\lambda,\quad 
T_{v_0}^1(\xi)\delta_\lambda
=\sum_{e\in d^{-1}(r(\lambda))}\xi(e)\delta_{e\lambda}$$ 
for $f\in C_0(E^0)$, $\xi\in C_d(E^1)$ and $\lambda\in\Lambda_{v_0}$. 
\end{definition}

It is not difficult to see that $T_{v_0}^1$ is 
a well defined norm-decreasing linear map. 
We will show that $T_{v_0}=(T_{v_0}^0,T_{v_0}^1)$ is 
a Cuntz-Krieger $E$-pair. 

\begin{lemma}
For $\xi\in C_d(E^1)$ and $(e,\lambda)\in E^1\drtimes \Lambda_{v_0}$, 
we have 
$$T_{v_0}^1(\xi)^*\delta_{v_0}=0,\quad 
T_{v_0}^1(\xi)^*\delta_{e\lambda}=\overline{\xi(e)}\delta_\lambda.$$ 
\end{lemma}

\begin{proof}
Straightforward. 
\end{proof}

\begin{lemma}\label{LemT1fin}
For $\xi,\eta\in C_d(E^1)$, 
we have $T_{v_0}^1(\xi)^*T_{v_0}^1(\eta)=T_{v_0}^0(\ip{\xi}{\eta})$. 
\end{lemma}

\begin{proof}
For $\lambda\in\Lambda_{v_0}$, 
we have 
\begin{align*}
T_{v_0}^1(\xi)^*T_{v_0}^1(\eta)\delta_{\lambda}
&=T_{v_0}^1(\xi)^*
\bigg(\sum_{e\in d^{-1}(r(\lambda))}\eta(e)\delta_{e\lambda}\bigg)\\
&=\sum_{e\in d^{-1}(r(\lambda))}\overline{\xi(e)}\eta(e)\delta_{\lambda}\\
&=\ip{\xi}{\eta}(r(\lambda))\delta_{\lambda}\\
&=T_{v_0}^0(\ip{\xi}{\eta})\delta_{\lambda}.
\end{align*}
This shows $T_{v_0}^1(\xi)^*T_{v_0}^1(\eta)=T_{v_0}^0(\ip{\xi}{\eta})$. 
\end{proof}

\begin{lemma}\label{LemT2fin}
For $f\in C_0(E^0)$ and $\xi\in C_d(E^1)$, 
we have $T_{v_0}^0(f)T_{v_0}^1(\xi)=T_{v_0}^1(\pi_r(f)\xi)$. 
\end{lemma}

\begin{proof}
This follows from the computation 
\begin{align*}
T_{v_0}^0(f)T_{v_0}^1(\xi)\delta_{\lambda}
&=T_{v_0}^0(f)\bigg(\sum_{e\in d^{-1}(r(\lambda))}\xi(e)\delta_{e\lambda}\bigg)\\
&=\sum_{e\in d^{-1}(r(\lambda))}f(r(e\lambda))\xi(e)\delta_{e\lambda}\\
&=\sum_{e\in d^{-1}(r(\lambda))}(\pi_r(f)\xi)(e)\delta_{e\lambda}\\
&=T_{v_0}^1(\pi_r(f)\xi)\delta_{\lambda}, 
\end{align*}
for $\lambda\in \Lambda_{v_0}$. 
\end{proof}

By Lemma~\ref{LemT1fin} and Lemma~\ref{LemT2fin}, 
the pair $T_{v_0}=(T_{v_0}^0,T_{v_0}^1)$ is a Toeplitz $E$-pair. 
We will show that it is a Cuntz-Krieger $E$-pair. 
Let $\varPhi_{v_0}\colon \cK(C_d(E^1))\to B(H_{v_0})$ 
be the $*$-homomorphism defined 
by $\varPhi_{v_0}(\theta_{\xi,\eta})=T_{v_0}^1(\xi)T_{v_0}^1(\eta)^*$ 
for $\xi,\eta\in C_d(E^1)$. 
Recall that the left action 
$\pi_r\colon C_0(E^0)\to \cL(C_d(E^1))$ 
is defined by $\pi_r(f)=\pi(f\circ r)$ for $f\in C_0(E^0)$, 
where $\pi\colon C_0(E^1)\to \cK(C_d(E^1))$ 
is defined by $(\pi(F)\xi)(e)=F(e)\xi(e)$ 
for $F\in C_0(E^1)$, $\xi\in C_d(E^1)$ and $e\in E^1$. 

\begin{lemma}\label{LemT3fin}
For $F\in C_0(E^1)$ and $(e,\lambda)\in E^1\drtimes \Lambda_{v_0}$, 
we have 
$$\varPhi_{v_0}(\pi(F))\delta_{v_0}=0,\quad 
\varPhi_{v_0}(\pi(F))\delta_{e\lambda}=F(e)\delta_{e\lambda}.$$
\end{lemma}

\begin{proof}
Let us take $\xi,\eta\in C_d(E^1)$ 
such that $\xi(e)\overline{\eta}(e')=0$ for $e,e'\in E^1$ 
with $e\neq e'$ and $d(e)=d(e')$. 
We set $F=\xi\overline{\eta}\in C_0(E^1)$. 
Since the linear span of such $F$ is dense 
in $C_0(E^1)$ by \cite[Lemma 1.16]{Ka1}, 
it suffices to show the equalities for this $F$. 
We have $\pi(F)=\theta_{\xi,\eta}$ by \cite[Lemma 1.15]{Ka1}. 
Hence we have 
$$\varPhi_{v_0}(\pi(F))\delta_{v_0}
=T_{v_0}^1(\xi)T_{v_0}^1(\eta)^*\delta_{v_0}=0,$$
and 
\begin{align*}
\varPhi_{v_0}(\pi(F))\delta_{e\lambda}
&=T_{v_0}^1(\xi)T_{v_0}^1(\eta)^*\delta_{e\lambda}\\
&=T_{v_0}^1(\xi)\overline{\eta(e)}\delta_{\lambda}\\
&=\sum_{e'\in d^{-1}(r(\lambda))}\xi(e')\overline{\eta(e)}
\delta_{e'\lambda}\\
&=\xi(e)\overline{\eta(e)}\delta_{e\lambda}\\
&=F(e)\delta_{e\lambda}.
\end{align*}
for $(e,\lambda)\in E^1\drtimes \Lambda_{v_0}$. 
We are done. 
\end{proof}

\begin{proposition}\label{IrrepFin}
The pair $T_{v_0}=(T_{v_0}^0,T_{v_0}^1)$ is a Cuntz-Krieger $E$-pair. 
\end{proposition}

\begin{proof}
By Lemma~\ref{LemT1fin} and Lemma~\ref{LemT2fin}, 
$T_{v_0}$ is a Toeplitz $E$-pair.
Take $f\in C_0(\s{E}{0}{rg})$. 
We have $T_{v_0}^0(f)\delta_{v_0}=f(v_0)\delta_{v_0}=0$ 
because $v_0\in\s{E}{0}{sg}$. 
We also have 
$$\varPhi_{v_0}(\pi_r(f))\delta_{v_0}
=\varPhi_{v_0}(\pi(f\circ r))\delta_{v_0}=0$$
by Lemma~\ref{LemT3fin} because $f\circ r\in C_0(E^1)$. 
Lemma~\ref{LemT3fin} also gives 
\begin{align*}
\varPhi_{v_0}(\pi_r(f))\delta_{e\lambda}
&=\varPhi_{v_0}(\pi(f\circ r))\delta_{e\lambda}\\
&=f(r(e))\delta_{e\lambda}\\
&=T_{v_0}^0(f)\delta_{e\lambda}
\end{align*}
for $(e,\lambda)\in E^1\drtimes \Lambda_{v_0}$. 
Thus we have $\varPhi_{v_0}(\pi_r(f))=T_{v_0}^0(f)$. 
This shows that the pair 
$T_{v_0}=(T_{v_0}^0,T_{v_0}^1)$ is a Cuntz-Krieger $E$-pair. 
\end{proof}

By Proposition~\ref{IrrepFin}, 
we get a representation $\psi_{v_0}\colon\cO(E)\to B(H_{v_0})$ 
such that $T_{v_0}^i=\psi_{v_0}\circ t^i$ for $i=0,1$. 

\begin{proposition}
The representation $\psi_{v_0}$ is irreducible. 
\end{proposition}

\begin{proof}
We will prove that the weak closure of $\psi_{v_0}(\cO(E))$ 
is $B(H_{v_0})$. 
For an approximate unit $\{f_\nu\}$ of $C_0(E^0)$, 
the net $\{T_{v_0}^0(f_\nu)\}$ converges 
to the identity $1\in B(H_{v_0})$ weakly. 
For an approximate unit $\{F_\nu\}$ of $C_0(E^1)$, 
the net $\{\varPhi_{v_0}(\pi(F_\nu))\}$ converges 
to the projection onto the orthogonal complement of $\C\delta_{v_0}$ weakly. 
Hence the rank-one projection $p_{v_0}\in B(H_{v_0})$ onto $\C\delta_{v_0}$ 
is in the weak closure of $\psi_{v_0}(\cO(E))$. 
For $\lambda\in\Lambda_{v_0}$ with $|\lambda|=1$, 
we can find $\xi\in C_d(E^1)$ with $\xi(\lambda)=1$ 
and $\xi(e)=0$ for $e\in d^{-1}(v_0)$ with $e\neq\lambda$. 
Then $T_{v_0}^1(\xi)\delta_{v_0}=\delta_{\lambda}$. 
Similarly, for each $(e,\lambda)\in E^1\drtimes \Lambda_{v_0}$, 
we can find $\xi\in C_d(E^1)$ 
such that $T_{v_0}^1(\xi)\delta_{\lambda}=\delta_{e\lambda}$. 
Hence for every $\lambda\in\Lambda_{v_0}$ 
with $k=|\lambda|\geq 1$, 
there exists $x_\lambda=T_{v_0}^1(\xi_1)T_{v_0}^1(\xi_2)
\ldots T_{v_0}^1(\xi_k)\in \psi_{v_0}(\cO(E))$ 
such that $x_\lambda\delta_{v_0}=\delta_{\lambda}$. 
Let us set $u_{v_0}=p_{v_0}$ and $u_\lambda=x_\lambda p_{v_0}$ 
for $\lambda\in\Lambda_{v_0}\setminus\{v_0\}$, 
which are in the weak closure of $\psi_{v_0}(\cO(E))$. 
Since the von Neumann algebra generated 
by $\{u_\lambda\}_{\lambda\in\Lambda_{v_0}}$ is $B(H_{v_0})$, 
the weak closure of $\psi_{v_0}(\cO(E))$ is $B(H_{v_0})$. 
We are done. 
\end{proof}

\begin{proof}[Proof of Proposition~\ref{PropIrrepFin}]
To finish the proof of Proposition~\ref{PropIrrepFin}, 
it remains to check 
$\rho_{\ker\psi_{v_0}}=(X^0,\s{X}{0}{sg}\cup\{v_0\})$ 
where $X^0=\overline{\Orb^+(v_0)}$. 

From $r(\Lambda_{v_0})=\Orb^+(v_0)$, 
we have $\ker T_{v_0}^0
=C_0(E^0\setminus\overline{\Orb^+(v_0)})=C_0(E^0\setminus X^0)$. 
This shows $X_{\ker\psi_{v_0}}^0=X^0$. 
By Proposition~\ref{rhoI}, 
we have $\s{X}{0}{sg}\subset Z_{\ker\psi_{v_0}}$. 
By the same way to the proof of Lemma~\ref{LemT3fin}, 
we have $\varPhi_{v_0}(k)\delta_{v_0}=0$ for all $k\in\cK(C_d(E^1))$. 
Hence $f(v_0)=0$ for all $f\in C_0(E^0)$ 
with $T_{v_0}^0(f)\in \varPhi_{v_0}\big(\cK(C_d(E^1))\big)$. 
This shows $v_0\in Z_{\ker\psi_{v_0}}$. 
To prove the other inclusion 
$Z_{\ker\psi_{v_0}}\subset \s{X}{0}{sg}\cup\{v_0\}$, 
take 
$$f\in C_0\big(E^0\setminus (\s{X}{0}{sg}\cup\{v_0\})\big),$$ 
and we will prove $T_{v_0}^0(f)\in\varPhi_{v_0}\big(\cK(C_d(E^1))\big)$. 
We have $(f\circ r)|_{X^1}\in C_0(X^1)$ 
because $f|_{X^0}\in C_0(\s{X}{0}{rg})$. 
Hence there exists $F\in C_0(E^1)$ 
such that $F(e)=f(r(e))$ for $e\in X^1$. 
For $(e,\lambda)\in E^1\drtimes \Lambda_{v_0}$, 
we have 
$$T_{v_0}^0(f)\delta_{e\lambda}=f(r(e))\delta_{e\lambda}
=F(e)\delta_{e\lambda}=\varPhi_{v_0}(\pi(F))\delta_{e\lambda}$$
because $e\in X^1$.
Since $f(v_0)=0$, 
we have $T_{v_0}^0(f)\delta_{v_0}=0=\varPhi_{v_0}(\pi(F))\delta_{v_0}$. 
Hence $T_{v_0}^0(f)=\varPhi_{v_0}(\pi(F))
\in\varPhi_{v_0}\big(\cK(C_d(E^1))\big)$. 
Therefore we get 
$\rho_{\ker\psi_{v_0}}=(X^0,\s{X}{0}{sg}\cup\{v_0\})$. 
\end{proof}

Next we will prove Proposition~\ref{PropIrrepInf}. 
The method is almost identical to the proof of 
Proposition~\ref{PropIrrepFin}, 
hence we just sketch it. 
In fact, the proof of Proposition~\ref{PropIrrepInf} 
is easier than the one of Proposition~\ref{PropIrrepFin} 
because in the case of Proposition~\ref{PropIrrepFin} 
there exists a special element $v_0\in \Lambda_{v_0}$ 
which cannot be expressed as $e\lambda$, 
but no such elements appear 
in the case of Proposition~\ref{PropIrrepInf}. 

Let us take $l=(l_1,l_2,\ldots,l_n,\ldots)\in E^\infty$. 
We denote by $\Lambda_{l}$ 
the set of 
$$\lambda=(\lambda_1,\lambda_2,\ldots,\lambda_n,\ldots)\in E^\infty$$ 
satisfying that there exist $k\in\Z$ and $N\in\N$ 
such that $\lambda_{n+k}=l_n$ for all $n\geq N$. 
We define $E^1\drtimes \Lambda_{l}$, and the bijective map 
$$E^1\drtimes \Lambda_{l}\ni (e,\lambda)\mapsto 
e\lambda \in \Lambda_{l}$$
in a similar way as before. 
Let $H_{l}$ be the Hilbert space whose complete orthonormal system 
is given by $\{\delta_{\lambda}\}_{\lambda\in\Lambda_{l}}$. 

\begin{definition}
We define a $*$-homomorphism 
$T_l^0\colon C_0(E^0)\to B(H_{l})$ and 
a linear map 
$T_l^1\colon C_d(E^1)\to B(H_{l})$ 
by 
$$T_l^0(f)\delta_\lambda=f(r(\lambda))\delta_\lambda,\quad 
T_l^1(\xi)\delta_\lambda
=\sum_{e\in d^{-1}(r(\lambda))}\xi(e)\delta_{e\lambda}$$ 
for $f\in C_0(E^0)$, $\xi\in C_d(E^1)$ and $\lambda\in\Lambda_{l}$. 
\end{definition}

In a similar way to the proof of Proposition~\ref{IrrepFin}, 
we can show that $T_l=(T_l^0,T_l^1)$ is 
a Cuntz-Krieger $E$-pair. 
Thus we get a representation $\psi_l\colon\cO(E)\to B(H_{l})$ 
such that $T_l^i=\psi_l\circ t^i$ for $i=0,1$. 

\begin{proposition}
The representation $\psi_l$ is irreducible. 
\end{proposition}

\begin{proof}
First note that 
we can define the set $E^n\drtimes \Lambda_{l}$ and the bijective map 
$$E^n\drtimes \Lambda_{l}\ni (e,\lambda)\mapsto 
e\lambda \in \Lambda_{l}$$
for a positive integer $n$. 
For $\xi\in C_d(E^n)$ and $(e,\lambda)\in E^n\drtimes \Lambda_{l}$, 
we have $T_{l}^n(\xi)\delta_{e\lambda}=\xi(e)\delta_{\lambda}$. 
Hence if $\xi\in C_d(E^n)$ satisfies $\xi(e)=1$ 
and $\xi(e')=0$ for $e'\in E^n$ with $d(e')=d(e)$ and $e'\neq e$, 
then $T_l^n(\xi)\delta_{\lambda}=\delta_{e\lambda}$ 
and $T_l^n(\xi)^*\delta_{e\lambda}=\delta_{\lambda}$. 
Hence for each $\lambda\in \Lambda_{l}$ 
we can find $x_\lambda\in \psi_l(\cO(E))$ 
such that $x_\lambda\delta_{l}=\delta_{\lambda}$. 

Let $p_{l}\in B(H_{l})$ be 
the rank-one projection onto $\C\delta_{l}$. 
To prove that $p_l$ is 
in the weak closure of $\psi_l(\cO(E))$, 
it suffices to show that 
for each finite subset $Y\subset \Lambda_{l}\setminus\{l\}$ 
there exists $x\in \psi_l(\cO(E))$ 
such that $x\delta_{l}=\delta_{l}$ 
and $x\delta_{\lambda}=0$ for $\lambda\in Y$. 
Take a finite subset $Y\subset \Lambda_{l}\setminus\{l\}$. 
For a positive integer $n$, 
we define $l_{[1,n]}=(l_1,l_2,\ldots,l_n)\in E^n$ and 
$$Y_{[1,n]}=\{e\in E^n\mid (e,\lambda)\in E^n\drtimes \Lambda_{l} 
\text{ with }e\lambda\in Y\}.$$
Since $l\notin Y$, 
we have $l_{[1,n]}\notin Y_{[1,n]}$ 
for a sufficiently large integer $n$. 
Choose $\xi\in C_d(E^{n})$ such that $\xi(l_{[1,n]})=1$ and 
$\xi(e)=0$ for every $e\in Y_{[1,n]}$. 
Then $x=T_l^n(\xi)T_l^n(\xi)^*\in \psi_l(\cO(E))$ 
satisfies $x\delta_{l}=\delta_{l}$ 
and $x\delta_{\lambda}=0$ for $\lambda\in Y$. 
This shows that $p_l$ is 
in the weak closure of $\psi_l(\cO(E))$. 
Since the von Neumann algebra generated 
by $\{x_\lambda p_l\}_{\lambda\in\Lambda_{l}}$ is $B(H_{l})$, 
the weak closure of $\psi_l(\cO(E))$ is $B(H_{l})$. 
The proof is completed. 
\end{proof}

\begin{proof}[Proof of Proposition~\ref{PropIrrepInf}]
To finish the proof of Proposition~\ref{PropIrrepInf}, 
it remains to check $\rho_{\ker\psi_l}=(X^0,\s{X}{0}{sg})$ 
where $X^0=\overline{\Orb(r(l),l)}$. 
From $r(\Lambda_{l})=\Orb(r(l),l)$, 
we have $\ker T_l^0=C_0(E^0\setminus X^0)$. 
This shows $X_{\ker\psi_l}^0=X^0$. 
The proof of $Z_{\ker\psi_l}=\s{X}{0}{sg}$ 
is the same as the proof of Proposition~\ref{PropIrrepFin}. 
\end{proof}

Thus we have completed the proof of Theorem~\ref{ThmPrimitive}. 
In the proof, we used Proposition~\ref{Pxuniq} 
which depends heavily on the Cuntz-Krieger Uniqueness Theorem 
(Proposition~\ref{CKUT}). 
In the following, 
we give the direct proof of the gauge-invariance 
of the kernels of the irreducible representations we consider, 
so that we can use the Gauge Invariant Uniqueness Theorem 
(\cite[Theorem 4.5]{Ka1}) instead of Proposition~\ref{Pxuniq}. 
Note that the proof of the Gauge Invariant Uniqueness Theorem 
is much shorter and easier than the one of 
the Cuntz-Krieger Uniqueness Theorem. 
We also analyze the primitive ideal which is not gauge-invariant 
in the detail. 

\begin{lemma}\label{GaugeFin}
For $v_0\in\s{E}{0}{sg}$, 
the primitive ideal $\ker\psi_{v_0}$ is gauge-invariant. 
\end{lemma}

\begin{proof}
It suffices to see that the Cuntz-Krieger $E$-pair 
$T_{v_0}=(T_{v_0}^0,T_{v_0}^1)$ 
admits a gauge action. 
For $z\in\T$, 
we define a unitary $u_z\in B(H_{v_0})$ 
by $u_z\delta_{\lambda}=z^{|\lambda|}\delta_{\lambda}$ 
for $\lambda\in \Lambda_{v_0}$. 
Then it is easy to see that the automorphism $\Ad(u_z)$ 
of $B(H_{v_0})$ defined by $\Ad(u_z)(x)=u_zxu_z^*$ for $z\in\T$ 
is a gauge action for $T_{v_0}$. 
We are done. 
\end{proof}

\begin{definition}
An infinite path $l=(l_1,l_2,\ldots,l_n,\ldots)\in E^\infty$ 
is said to be {\em periodic} 
if there exist positive integers $k$ and $N$ 
such that $l_{n+k}=l_n$ for all $n\geq N$. 
Otherwise $l\in E^\infty$ is said to be {\em aperiodic}. 
\end{definition}

\begin{lemma}\label{GaugeInf}
When $l\in E^\infty$ is aperiodic, 
the primitive ideal $\ker\psi_l$ is gauge-invariant. 
\end{lemma}

\begin{proof}
It suffices to see that the Cuntz-Krieger $E$-pair $T_l=(T_l^0,T_l^1)$ 
admits a gauge action. 
When $l$ is aperiodic, 
for each $\lambda\in\Lambda_{l}$, 
an integer $k\in\Z$ 
satisfying $\lambda_{n+k}=l_n$ for large $n$ 
is unique. 
Hence we can write $c_\lambda=k\in\Z$. 
It is easy to see that $c_{e\lambda}=c_\lambda+1$ 
for $(e,\lambda)\in E^1\drtimes \Lambda_{l}$. 
Now we define a unitary $u_z\in B(H_{l})$ 
for $z\in\T$ by $u_z\delta_{\lambda}=z^{c_\lambda}\delta_{\lambda}$ 
for $\lambda\in \Lambda_{l}$. 
Then it is easy to see that the automorphism $\Ad(u_z)$ 
of $B(H_{l})$ defined by $\Ad(u_z)(x)=u_zxu_z^*$ for $z\in\T$ 
is a gauge action for $T_l$. 
We are done. 
\end{proof}

Let $l$ be a periodic infinite path. 
Then there exists a simple loop 
$l'=(e_1,e_2,\ldots,e_n)$ 
such that $l=(e', l',l',\ldots,l',\ldots)$ for some $e'\in E^*$. 
Set $v_0=r(l')\in E^0$. 
Then the closed set $X^0=\overline{\Orb(r(l),l)}$ coincides 
with $\overline{\Orb^+(v_0)}$. 
Let us define $U=\{e_1,e_2,\ldots,e_n\}\subset E^1$. 
For $\lambda=(\lambda_1,\lambda_2,\ldots,\lambda_n,\ldots)\in \Lambda_{l}$, 
we define $|\lambda|\in\N$ 
by $|\lambda|=0$ if $\lambda_k\in U$ for all $k$, 
and $|\lambda|=\max\{k \mid \lambda_k\notin U\}$ 
otherwise. 
For $(e,\lambda)\in E^1\drtimes \Lambda_{l}$, 
we have $|e\lambda|=0$ if and only if $e\in U$ and $|\lambda|=0$. 
Otherwise, we have $|e\lambda|=|\lambda|+1$. 
Note that the number of elements $\lambda\in \Lambda_{l}$ with $|\lambda|=0$ 
is $n$. 

Let $\psi_l'\colon\cO(E)\to B(H_{l})/K(H_{l})$ 
be the composition of the representation $\psi_l\colon\cO(E)\to B(H_{l})$ 
and the natural surjection $\varPi\colon B(H_{l})\to B(H_{l})/K(H_{l})$. 
We have $\ker\psi_l\subset\ker\psi_l'$. 

\begin{proposition}\label{psi'gauge}
The ideal $\ker\psi_l'$ is gauge-invariant. 
\end{proposition}

\begin{proof}
It suffices to see that $(\varPi\circ T_l^0,\varPi\circ T_l^1)$ admits a gauge action. 
For each $z\in\T$, 
define a unitary $u_z\in B(H_{l})$ 
by $u_z\delta_{\lambda}=z^{|\lambda|}\delta_{\lambda}$ 
for $\lambda\in\Lambda_{l}$. 
Then we have $u_zT_l^0(f)u_z^*=T_l^0(f)$ for all $f\in C_0(E^0)$. 
For $\xi\in C_d(E^1)$, 
we have $u_zT_l^1(\xi)u_z^*=zT_l^1(\xi)$ 
on the subspace generated by $\{\delta_{\lambda}\}_{|\lambda|\geq 1}$ 
which is codimension finite. 
Hence if we define an automorphism $\beta'_z$ of $B(H_{l})/K(H_{l})$ 
by $\beta'_z(x)=\varPi(u_z)x\varPi(u_z)^*$ for $z\in\T$, 
then $\beta'$ is a gauge action for $(\varPi\circ T_l^0,\varPi\circ T_l^1)$. 
This completes the proof. 
\end{proof}

\begin{lemma}\label{X0psi'}
We have $X^0_{\ker\psi_l'}\subset X^0$, 
and $v\in X^0\setminus X^0_{\ker\psi_l'}$ 
if and only if 
$\{v\}$ is open in $X^0$ and 
the set $\{\lambda\in\Lambda_l\mid r(\lambda)=v\}$ is finite. 
\end{lemma}

\begin{proof}
Since $\ker\psi_l\subset\ker\psi_l'$, 
we have $X^0_{\ker\psi_l'}\subset X^0_{\ker\psi_l}=X^0$. 
Take $v\in X^0\setminus X^0_{\ker\psi_l'}$. 
Then there exists $f\in C_0(E^0)$ such that 
$t^0(f)\in\ker \psi_l'$ and $f(v)=1$. 
Set $$V=\{v\in X^0\mid f(v)\geq 1/2\}$$ 
which is a neighborhood of $v\in X^0$. 
Since $t^0(f)\in\ker \psi_l'$, we have $T^0(f)=\psi_l(t^0(f))\in K(H_l)$. 
Hence the projection $p=\chi_{[1/2,\infty)}(T^0(f))$ is of finite rank. 
By the definition of $T^0$, 
$p$ is the projection onto the subspace spanned by 
$$\{\delta_{\lambda}\mid f(r(\lambda))\geq 1/2\}
=\{\delta_{\lambda}\mid r(\lambda)\in V\}.$$
Thus $\{\lambda\in\Lambda_l\mid r(\lambda)\in V\}$ is finite. 
This shows that $\{\lambda\in\Lambda_l\mid r(\lambda)=v\}$ is finite. 
Since $\Orb^+(v_0)=\Orb(r(l),l)$ 
is the image of the map $r\colon \Lambda_l\to E^0$, 
$\Orb^+(v_0)\cap V$ is a finite set. 
Since $\Orb^+(v_0)$ is dense in $X^0$, 
$V$ is a finite subset of $\Orb^+(v_0)$. 
This shows that $\{v\}$ is open. 

Conversely suppose that $\{v\}$ is open in $X^0$ and 
the set $\{\lambda\in\Lambda_l\mid r(\lambda)=v\}$ is finite. 
We can find $f\in C_0(E^0)$ 
such that $f(v)=1$ and $f(v')=0$ for $v'\in X^0\setminus\{v\}$. 
Then $T^0(f)$ is a projection onto the subspace spanned by 
the finite set $\{\delta_{\lambda}\mid r(\lambda)=v\}$. 
Hence $t^0(f)\in\ker \psi_l'$. 
This shows $v\notin X^0_{\ker\psi_l'}$. 
We are done. 
\end{proof}

\begin{proposition}
When $v_0 \in\Aper(E)$, 
we have $\ker\psi_l=\ker\psi_l'$ 
and hence the primitive ideal $\ker\psi_l$ is gauge-invariant. 
\end{proposition}

\begin{proof}
By Lemma~\ref{X0psi'}, 
$v_0\in\Aper(E)$ implies $v_0\in X^0_{\ker\psi_l'}$. 
Hence 
$$X^0=\overline{\Orb^+(v_0)}\subset X^0_{\ker\psi_l'}\subset X^0.$$
This shows $X^0_{\ker\psi_l'}=X^0$. 
Hence $\s{X}{0}{sg}=(X^0_{\ker\psi_l'})_{\rs{sg}}\subset Z_{\ker\psi_l'}$. 
Since $\ker\psi_l\subset\ker\psi_l'$, 
we have $Z_{\ker\psi_l'}\subset \s{X}{0}{sg}$. 
Thus $Z_{\ker\psi_l'}=\s{X}{0}{sg}$. 
We have shown that $\rho_{\ker\psi_l'}=\rho_{\ker\psi_l}$. 
By Proposition~\ref{psi'gauge}, 
we have $\ker\psi_l'=I_{\rho_{\ker\psi_l'}}
=I_{\rho_{\ker\psi_l}}\subset \ker\psi_l$. 
Therefore we have $\ker\psi_l=\ker\psi_l'$. 
\end{proof}

\begin{proposition}\label{P_{V,1}}
When $v_0\in\Per(E)$, 
we have $\ker\psi_l=P_{V,1}$ 
and hence the primitive ideal $\ker\psi_l$ is not gauge-invariant. 
\end{proposition}

\begin{proof}
We consider $X^n$ as a closed subset of $E^n$. 
Since $v_0\in\Per(E)$, 
$\{l'\}$ is open in $X^n$. 
Hence there exists $\xi\in C_d(E^n)$ such that $\xi(l')=1$ 
and $\xi(e)=0$ for $e\in X^n\setminus\{l'\}$. 
Note that the image of $t^n(\xi)\in\cO(E)$ in $\cO(E)/I_{\rho_{X^0}}$ 
is the element $u_0$ defined in the previous section 
(see the remarks before Definition~\ref{PVw}). 
We can see that $T_l^n(\xi)\delta_{\lambda}=\delta_{\lambda}$ 
when $|\lambda|=0$, and $T_l^n(\xi)\delta_{\lambda}=0$ 
when $|\lambda|\geq 1$. 
Hence $T_l^n(\xi)\in B(H_l)$ is the projection onto 
the subspace spanned by $\{\delta_{\lambda}\}_{|\lambda|=0}$. 
This shows $t^n(\xi)-t^n(\xi)^*t^n(\xi)\in \ker\psi_l$. 
Thus $\ker\psi_l$ is a prime ideal with $\rho_{\ker\psi_l}=\rho_{X^0}$ and 
$u_0-u_0^*u_0\in \ker\psi_l/I_{\rho_{X^0}}$. 
Therefore by Proposition~\ref{QDef} we have $\ker\psi_l=P_{V,1}$. 
\end{proof}

When $v_0\in\Per(E)$, 
the admissible pair $\rho_{X^0}'$ 
was defined by $\rho_{X^0}'=(X^0\setminus S(V),\s{X}{0}{sg})$ 
where $V=\{d(e_1),d(e_2),\ldots,d(e_n)\}$. 

\begin{proposition}\label{I_{r_{X^0}'}}
When $v_0\in\Per(E)$, 
we have $\rho_{\ker\psi_l'}=\rho_{X^0}'$ 
and hence $\ker\psi_l'=I_{\rho_{X^0}'}$. 
\end{proposition}

\begin{proof}
By Proposition~\ref{SHV} and Lemma~\ref{X0psi'}, 
we have $X^0\setminus X^0_{\ker\psi_l'}=S(V)$. 
Hence $X^0_{\ker\psi_l'}=X^0\setminus S(V)$. 
Since $\ker\psi_l\subset\ker\psi_l'$, 
we have $Z_{\ker\psi_l'}\subset Z_{\ker\psi_l}=\s{X}{0}{sg}$. 
We will show the other inclusion  $\s{X}{0}{sg}\subset Z_{\ker\psi_l'}$. 
To do so, 
it suffices to show $T_l^0(f)\notin K(H_l)+\varPhi_l\big(\cK(C_d(E^1))\big)$ 
for $f\in C_0(E^0)$ with $f(v)\neq 0$ for some $v\in\s{X}{0}{sg}$, 
where $\varPhi_l\colon \cK(C_d(E^1))\to B(H_l)$ is 
defined by $\varPhi_l(\theta_{\xi,\eta})=T_l^1(\xi)T_l^1(\eta)^*$ 
for $\xi,\eta\in C_d(E^1)$. 
For $x\in K(H_l)$ and $\e>0$, 
$$\big\{e\in E^1\ \big|\ 
\|x\delta_{e\lambda}\|>\e, 
(e,\lambda)\in E^1\drtimes \Lambda_{l}\big\}$$
is finite. 
For $\xi,\eta\in C_d(E^1)$ and 
$(e,\lambda)\in E^1\drtimes \Lambda_{l}$, 
we have 
\begin{align*}
\big\|\varPhi_l(\theta_{\xi,\eta})\delta_{e\lambda}\big\|
&=\big\|\big(T_l^1(\xi)T_l^1(\eta)^*\big)\delta_{e\lambda}\big\|\\
&=\big\|T_l^1(\xi)(\overline{\eta(e)}\delta_{\lambda})\big\|\\
&=\big|\overline{\eta(e)}\big|
\bigg\|\sum_{e'\in d^{-1}(r(\lambda))}\xi(e')\delta_{e'\lambda}\bigg\|\\
&=|\eta(e)|\bigg(\sum_{e'\in d^{-1}(d(e))}
\big|\xi(e')\big|^2\bigg)^{1/2}\\
&=\big|(\eta f)(e)\big|
\end{align*}
where $f=\ip{\xi}{\xi}^{1/2}\in C_0(E^0)$. 
Hence 
for all $x\in \varPhi_l\big(\cK(C_d(E^1))\big)$ and $\e>0$, the set 
$$\big\{e\in E^1\ \big|\ 
\|x\delta_{e\lambda}\|>\e, 
(e,\lambda)\in E^1\drtimes \Lambda_{l}\big\}$$
has a compact closure. 
Thus this is true for all $x\in K(H_l)+\varPhi_l\big(\cK(C_d(E^1))\big)$. 

Take $f\in C_0(E^0)$ with $f(v)\neq 0$ for some $v\in\s{X}{0}{sg}$. 
Since $\Orb^+(v_0)$ is dense in $X^0$, 
we have $\s{X}{0}{sce}=\emptyset$. 
Hence $v\in \s{X}{0}{inf}$. 
There exists $\e>0$ such that 
$W=\{w\in E^0\mid |f(w)|>\e\}$ 
is a neighborhood of $v$. 
Since $v\in \s{X}{0}{inf}$, 
the closure of $r^{-1}(W)\cap X^1$ is not compact. 
Since $\|T_l^0(f)\delta_{e\lambda}\|=|f(r(e\lambda))|=|f(r(e))|$, 
we have 
\begin{align*}
\big\{e\in E^1\ &\big|\ 
\|T_l^0(f)\delta_{e\lambda}\|>\e, 
(e,\lambda)\in E^1\drtimes \Lambda_{l}\big\}\\
&=\big\{e\in E^1\ \big|\ 
|f(r(e))|>\e, d(e)\in \Orb^+(v_0)\big\}=r^{-1}(W)\cap 
d^{-1}(\Orb^+(v_0))
\end{align*}
whose closure is not compact. 
This shows $T_l^0(f)\notin K(H_l)+\varPhi_l\big(\cK(C_d(E^1))\big)$. 
Hence we have $Z_{\ker\psi_l'}=\s{X}{0}{sg}$. 
Therefore we get $\rho_{\ker\psi_l'}=\rho_{X^0}'$. 

Since $\ker\psi_l'$ is gauge-invariant by Proposition~\ref{psi'gauge}, 
we have $\ker\psi_l'=I_{\rho_{X^0}'}$. 
\end{proof}

\section{Primitivity of $\cO(E)$}\label{SecPrimitive}

A $C^*$-algebra is said to be {\em primitive} 
if $0$ is a primitive ideal. 
Equivalently, a $C^*$-algebra is primitive 
if and only if it has a faithful irreducible representation. 
On the primitivity of $\cO(E)$, 
we have the following. 

\begin{proposition}\label{PropPrimitive}
For a topological graph $E$, 
consider the following three conditions. 
\benu
\item $E$ is topologically free and $E^0=\overline{\Orb(v,e)}$ 
for some $v\in E^0$ and a negative orbit $e$ of $v$. 
\item the $C^*$-algebra $\cO(E)$ is primitive. 
\item $E$ is topologically free and topologically transitive. 
\eenu
Then we have {\rm (i)}$\Rightarrow${\rm (ii)}$\Rightarrow${\rm (iii)}. 
When $E^0$ is second countable, 
the three conditions are equivalent. 
\end{proposition}

\begin{proof}
The condition (i) is equivalent to $E^0\in\dcM_{\aper}(E)$ 
by Proposition~\ref{NotTopFree}. 
Hence we have (i)$\Rightarrow$(ii) 
by Theorem~\ref{ThmPrimitive}. 
Since a primitive $C^*$-algebra is prime, 
the implication (ii)$\Rightarrow$(iii) follows from 
Theorem~\ref{C*prime}. 
When $E^0$ is second countable, 
(iii) implies (i) by Proposition~\ref{TopTran}. 
\end{proof}

The converses of the two implications 
in Proposition~\ref{PropPrimitive} 
are not true in general 
as we will see the following two examples. 

\begin{example}[A topological graph satisfying (ii) but not (i)]
\label{ExErgode}\ 

Let $\mu$ be the Haar measure on $\T$. 
An irrational rotation $\alpha$ on $\T$ preserves 
the measure $\mu$. 
Hence $\alpha$ induces the automorphism $\bar{\alpha}$ on 
the commutative von Neumann algebra $L^\infty(\T,\mu)$. 
Let $X$ be the spectrum of $L^\infty(\T,\mu)$ 
which is considered as a commutative $C^*$-algebra. 
Thus $X$ is a compact hyperstonean space 
such that $C(X)\cong L^\infty(\T,\mu)$. 
The automorphism $\bar{\alpha}$ of $L^\infty(\T,\mu)$ gives us 
a homeomorphism $\sigma$ on $X$. 
Thus we get a dynamical system 
$\Sigma=(X,\sigma)$. 
Since $\mu$ is non-atomic, 
every orbit of $\Sigma$ is not dense in $X$ 
by \cite[Proposition~1.2~(1)]{T2}. 
Hence the topological graph 
$E_\Sigma=(X,X,\id_X,\sigma)$ does not satisfy (i). 
We will show that $\cO(E_\Sigma)$ is primitive. 

Let us consider the covariant representation $\{\pi,u\}$ 
of $\Sigma=(X,\sigma)$ 
on $L^2(\T,\mu)$ 
such that $\pi\colon C(X)\cong L^\infty(\T,\mu)\to B(L^2(\T,\mu))$ 
is defined by a multiplication, 
and the unitary $u\in B(L^2(\T,\mu))$ is defined from 
the irrational rotation $\alpha$ on $\T$. 
Since the irrational rotation $\alpha$ on $(\T,\mu)$ 
is free, 
the dynamical system $\Sigma=(X,\sigma)$ is topologically free 
by \cite[Proposition~1.2~(3)]{T2}. 
Hence the covariant representation $\{\pi,u\}$ 
gives a faithful representation 
$\psi\colon \cO(E_\Sigma)\to B(L^2(\T,\mu))$. 
Since the irrational rotation $\alpha$ is ergodic, 
$\psi$ is irreducible. 
Thus the $C^*$-algebra $\cO(E_\Sigma)$ is primitive. 
\end{example}

\begin{example}[A topological graph satisfying (iii) but not (ii)]\ 

Let $E=(E^0,E^1,d,r)$ be 
the discrete graph in Example~\ref{NiceExample}. 
Namely $E^0$ is the set of all finite subsets 
of an uncountable set $X$, 
$$E^1=\big\{(x;v)\ \big|\ 
\text{$v\in E^0$ and $x\in v$}\big\},$$
$d((x;v))=v$ and $r((x;v))=v\setminus \{x\}$ for $(x;v)\in E^1$. 
For a positive integer $n$ and $v\in E^0$, 
let $v^{(n)}$ 
be the set of $n$-tuples 
$x=(x_1,x_2,\ldots,x_n)\in v^n$ 
such that $x_k\neq x_l$ for $k\neq l$. 
Note that $v^{(1)}$ is identified with $v$. 
For $v\in E^0$, 
$|v|\in\N$ denotes the number of elements of $v$. 
When $|v|<n$ we have $v^{(n)}=\emptyset$, 
and when $|v|=m\geq n$ 
we have $|v^{(n)}|=m!/(m-n)!$. 
For $\big((x_1;v_1),(x_2;v_2),\ldots,(x_n;v_n)\big)\in E^n$, 
$x=(x_1,x_2,\ldots,x_n)$ is in $v^{(n)}$ 
where $v=v_n\in E^0$. 
Conversely for $v\in E^0$ and $x=(x_1,x_2,\ldots,x_n)\in v^{(n)}$, 
we have $\big((x_1;v_1),(x_2;v_2),\ldots,(x_n;v_n)\big)\in E^n$ 
where $v_k=v\setminus\{x_{k+1},\ldots,x_n\}$ 
for $k=1,2,\ldots,n$. 
By these correspondences, we will identify $E^n$ with the set 
$$\big\{(x;v)\ \big|\ 
\text{$v\in E^0$ and $x\in v^{(n)}$}\big\},$$
For $(x;v)\in E^n$ where $x=(x_1,x_2,\ldots,x_n)\in v^{(n)}$, 
we have $d^n((x;v))=v$ 
and $r^n((x;v))=v\setminus \{x_1,x_2,\ldots,x_n\}$. 
In order to save the notation, 
we set $v^{(0)}=\{\emptyset\}$ for $v\in E^0$ 
and identify $E^0$ with $\{(\emptyset;v)\mid v\in E^0\}$. 
For $v\in E^0$, 
we set $v^{(*)}=\coprod_{n\in\N}v^{(n)}=\coprod_{n=0}^{|v|}v^{(n)}$. 

We put $s_{(x;v)}=t^n(\delta_{(x;v)})\in \cO(E)$ 
for $(x;v)\in E^n$ and $n\in\N$, 
where $\delta_{(x;v)}\in C_d(E^n)$ is 
the characteristic function of $\{(x;v)\}$. 
The linear space 
$$\spa\big\{s_{(x;v)}s_{(y;v)}^*\ \big|\ 
\text{$v\in E^0$ and $x,y\in v^{(*)}$}\big\}$$ 
is dense in the $C^*$-algebra $\cO(E)$. 

\begin{proposition}\label{notprimitive}
The $C^*$-algebra $\cO(E)$ is prime, but not primitive. 
\end{proposition}

\begin{proof}
We had already seen that $E^0$ is a maximal head. 
Hence $E$ is topologically transitive by Proposition~\ref{TopTran}. 
Since $E$ has no loops, $E$ is topologically free. 
By Theorem~\ref{C*prime}, the $C^*$-algebra $\cO(E)$ is prime. 

Take an irreducible representation 
$\psi\colon\cO(E)\to B(H)$, 
and we will show that $\psi$ is not faithful. 
Choose $\xi\in H$ arbitrary. 
For $n\in\N$, 
we set $\Omega_n\subset E^0$ by 
$$\Omega_n=\big\{v\in E^0\ \big|\ 
\text{$\psi(s_{(x;v)}s_{(x;v)}^*)\xi\neq 0$ 
for some $x\in v^{(n)}$}\big\}.$$
For each $n\in\N$, 
$\Omega_n$ is countable 
because $\{s_{(x;v)}s_{(x;v)}^*\}_{v\in E^0,x\in v^{(n)}}$ 
is an orthogonal family of projections. 
Therefore $\Omega=\bigcup_{n=0}^\infty \Omega_n$ is also a countable 
subset of $E^0$. 
Hence we can find $x_0\in X$ such that $x_0\notin v$ 
for all $v\in \Omega$. 
Let $I\subset \cO(E)$ be the closure of 
$$\spa\big\{s_{(x;v)}s_{(y;v)}^* \ \big|\ 
\text{$v\in E^0$ with $x_0\in v$, and 
$x,y\in v^{(*)}$}\big\}.$$
By noting that the set $\{v\in E^0\mid x_0\in v\}$ 
is hereditary, 
we can show that $I$ is an ideal 
(cf.\ the proof of Proposition~\ref{IisGaugeInv}). 
Since $x_0\in v$ implies $v\notin \Omega$, 
we have $\psi(s_{(y;v)}^*)\xi=0$ for $v\in E^0$ with $x_0\in v$. 
Hence $\psi(a)\xi=0$ for all $a\in I$. 
Since $\xi$ is a cyclic vector for the representation $\psi$, 
we have $I\subset \ker\psi$. 
Thus $\psi$ is not injective. 
\end{proof}

\begin{remark}
By the proof of Proposition~\ref{notprimitive}, 
we can see that the $C^*$-algebra $\cO(E)$ does not have 
a faithful cyclic representation. 
This is the obstacle that N. Weaver used in \cite{We}. 
We do not know whether this is the only obstacle 
for prime $C^*$-algebras to become primitive. 
Namely the following problem is still open. 
\end{remark}

\begin{problem}
Is a prime $C^*$-algebra primitive 
if it has a faithful cyclic representation? 
\end{problem}

We give two results on the $C^*$-algebra $\cO(E)$ 
which suggest that $\cO(E)$ is not an 
``exotic'' $C^*$-algebra. 

\begin{proposition}
The $C^*$-algebra $\cO(E)$ is an inductive limit of 
finite dimensional $C^*$-algebras. 
\end{proposition}

\begin{proof}
For $v\in E^0$ we define $A_{v}\subset \cO(E)$ by 
$$A_{v}=\spa\{s_{(x;w)}s_{(y;w)}^*\mid 
\text{$w\subset v$ and $x,y\in w^{(*)}$}\}.$$
It is not difficult to see that $A_{v}$ is a finite dimensional $C^*$-algebra 
(see Remark~\ref{computation}). 
It is also easy to see that $v_1\subset v_2$ implies 
$A_{v_1}\subset A_{v_2}$ and 
that $\bigcup_{v\in E^0}A_v$ is dense in $\cO(E)$. 
Thus $\cO(E)$ is an inductive limit of 
finite dimensional $C^*$-algebras $A_{v}$. 
\end{proof}

\begin{remark}\label{computation}
For $m\in\N$, 
we define $a(m)\in\N$ by $a(m)=\sum_{k=0}^m m!/k!$. 
The sequence $\{a(m)\}_{m\in\N}$ is determined 
by $a(0)=1$ and $a(m)=ma(m-1)+1$. 
For $v\in E^0$, we have $|v^{(*)}|=a(|v|)$. 
We have 
$$A_{v_0}\cong\bigoplus_{v\subset v_0}\M_{a(|v|)}
=\bigoplus_{m=0}^n\big(
\underbrace{\M_{a(m)}\oplus\cdots\oplus \M_{a(m)}
}_{\frac{n!}{m!(n-m)!}\text{times}}
\big)$$
for $v_0\in E^0$ where $n=|v_0|$. 
For each $v\subset v_0$, 
the matrix units of $\M_{a(|v|)}$ are given by 
$$\Big\{s_{(y;v)}
\Big(s_{(\emptyset,v)}-
\sum_{x\in v_0\setminus v}s_{(x;v\cup\{x\})}s_{(x;v\cup\{x\})}^*
\Big)s_{(z;v)}^*\Big\}_{y,z\in v^{(*)}}.$$
\end{remark}

\begin{remark}
For each $v_0\in E^0$, 
we can define a Cuntz-Krieger $E$-family $T_{v_0}=(T_{v_0}^0,T_{v_0}^1)$ 
on $A_{v_0}$ by 
$$T_{v_0}^0(\delta_v)=\begin{cases}
s_{(\emptyset,v)} & \text{if }v\subset v_0,\\
0 & \text{otherwise},
\end{cases}\qquad 
T_{v_0}^1(\delta_{(x;v)})=\begin{cases}
s_{(x;v)} & \text{if }v\subset v_0,\\
0 & \text{otherwise}.
\end{cases}$$
This gives us a $*$-homomorphism $\psi_{v_0}\colon \cO(E)\to A_{v_0}$ 
such that $\psi_{v_0}(x)=x$ for $x\in A_{v_0}$. 
This proves the next proposition. 
\end{remark}

\begin{proposition}
The $C^*$-algebra $\cO(E)$ is residually finite dimensional. 
\end{proposition}
\end{example}

\end{document}